\documentclass{article}

\usepackage{arxiv}

\usepackage{tikz}
\graphicspath{{images/}}
\usepackage{amsmath}
\usepackage{latexsym, amssymb}
\usepackage{graphicx}
\usepackage{amsmath, amsbsy}
\usepackage{amsopn, amstext}
\usepackage{ifpdf,hyperref}
\usepackage{cancel, color}
\usepackage{epstopdf}

\usepackage[utf8]{inputenc} % allow utf-8 input
\usepackage[T1]{fontenc}    % use 8-bit T1 fonts
\usepackage{hyperref}       % hyperlinks
\usepackage{url}            % simple URL typesetting
\usepackage{booktabs}       % professional-quality tables
\usepackage{amsfonts}       % blackboard math symbols
\usepackage{nicefrac}       % compact symbols for 1/2, etc.
\usepackage{microtype}      % microtypography
\usepackage{lipsum}
\usepackage{graphicx}
\graphicspath{ {./images/} }

\def\ttimes{\,\rotatebox[]{-90}{$\ltimes$}\,}

\def\lvtimes{\vec{\ltimes}}

\def\J{{\bf 1}}

\def\GL{GL}
\def\gl{gl}

\DeclareMathOperator{\Span}{Span}
\DeclareMathOperator{\Col}{Col}
\DeclareMathOperator{\Row}{Row}

\DeclareMathOperator{\lcm}{lcm}

\def\cal{\mathcal}

\def\argmin{argmin}
\def\ra{\rightarrow}

\def\lra{\leftrightarrow}

\def\a{\alpha}
\def\b{\beta}
\def\d{\delta}
\def\D{\Delta}
\def\Tr{Tr}
\def\E{{\cal E}}
\def\Im{Im}

\def\argmin{argmin}

\newcommand{\R}{{\mathbb R}}
\newcommand{\Q}{{\mathbb Q}}
\newcommand{\Z}{{\mathbb Z}}
\newcommand{\F}{{\mathbb F}}

\def\dsum{\mathop{\sum}\limits}

\newtheorem{thm}{Theorem}[section]
\newtheorem{dfn}[thm]{Definition}
\newtheorem{prp}[thm]{Proposition}
\newtheorem{exa}[thm]{Example}
\newtheorem{lem}[thm]{Lemma}
\newtheorem{cor}[thm]{Corollary}
\newtheorem{rem}[thm]{Remark}

\title{Cross-Dimensional Mathematics:\\
A Foundation for Semi-Tensor Product/Semi-Tensor Addition}

\author{
Daizhan Cheng \\
  Key Laboratory of Systems and Control,\\Academy of Mathematics and Systems Sciences,\\
  Chinese Academy of Sciences\\
  Beijing 100190, China\\
  \texttt{dcheng@iss.ac.cn} \\
}

\begin{document}
\maketitle
\begin{abstract}
A new mathematical structure, called the cross-dimensional mathematics (CDM), is proposed.
The CDM considered in this paper consists of three parts: hyper algebra, hyper geometry, and hyper Lie group/Lie algebra.
Hyper algebra proposes some new algebraic structures such as hyper group, hyper ring, and hyper module over matrices and vectors with mixed dimensions (MVMDs). They have sets of classical groups, rings, and modules as their components and cross-dimensional connections among their components. Their basic properties are investigated. Hyper geometry starts from mixed dimensional Euclidian space, and hyper vector space. Then the hyper topological vector space, hyper inner product space, and hyper manifold are constructed. They have a joined cross-dimensional geometric structure.
Finally, hyper metric space, topological hyper group and hyper Lie algebra are built gradually, and finally, the corresponding hyper Lie group is introduced. All these concepts are built over MVMDs, and a couple of most general STP and STA are introduced. Some existing structures/results about STPs/STAs have also been resumed and integrated into this CDM.
\end{abstract}

% keywords can be removed
\keywords{Semi-tensor product/addition (STP/STA) \and hyper manifold \and hyper linear algebra}

\section{Introduction}

The semi-tensor product (STP) of matrices was initiated in 2021\cite{che01}. It is a generalization of classical matrix product by removing the dimension restriction of the classical matrix product. Later on, the semi-tensor addition (STA) was also proposed to weaken the dimension restriction of matrix addition \cite{che19b}.  The past two decades have witnessed the development of STP/STA\cite{for16,muh16}. They have found many applications, including Boolean networks \cite{lu17}, finite games \cite{che21}, dimension-varying systems \cite{che19}, engineering problems \cite{li18}, finite automata \cite{yan22}, coding \cite{zho18}, etc.
In addition to thousands of papers, there are also many monographs  on  STP/STA \cite{che07,che11,che12,che16,che19,che20,che22,che22b,che23b,che24,li18b,liu23,mei10,zha20}, and books with STP chapter or appendix \cite{aku18,vla23}.

Roughly speaking, the STPs and STAs are cross-dimensional operators, which means they are operators acting on objects with different dimensions, i.e., MVMDs,
which can formally be described as follows:
\begin{itemize}
	\item[(i)] Vectors with mixed dimensions: $\R^{\infty}:=\bigcup_{n=1}^{\infty}\R^n$,
	where $\R$ may be replaced by $\F$, which is a field with characteristic number $0$. To make the discussions easy and simplify notations, throughout this paper $\F$ is assumed to be $\R$ unless elsewhere stated. It is easy to see that in the following discussions replacing $\R$ by $\mathbb{C}$ is mostly straightforward.
	
	\item[(ii)] Matrices with mixed dimensions:
	${\cal M}=\bigcup_{m=1}^{\infty}\bigcup_{n=1}^{\infty}{\cal M}_{m\times n}$.
\end{itemize}

In the long history of development of matrix theory, the algebraic structure of matrices, including group, ring, module, general linear algebra, etc., and their geometric structure, including vector space, inner product, metric, general linear group, etc., played  fundamental roles in both theoretical investigations and practical applications. Since initiated two decades ago,  the research on STP/STA has been mainly concentrated on its applications. When investigating the set of MVMDs, which are posed by cross-dimensional operators: STPs and STAs, the existing mathematical tools seem not very efficient.
Hence, it is now the time to explore the  mathematical foundation of STP/STA.

When we tried to pose the existing algebraic and geometric structures on STPs/STAs with MVMDs, we found the existing mathematical framework is not suitable for them. The most critical contradiction lies on the fact that the existing algebraic or geometric structure, such as group or manifold, is for a set of objects with fixed dimensions or fixed sizes, while the STPs/STAs are cross-dimensional operators for MVMDs.

Observing carefully these objects with mixed dimensions, we found that the main difficult point for posing an algebraic or geometric structure on them lies on the trouble of ``multiply unitary".  For example, if we consider the group structure, then in $\R^{\infty}$ there are infinite zeros, precisely speaking, there are $0_n\in \R^n$,  $n=1,2,\cdots$. Similarly, ${\cal M}$ also have infinite $0_{m\times n}$, $m,n=1,2,\cdots$. This fact makes it almost impossible to pose a group structure to these objects (with respect to addition).  If we consider the set of square matrices, denoted by ${\cal M}_1$, then it has also infinite identities, $I_n$, $n=1,2,\cdots$, which makes the group structure of matrices (with respect to product) also very difficult to define.

When we consider the geometric structure, say, vector space structure, we are also facing this ``multiple zero" puzzle.

An inspiration for solving this problem comes from mechanics, where a pseudo-metric space is used.

\begin{dfn}\label{d1.1} \cite{abr78} Consider ${\cal X}=(X,d)\}$, where $X$ is a set,  $d:X\times X\ra \R$ is a mapping.
	\begin{enumerate}
		\item[(1)] ${\cal X}$ is called a metric space if $d$ satisfies
		\begin{itemize}
			\item[(i)] $d(x,y)\geq 0$, $x,y\in X$. Moreover, $d(x,y)=0$, if and only if, $x=y$;
			\item[(ii)] Symmetry:
				$d(x,y)=d(y,x),\quad x,y\in X$;
			\item[(iii)] Triangular Inequality: $d(x,y)+d(y,z)\geq d(x,z),\quad x,y,z\in X$.
		\end{itemize}
		\item[(2)]  ${\cal X}$ is a pseudo-metric space, if $d$ satisfies (ii), (iii), and
		\begin{itemize}
			\item[(i')] $d(x,y)\geq 0$, $x,y\in X$. Moreover, there exists a zero set ${\bf 0}$ such that $d(x,x)=0,\quad x\in {\bf 0}$.
		\end{itemize}
	\end{enumerate}
\end{dfn}

It seems that the concept of ``pseudo-metric space"  can provide a way to address the distance of points with different dimensions. This stimulates us to define a new ``group" which contains multiple identities. Instead of ``pseudo-group" such a multi-identity group  is called a ``hyper group" later in this paper. The reason is that a hyper group contains a set of many (might be countable) groups as its components.

Motivated by this fact, we may modify the existing algebraic and geometric structures to suit the MVMDs with respect to (briefly, with respect to)  STPs/STAs.  The purpose of this paper is to propose a new mathematical framework, called the CDM, for objects with mixed dimensions. We hope that it can provide a solid mathematical foundation for STPs/STAs. Meanwhile, it may result in the development of  a new cross-dimensional algebra and geometry.

The  CDM proposed in this paper consists of three parts: (1) hyper algebra; (2) hyper geometry; (3) hyper Lie algebra and hyper Lie group. The hyper algebra contains hyper group, hyper ring, and hyper module; the hyper geometry includes hyper vector space, hyper inner product space, hyper manifold; and the hyper Lie algebra and hyper Lie group start from mixed dimensional Euclidian space, through constructing topological hyper group,  and finally lead to hyper linear group and hyper linear algebra.

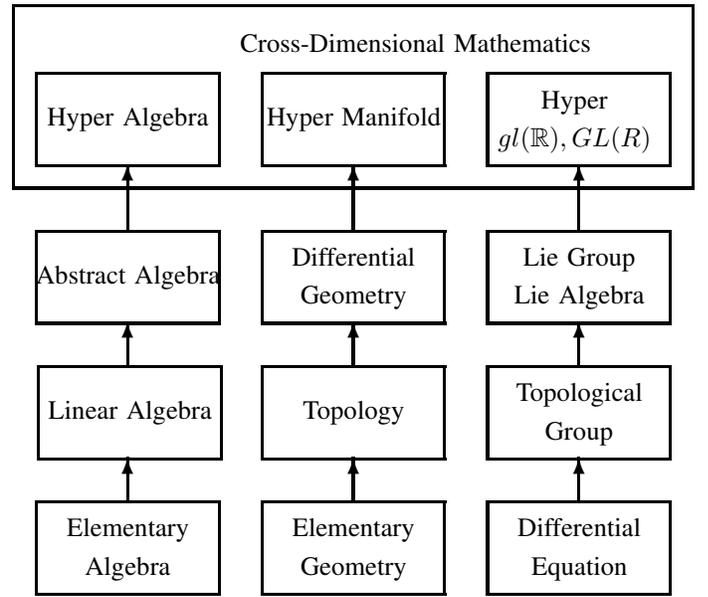
\begin{figure}
	\centering
	\setlength{\unitlength}{8mm}
	\begin{picture}(15,15)\thicklines
    \put(7.3,13){$\oplus$}
    \put(2.5,11.5){\line(0,1){1.5}}
    \put(7.5,11.5){\vector(0,1){1.4}}
    \put(2.5,13){\vector(1,0){4.8}}
     \put(7.7,13){\line(1,0){4.8}}
      \put(12.5,13){\vector(0,-1){1.5}}
		\put(0.5,0){\framebox(4,2){$\begin{array}{c}
					\mbox{Elementary}\\\mbox{Algebra}\end{array}$}}
		\put(0.55,3){\framebox(4,2){Linear Algebra}}
		\put(0.5,6){\framebox(4,2){Abstract Algebra}}
		\put(0.5,9.5){\framebox(4,2){Hyper Algebra}}
		\put(2.5,2){\vector(0,1){1}}
		\put(2.5,5){\vector(0,1){1}}
		\put(2.5,8){\vector(0,1){1.5}}
		\put(0,9){\framebox(15,6){ }}
		\put(5,14){Cross-Dimensional Mathematics}
		\put(5.5,0){\framebox(4,2){$\begin{array}{c}
					\mbox{Elementary}\\\mbox{Geometry}\end{array}$}}
		\put(5.5,3){\framebox(4,2){Topology}}
		\put(5.5,6){\framebox(4,2){$\begin{array}{c}
					\mbox{Differential}\\\mbox{Geometry}\end{array}$}}
		\put(5.5,9.5){\framebox(4,2){Hyper Manifold}}
		\put(7.5,2){\vector(0,1){1}}
		\put(7.5,5){\vector(0,1){1}}
		\put(7.5,8){\vector(0,1){1.5}}
		\put(10.5,0){\framebox(4,2){$\begin{array}{c}
					\mbox{Differential}\\\mbox{Equation}\end{array}$}}
		\put(10.5,3){\framebox(4,2){$\begin{array}{c}
					\mbox{Topological}\\\mbox{Group}\end{array}$}}
		\put(10.5,6){\framebox(4,2){$\begin{array}{c}
					\mbox{Lie Group}\\\mbox{Lie Algebra}\end{array}$}}
		\put(10.5,9.5){\framebox(4,2){$\begin{array}{c}
					\mbox{Hyper}\\\gl(\R), \GL(R) \end{array}$ }}
		\put(12.5,2){\vector(0,1){1}}
		\put(12.5,5){\vector(0,1){1}}
		\put(12.5,8){\vector(0,1){1.5}}
	\end{picture}
	
	\caption{Cross-Dimensional Mathematics vs. Existing Math.\label{Fig1.1}}
\end{figure}

The relationship of CDM with existing mathematical structures is depicted in Figure \ref{Fig1.1}. Figure \ref{Fig1.1} demonstrates also the relationship among the components of CDM.

Some key issues addressed in this paper are emphasized as follows:
\begin{itemize}
	\item[(i)] The most general matrix addition ($\hat{+}$) and matrix product ($\hat{\ttimes}$) are proposed, which make the set of all matrices a hyper ring ($({\cal M},\hat{+},\hat{\ttimes})$). Particularly, the hyper group structure of cross-dimensional permutations ${\bf S}=\bigcup_{n=1}^{\infty}{\bf S}_n$ is obtained.
	
	\item[(ii)] The most general Euclidian hyper vector space ($\R^{\infty}$) is constructed. The vector space structure and the differentiable structure are proposed, which provide a coordinate frame and differentiable structure for hyper manifold. (This function is similar to $\R^n$ for $n$-dimensional manifold.)
	
	\item[(iii)] The Lie group ($\GL(m\times n,\R)$) and Lie algebra ($\gl(m\times n,\R)$) of the set of $m\times n$ dimensional matrices is constructed. Their relationship with general linear group ($\GL(n,\R)$) and general linear algebra  ($\gl(n,\R)$)  is revealed. They provide a foundation for overall hyper matrix Lie group ($\GL(\R)$ and hyper Lie algebra ($\gl(\R)$).
\end{itemize}

The rest part of this paper is  outlined as follows:

Section 2 provides some preliminaries including (i) algebraic structures of matrices/vectors such as matrix group, ring, module, and matrix Lie algebra; (ii) some commonly used STPs and STAs; (iii) dimension-keeping (DK-) STP and pseudo-STP, where pseudo-STP is newly defined. Sections 3-7 are concentrated on the hyper algebraic structure of MVMDs. Starting from lattice, Section 3 introduces hyper group, particularly, matrix/vector hyper groups, and their sub hyper groups. The group decomposition of hyper group into component groups and equivalence group is also presented. Section 4 considers the hyper group structure of permutations, including left and right permutation hyper groups. In Section 5 the hyper group homomorphisms are investigated. Three isomorphism theorems of groups have been extended to hyper groups. Symmetric and skew symmetric and squaring endomorphisms of matrix groups are defined and discussed, which extend the (skew) symmetry and eigenvalue/eigenvector of square matrices to non-square matrices. In Section 6 the hyper ring of matrices is proposed. Its sub hyper rings, homomorphisms, and some properties are investigated. Section 7 considers the hyper module of matrices. The general linear (control) systems are reconsidered under the matrix hyper module framework. Sections 8-11 consider the hyper geometric structure of MVMDs. In Section 8 the hyper vector space is introduced. The inner product of vectors with different dimensions is introduced. Then the norm, the distance are also obtained. In Section 9 a metric space structure is posed on dimension-free Euclidian space $\R^{\infty}$, which turns this cross dimensional Euclidian space a path-wise connected topological space. Furthermore, in Section 10 a differentiable structure is also proposed to $\R^{\infty}$. The smooth functions, (co-) vector fields, distributions, and tensor fields are all introduced for $\R^{\infty}$. Section 11 considers the hyper manifold, where the topological and differentiable structure of $\R^{\infty}$ is transferred to a topological manifold through homeomorphism, similarly to $n$-dimensional manifold with respect to $\R^n$.
Sections 12-15 consider the hyper Lie group and hyper Lie algebra, which are emphasized on their expression on MVMDs as hyper  linear group and hyper linear algebra. Section 12 constructs the general linear algebra ($\gl(m\times n, \R)$) and general linear group ($\GL(m\times n, \R)$) for non-square matrices. Their homomorphisms to classical general linear group ($\GL(t, \R)$) and general linear algebra ($\gl(t,\R)$) are also revealed. Section 13 considers how to construct a hyper Lie group from a topological hyper group. Section 14 investigates the structure of hyper general linear algebras based on hyper vector space and hyper rings. In Section 15 the hyper linear group is proposed. Its relationship with hyper linear algebra is revealed. Finally, Section 16 is a brief summary with some remarks.

Before ending this section we give a list of notations.

\begin{itemize}
	
	\item $\R$: The field of real numbers.
	
	\item $\mathbb{C}$: The field of complex numbers.
	
	\item $\Q$: The field of rational numbers.
	
	\item $\mathbb{F}$: A field of characteristic number $0$,  (in most cases $\mathbb{F}\in\{\mathbb{R}.\mathbb{C}\}$) .
	
	\item $\Z$ ($\Z^+$): Set of (positive) integers.
	
	\item $t|r$: $t,r\in \Z^+$ and $r/t\in \Z^+$.

	\item $\mathbb{F}^n$: $n$ dimensional Euclidian space over field $\F$.
	
	\item ${\cal M}^\mathbb{F}_{m\times n}$: set of $m\times n$ matrices over $\F$ ( the superscript $\F$ is mostly omitted).
	
	\item ${\cal M}_{\mu}:=\{A_{m\times n}\;|\; m/n=\mu\} $ where $0<\mu\in \Q$.
	
	\item $[a,b]$: the set of integers $\{a,a+1,\cdots,b\}$, where $a\leq b$.
	
	\item $\lcm(n,p)$ (or $n\vee p$) : the least common multiple of $n$ and $p$.
	
	\item ${\cal D}=\{1,2,\cdots,n\}$.
	
	\item $\gcd(n,p)$ (or $n\wedge p$) : the greatest common division of $n$ and $p$.
	
	\item $\d_n^i$: the $i$-th column of the identity matrix $I_n$.
	
	\item $\D_n:=\{\d_n^i\mid i=1,\cdots,n\}$.
	
	\item $\J_k:=(\underbrace{1,\cdots,1}_k)^{\rm T}$ ($\J_{m\times n}=\J_m^{\rm T}\J_n$).
	
	\item $J_n:=\frac{1}{n}\J_{n\times n}.$

	\item $\E_n:=\frac{1}{\sqrt{n}}\J_n$.
	
	\item $\E_{m\times n}:=E_m E^{\rm T}_n$.
	
	\item ${\cal L}_{m\times n}$: the set of logical matrices, i.e.
	${\cal L}_{m\times n}=\{ A\in {\cal M}_{m\times n}\;|\; \Col(A)\subset \D_m\}$.
	
\item $\d_n[i_1,\cdots,i_m]:=[\d_n^{i_1},\cdots,\d_n^{i_m}]$.

	\item $\Col(A)$: the set of columns of matrix $A$, and $\Col_i(A)$ is the $i$-th column of $A$.
	
	\item $\Row(A)$: the set of rows of matrix $A$, and $\Row_i(A)$ is the $i$-th row of $A$.
	
	\item $\Span(\cdot)$: subspace or dual subspace generated by $\cdot$.
	
	\item ${\bf S}_n$ (${\bf S}$): order $n$ permutation group (hyper permutation group).
	
	\item ${\bf A}_n$ (${\bf A}$): order $n$ alternative group (hyper alternative group).
	
	\item $\perp$: perpendicular relation between $x\in \R^n$ and $h\in \R^n_*$.
	
	\item $\bar{+}$: STA with respect to $I_{n\times n}$.
	
	\item $\hat{+}$: STA with respect to $\J_n$ or $\E_n$ ($\E_{m\times n}$).
	
	\item $\ltimes$: (left) matrix-matrix semi-tensor product, which is the default one.
	
	\item $\rtimes$: right semi-tensor product of matrices.
	
	\item $\circ$: second (left) matrix-matrix semi-tensor product.
	
	\item $\lvtimes$: matrix-vector semi-tensor product.
	
	\item $\vec{\circ}$: second matrix-vector semi-tensor product.
	
	\item $\ttimes$: dimension-keeping semi-tensor product of matrices.
	
	\item $\hat{\ttimes}$: pseudo-STP.
	
	\item $\sup(a,b)$  (or $a\vee b$): least common upper bound of $a$ and $b$.
	
	\item $\inf(a,b)$  (or $a\wedge b$): greatest  common lower bound of $a$ and $b$.
	
	\item $\Box(A)$: squaring matrix of $A$.
	
	\item $\GL(n,\R)$ ($\GL(\R)$): (hyper) general linear group.
	
	\item $\gl(n,\R)$ ($\gl(\R)$): (hyper) general linear algebra.
	
	\item $H<G$ ($H\lhd G$): $H$ is a (normal) subgroup of $G$.
	
	\item $H<_h G$ ($H\lhd_h G$): $H$ is a (normal) hyper-subgroup of $G$.
	
	\item ${\rm Id}$: identity mapping.
	
	\item ${\rm Im}$: image set.
	
	\item ${\rm Ker}$: kernel set.
\end{itemize}

\section{Preliminaries}

\subsection{Algebraic Structure on Matrices/Vectors with proper dimensions}

We first review some algebraic structures over matrices/vectors.

\begin{dfn}\label{d2.1.1} Consider a couple ${\cal G}=(G,*)$, where $G$ is a set, and $*:G\times G\ra G$ is a binary operator.
	\begin{itemize}
		\item[(i)] ${\cal G}$ is a semi-group, if $(a*b)*c=a*(b*c),\quad a,b,c\in G$.
		\item[(ii)] A semi-group ${\cal G}$ is a monoid (i.e., a semi-group with an identity), if
		there exists an identity $e\in G$, such that $a*e=e*a=a,\quad \forall a\in G$.
		\item[(iii)] A monoid  ${\cal G}$ is a group, if
		for each $a\in G$, there exists an $a^{-1}$, such that $a*a^{-1}=a^{-1}*a=e,\quad a\in {\cal G}$
	\end{itemize}
	In addition, if $a*b=b*a,\quad \forall a,b\in G$, then it is called an abelian group.
\end{dfn}

\begin{exa}\label{e2.1.2}
	\begin{itemize}
		\item[(i)] $(\R^n,+)$  is an abelian group with
		the identity $0_n\in \R^n$.
		
		\item[(ii)] $({\cal M}_{m\times n}, +)$ is an abelian group with
		the identity $0_{m\times n}\in {\cal M}_{m\times n}$.
		
		\item[(iii)] $({\cal M}_{n\times n},+)$
		is an abelian group with the identity $0_n$.
		
		\begin{itemize}
			\item Set
			${\cal S}_{n}:=\{A\in {\cal M}_{n\times n}\;|\; A^{^{\rm T}}=A \}$,
			$({\cal S}_{n},+)< ({\cal M}_{n\times n},+)$ is the symmetric subgroup.
			\item Set
			${\cal K}_{n}:=\{A\in {\cal M}_{n\times n}\;|\; A^{^{\rm T}}=-A \}$,
			$({\cal K}_{n},+)< ({\cal M}_{n\times n},+)$ is the skew-symmetric subgroup.
			
			\item Set
			${\cal Q}_{n}:=\{A\in {\cal M}_{n\times n}\;|\; \Tr(A)=0\}$,
			$({\cal Q}_{n},+)< ({\cal M}_{n\times n},+)$ is the zero-trace subgroup.
			
			(There are also some other subgroups such as diagonal, upper-triangular, lower-triangular, etc.)
		\end{itemize}
		
		\item[(iv)] $({\cal M}_{n\times n},\times)$ is a monoid with the identity $I_n$.
		
		\item[(v)] ${\cal T}_n:=\{ A\in {\cal M}_{n\times n} \;|\; A~\mbox{is invertible.}\}$. Then
		$({\cal T}_n, \times)$ is a group with identity $I_n$.
		
		\begin{itemize}
			\item Set
			${\cal O}_{n}:=\{A\in {\cal T}_{n}\;|\; A^{^{\rm T}}=A^{-1} \}$,
			$({\cal O}_{n},\times)< ({\cal T}_{n},\times)$ is the orthogonal subgroup.
			\item Set
			${\cal H}_{n}:=\{A\in {\cal T}_n \; \det(A)=1\}$,
			$({\cal H}_{n},\times)< ({\cal T}_{n},\times)$ is the special linear subgroup.
			\item Set
			${\cal P}_{n}:={\cal T}_{n}\bigcap {\cal L}_{n\times n}$,
			where ${\cal L}_{n\times n}$ is the set of logical matrices ($A$ is a logical matrix if $\Col(A)\subset \D_n$), then $({\cal P}_{n},\times)< ({\cal T}_{n},\times)$ is the permutation  sub-group.
			
			(there are also some other sub-groups such as non-singular diagonal, non-singular upper-triangular, non-singular lower-triangular, etc.)
		\end{itemize}
	\end{itemize}
\end{exa}

\begin{dfn}\label{d2.1.3} Triple ${\cal R}=(R,+,*)$ is a ring, if
	\begin{itemize}
		\item[(i)] $(R,+)$ is an abelian group;
		\item[(ii)] $(R,*)$ is a semi-group;
		\item[(iii)] (Distributivity) $(a+b)*c=a*c+b*c$, $a*(b+c)=a*b+a*c$, $a,b,c\in R$.
	\end{itemize}
	In addition, if $(R,*)$ is a monoid, ${\cal R}$ is a unitary ring; if $a*b=b*a,\quad \forall a,b\in R$,   ${\cal R}$ is a commutative
	ring.
\end{dfn}

\begin{dfn}\label{d2.1.301} Assume ${\cal R}=(R,+,*)$ is a ring and $H\subset R$.  If
	$(H,+,*)$ is also a ring,  $H$ is called a sub-ring of $R$.
\end{dfn}

\begin{exa}\label{e2.1.4}
	$({\cal M}_{n\times n},+,\times)$ is a unitary ring. It has some sub-rings as follows.
	
	\begin{itemize}
		\item[(i)] ${\cal D}_n:=\{A\in {\cal M}_{n\times n}\;|\; A~\mbox{is diagonal.}\}$
		Then $({\cal D}_n, +,\times)$ is a sub-ring of $({\cal M}_{n\times n},+,\times)$ .
		\item[(ii)] ${\cal U}_n:=\{A\in {\cal M}_{n\times n}\;|\; A~\mbox{is upper triangular.}\}$.
		Then $({\cal U}_n, +,\times)$ is a sub-ring of $({\cal M}_{n\times n},+,\times)$ .
		\item[(iii)] ${\cal B}_n:=\{A\in {\cal M}_{n\times n}\;|\; A~\mbox{is lower triangular.}\}$.
		Then $({\cal B}_n, +,\times)$ is a sub-ring of $({\cal M}_{n\times n},+,\times)$ .
	\end{itemize}
\end{exa}

\begin{dfn}\label{d2.1.5}  Let ${\cal R}=(R,+,*)$ be a ring, and ${\cal G}=(G,+)$ be an abelian group. ${\cal R}$, ${\cal G}$ with a mapping $\pi:{\cal R}\times {\cal G} \ra {\cal G}$ is a (left) ${\cal R}-$ module, if
	\begin{itemize}
		\item[(i)] $r(a+b)=ra+rb,\quad r\in {\cal R},~a,b\in {\cal G}$;
		\item[(ii)]
		$(r+s)a=ra+sa,\quad r,s\in {\cal R},~a,b\in {\cal G}$;
		\item[(iii)]
		$r(sa)=(rs)a$.
	\end{itemize}
	In addition, if there exists $e\in {\cal R}$ such that for all  $a\in G$
	$ea=a,\quad \forall a\in G$,
	then ${\cal G}$ is called a unitary ${\cal R}-$ module.
\end{dfn}

\begin{exa}\label{e2.1.6}
	\begin{itemize}
		\item[(i)]
		$({\cal M}_{n\times n},+,\times)$ is a ring with identity.
		\item[(ii)]
		$(\R^n,+)$ is an abelian group.
		\item[(iii)] Let $\pi: {\cal M}_{n\times n}\times \R^n\ra \R^n$ be the classical matrix-vector product.
	\end{itemize}
	Then
	$(({\cal M}_{n\times n},+,\times), (\R^n,+), \pi)$
	is a module.
\end{exa}

\begin{dfn}\label{d2.1.7}\cite{boo86} Let $V$ be a vector space over $\R$ and there is a mapping, called the Lie bracket and denoted by $[\cdot,\cdot]: V\times V\ra V$. $(V, [\cdot,\cdot])$ is called a Lie algebra, if the followings are satisfied.
	\begin{itemize}
		\item[(i)] Skew-Symmetry: $[x,y]=-[y,x],\quad x,y\in V$.
		\item[(ii)] Linearity: $[ax+by,z]=a[x,z]+b[y,z]$, $ a,b\in \R;~x,y.z\in V$.
		\item[(iii)] Jacobi Identity: $[[x,y],z]+[[y,z],x]+[[z,x],y]=0$.
	\end{itemize}
\end{dfn}

\begin{exa}\label{e2.1.8} Consider ${\cal M}_{n\times n}$. Define $[A,B]=AB-BA, \quad A,B\in {\cal M}_{n\times n}$.
	Then $({\cal M}_{n\times n},[\cdot,\cdot])$ becomes a Lie algebra, called the general linear algebra and denoted by $\gl(n,\R)$.
\end{exa}

Note that all the above algebraic structures over matrices/vectors are restricted to a particular set of matrices with certain fixed dimensions. We refer to \cite{hun74} or any other text book of abstract algebra for the related concepts. These concepts will be extended to MVMDs in this paper.

\subsection{STP and STA of Matrices/Vectors}

The main objects of linear algebra are vectors and matrices. In linear algebra both the addition and the product of  two matrices, or matrix with vector, or two vectors, require the factors to meet certain dimension matching condition, here the product of two vectors means the inner product.

STP/STA are proposed to deal with MVMDs. The sets of vectors with mixed dimensions and matrices with mixed dimensions are described by  $\R^{\infty}$ and ${\cal M}$ respectively. In addition we also need ${\cal M}_{\mu}=\{A\in {\cal M}_{m\times n}\;|\; m/n=\mu\},\quad \mu\in \Q_+$.

\begin{rem}\label{r2.2.0}
	Consider ${\cal M}_{\mu}$, where $\mu=\frac{\mu_y}{\mu_x}\in \Q_+$. Moreover, to be uniqueness, we assume
	$\gcd(\mu_y,\mu_x)=1$. Then
	\begin{itemize}
		\item[(i)]
		${\cal M}_1$ is the set of square matrices.
		\item[(ii)]
		${\cal M}=\bigcup_{\mu\in \Q_+}{\cal M}_{\mu}$
		is a partition.
	\end{itemize}
\end{rem}

The STP and STA have been proposed to extend the matrix (vector) product and addition to overall (or larger then the original) sets.
The followings are some commonly used ones.

\begin{dfn}\label{d2.2.1} \cite{che11,che19}  Let $A,~B\in {\cal M}$, specified by $A\in {\cal M}_{m\times n}$, $B\in {\cal M}_{p\times q}$, and $t=\lcm(n,p)$.
	\begin{itemize}
		\item[(i)] The first type of matrix-matrix (MM-) STP of $A,B$, denoted by $A\ltimes B$, is defined by
		\begin{align}\label{2.2.1}
			A\ltimes B:=(A\otimes I_{t/n}) (B\otimes I_{t/p})\in {\cal M}_{mt/n\times qt/p}.
		\end{align}
		\item[(ii)] The second type of MM-STP of $A,B$, denoted by $A\circ B$, is defined by
		\begin{align}\label{2.2.2}
			A\circ B:=(A\otimes J_{t/n}) (B\otimes J_{t/p}\in {\cal M}_{mt/n\times qt/p}.
		\end{align}
	\end{itemize}
\end{dfn}

\begin{dfn}\label{d2.2.2} \cite{che12}  Let $A\in {\cal M}$ and $x\in \R^{\infty}$, specified by $A\in {\cal M}_{m\times n}$, $x\in \R^{p}$, and $t=\lcm(n,p)$.
	\begin{itemize}
		\item[(i)] The first type of matrix-vector (MV-) STP of $A$ and $x$, denoted by $A\lvtimes x$, is defined by
		\begin{align}\label{2.2.3}
			A\lvtimes x:=(A\otimes I_{t/n}) (x\otimes \J_{t/p})\in \R^{mt/n}.
		\end{align}
		\item[(ii)] The second type of MV-STP of $A$ and $x$, denoted by $A\vec{\circ} x$, is defined by
		\begin{align}\label{2.2.4}
			A\vec{\circ} x:=(A\otimes J_{t/n\times t/n}) (x\otimes \J_{t/p})\in \R^{mt/n}.
		\end{align}
	\end{itemize}
\end{dfn}

\begin{dfn}\label{d2.2.3} \cite{che12}  Let $x,~y\in \R^{\infty}$, specified by $x\in \R^{m}$, $y\in \R^n$ ,  and $t=\lcm(m,n)$. The vector-vector (VV-) STP  of $x$ and $y$, denoted by $x~\bar{\cdot}~y$, is defined by
	\begin{align}\label{2.2.5}
		x~\bar{\cdot}~y:=\langle x\otimes \J_{t/m},  y\otimes \J_{t/n} \rangle \in \R,
	\end{align}
	where $\langle \cdot,\cdot\rangle$ is the classical inner product on Euclidian space.
\end{dfn}

\begin{rem}\label{r2.2.4}
	\begin{itemize}
		\item[(i)] All the MM-STPs, MV-STPs, and VV-STP  are the generalization of classical matrix product.
		That is, when the dimension-matching conditions for classical matrix-matrix (or matrix-vector or vector-vector) products are satisfied, they degenerate to classical products.
		\item[(ii)] All the MM-STPs and MV-STPs defined in (\ref{2.2.1}), (\ref{2.2.3}), and (\ref{2.2.5}) are called the left STPs, because the main objects ($A$, $B$, etc.) are lying on left. They have their corresponding right STPs. Say,
		\begin{align}\label{2.2.501}
			(A\rtimes B):=(I_{t/n}\otimes A)(I_{t/p}\otimes B) ,
		\end{align}
		etc.   Such right extensions are also available for all other STPs defined in the sequel. And a parallel discussion shows similar properties of right STPs corresponding to the default left STPs.  (We leave such discussions for the reader and will not mention this again.)
		
		\item[(iii)] All the MM-STPs, MV-STPs maintain most properties of the classical matrix product available. Particularly, the associativity and the distributivity hold. Say, $(A\ltimes B)\ltimes C=A\ltimes (B\ltimes C), \quad A,B,C\in {\cal M}$;
		and $(A+ B)\ltimes C=A\ltimes C+ B\ltimes C$; $C\ltimes (A+ B)=C\ltimes A+ C\ltimes B$, $A,B,C\in {\cal M}$.
	\end{itemize}
\end{rem}

\begin{dfn}\label{d2.2.5} Some STAs of matrices/vectors are defined as follows:
	\begin{itemize}
		\item[(i)] $\R^{\infty}$: Assume $x\in \R^m$, $y\in \R^n$, and $t=\lcm(m,n)$, then
		\begin{align}\label{2.2.6}
			x\bar{\pm}y:=(x\otimes \J_{t/m})\pm(y\otimes \J_{t/n}).
		\end{align}
	
		\item[(ii)]  ${\cal M}_1$: Assume $A\in {\cal M}_{m\times m}$, $B\in {\cal M}_{n\times n}$, and $t=\lcm(m,n)$, then
		\begin{align}\label{2.2.7}
			A\bar{\pm}B:=(A\otimes I_{t/m})\pm (B\otimes I_{t/n}).
		\end{align}
		
		An alternative STA is defined as
\begin{align}\label{2.2.8}
    A\hat{\pm}B:=(A\otimes J_{t/m})\pm (B\otimes J_{t/n}).
\end{align}
		\item[(iii)] ${\cal M}_{\mu}$: $A\in {\cal M}_{m\mu_y\times m\mu_x}$, $B\in {\cal M}_{n\mu_y\times n\mu_x}$, and $t=\lcm(m,n)$:
		Its additivity can be define by either $\bar{\pm}$, as in (\ref{2.2.7}), or  $\hat{\pm}$, as in (\ref{2.2.8}).
	\end{itemize}
\end{dfn}

To define ``addition" for more general case, we need a new set of matrices as ``weights".
Set
\begin{align}\label{2.2.9}
	\begin{array}{l}
		\E_{n}:=\frac{1}{\sqrt{n}}\J_n;\\
		\E_{m\times n}:=\E_{m}\E^{^{\rm T}}_n=\frac{1}{\sqrt{mn}}\J_{m\times n}, \quad m,n\in \Z_+.
	\end{array}
\end{align}

Then the following properties of this set of matrices are easily verifiable, but they are fundamental in further discussion.

\begin{prp}\label{p2.2.6}
	\begin{align}\label{2.2.10}
		\begin{cases}
			\E^{^{\rm T}}_{n}\E_n=1;\\
			\E_{n}\otimes \E_m=\E_m\otimes \E_n=\E_{mn};\\
			\E_{m\times n}\otimes \E_{p\times q}=\E_{p\times q}\otimes \E_{m\times n}=\E_{mp\times nq};\\
			\E^{^{\rm T}}_{m\times n}=\E_{n\times m};\\
			J_n=\E_{n\times n}.
		\end{cases}
	\end{align}
\end{prp}

Using this set of weights, some other STAs can be defined.

\begin{dfn}\label{d2.2.7}
	\begin{itemize}
		\item[(i)] $\R^{\infty}$: Assume $x\in \R^m$, $y\in \R^n$, and $t=\lcm(m,n)$, then
		\begin{align}\label{2.2.11}
			x\hat{\pm}y:=(x\otimes \E_{t/m})\pm (y\otimes \E_{t/n}),
		\end{align}
		
		\item[(ii)] Consider ${\cal M}$, and assume $A\in {\cal M}_{m\times n}$, $B\in {\cal M}_{p\times q}$, $s=\lcm(m,p)$, and $t=\lcm(n,q)$.
		Then we define the STA over ${\cal M}$ as follows:
		\begin{align}\label{2.2.12}
			A\hat{\pm}B:=(A\otimes \E_{s/m\times t/n})\pm (B\otimes \E_{s/p\times t/q}).
		\end{align}
	\end{itemize}
\end{dfn}

Note that all the STPs and STAs, defined up to this time are the generalizations of the classical matrix product/addition. That is, when the dimension matching conditions required by the classical products are satisfied, they degenerate to the classical matrix product/addition.

%%%%%%%%%%%%%%%%%%%%%%%%%%%%%%%%%%%
Before ending this subsection we give a formula for calculating equivalence transformation, which will be used in the sequel.

\begin{prp}\label{p2.2.8}
	Let $A\in {\cal M}_{m\times n}$, and $B\in {\cal M}_{p\times q}$. Then
	\begin{align}\label{2.2.13}
		\begin{array}{l}
			(A\otimes \E_{r\times kp}) (B\otimes \E_{kn\times s})= (A\otimes \E_{r\times p}) (B\otimes \E_{n\times s}),\\
			(A\otimes \E_{kr\times p}) (B\otimes \E_{n\times ks})= (A\otimes \E_{r\times p}) (B\otimes \E_{n\times s}).\\
		\end{array}
	\end{align}
\end{prp}

\noindent{\it Proof.} We prove the first equality only. The proof of the second one is similar.
$$
\begin{array}{l}
	\mbox{LHS of (\ref{2.2.13})}= (A\otimes \E_{r\times p}\otimes \E^{^{\rm T}}_{k}) (B\otimes \E_{n\times s}\otimes \E_k)\\
	~=[(A\otimes \E_{r\times p}) (B\otimes \E_{n\times s})]\otimes (\E^{^{\rm T}}_k\E_k)\\
	~=(A\otimes \E_{r\times p}) (B\otimes \E_{n\times s})=\mbox{RHS of (\ref{2.2.13})}.
\end{array}
$$
\hfill $\Box$
%%%%%%%%%%%%%%%%%%%%%%%%%%%%%%%%%%%%%%%%%%%%%%

\subsection{DK-STP and Pseudo-STP}

The dimension-keeping (DK-) STP is defined as follows.

\begin{dfn}\label{d2.3.1} Let $A\in {\cal M}_{m\times n}$, $B\in {\cal M}_{p\times q}$, and $t=\lcm(n,p)$. The DK-STP of $A$ and $B$, denoted as $A\ttimes B$, is defined by
	\begin{align}\label{2.3.1}
		A\ttimes B:=(A\otimes \E^{^{\rm T}}_{t/n})(B\otimes \E_{t/p}).
	\end{align}
\end{dfn}

\begin{rem}\label{r2.3.2}
	DK-STP was newly proposed in \cite{chepr} as follows:
	\begin{align}\label{2.3.2}
		A\ttimes B:=(A\otimes \J^{^{\rm T}}_{t/n}) (B\otimes \J_{t/p}).
	\end{align}
	The original one defined by (\ref{2.3.2}) has very clear physical meaning.  Duplicating the columns of the first matrix and the rows of the second matrix to least times, such that the enlarged factor matrices meet the dimension matching condition of the classical matrix product.
	
	In fact, the modified Definition \ref{d2.3.1} here uses certain weights for the duplications, which keep the ``norms" of the enlarged matrices unchanged. In the sequel you will see that  this modification is indeed meaningful. In \cite{chepr2} (\ref{2.3.1}) is called a weighted DK-STP.
\end{rem}

An obvious property of DK-STP is: if $A,B\in {\cal M}_{m\times n}$, the $A\ttimes B\in {\cal M}_{m\times n}$.
That is why it is named by ``dimension-keeping".

The following properties about DK-STP are mimic to the one in \cite{chepr}.

\begin{prp}\label{p2.3.3} Define a bridge matrix as
	\begin{align}\label{2.3.3}
		\Psi_{n\times p}:=(I_n\otimes \E^{^{\rm T}}_{t/n})(I_p\otimes \E_{t/p}),\quad n,p\in \Z_+,
	\end{align}
	where $t=\lcm(n,p)$, then
	\begin{align}\label{2.3.4}
		A\ttimes B=A\Psi_{n\times p}B,\quad A\in {\cal M}_{m\times n}, B\in {\cal M}_{p\times q}.
	\end{align}
\end{prp}

\begin{rem}\label{r2.3.301} Corresponding to (\ref{2.3.2}), the  bridge matrix becomes
	\begin{align}\label{2.3.301}
		\Psi'_{n\times p}:=(I_n\otimes \J^{^{\rm T}}_{t/n})(I_p\otimes \J_{t/p}),\quad n,p\in \Z_+,
	\end{align}
	which is only different from (\ref{2.3.3}) by a coefficient, i.e.,
	\begin{align}\label{2.3.302}
		\Psi_{n\times p}:=\frac{\sqrt{np}}{t}\Psi'_{n\times p}.
	\end{align}
\end{rem}

Using Proposition \ref{p2.3.3}, the following results are obvious.

\begin{prp}\label{p2.3.4}
	\begin{itemize}
		\item[(i)] $\ttimes$ is a generalization of classical matrix product, if $n=p$ it turns out that $A\ttimes B = AB$.
		\item[(ii)] (Distributivity)
		\begin{align}\label{2.3.5}
			\begin{array}{l}
				(A+B)\ttimes C=A\ttimes C+B\ttimes C, \\
				C\ttimes (A+B)= C\ttimes A+C\ttimes B,
			\end{array}
		\end{align}
		where $A$ and $B$ are of the same dimension.
		\item[(iii)] (Associativity)
		\begin{align}\label{2.3.6}
			(A\ttimes B)\ttimes C=A(B\ttimes C).
		\end{align}
		\item[(iv)]
		\begin{align}\label{2.3.7}
			(A\ttimes B)^{\rm T}=B^{\rm T}\ttimes A^{\rm T}.
		\end{align}
	\end{itemize}
\end{prp}

DK-STP brings a new algebraic structure to non-square matrices. The following result is interesting and useful.

\begin{prp}\label{p2.3.5} \cite{chepr} $({\cal M}_{m\times n}, +,\ttimes)$ is a ring.
\end{prp}

Using DK-STP, we can define a pseudo-STP as follows:

\begin{dfn}\label{d2.3.6} Consider ${\cal M}$, and assume $A\in {\cal M}_{m\times n}$, $B\in {\cal M}_{p\times q}$, $s=\lcm(m,p)$, and $t=\lcm(n,q)$.
	Then the pseudo-STP of $A$ and $B$ is defined as follows:
	\begin{align}\label{2.3.8}
		A\hat{\ttimes}B:=(A\otimes \E_{s/m\times t/n})\ttimes (B\otimes \E_{s/p\times t/q}).
	\end{align}
\end{dfn}

\begin{rem}\label{r2.3.7}
	\begin{itemize}
		\item[(i)] As a special case, we consider ${\cal M}_{\mu}$, and assume $A\in {\cal M}_{m\mu_y\times m\mu_x}$, $B\in {\cal M}_{n\mu_y\times n\mu_x}$, and $t=\lcm(m,n)$.
		Then the pseudo-STP for ${\cal M}_{\mu}$ becomes:
		\begin{align}\label{2.3.9}
			A\hat{\ttimes}B:=(A\otimes J_{t/m}) \ttimes (B\otimes J_{t/n}).
		\end{align}
		\item[(ii)] Comparing (\ref{2.3.8}) with (\ref{2.2.12}), one sees easily that in certain sense this pseudo-STP is ``consistent" with
		the STA defined by (\ref{2.2.12}). It will bee seen in the sequel that this consistence is very meaningful.
		\item[(iii)] The reason to call   (\ref{2.3.8}) or  (\ref{2.3.9}) a pseudo-STP is: they do not always consistent with classical matrix product even if the dimension matching condition is satisfied.
		So they do not meet the fundamental requirement of STP.
	\end{itemize}
\end{rem}

Next, we explore some basic properties of pseudo-STP.

\begin{prp}\label{p2.3.8} Let $A,~B\in {\cal M}$. Then
	\begin{itemize}
		\item[(i)]
		\begin{align}\label{2.3.10}
			(A\hat{\ttimes}B)^{^{\rm T}}=B^{^{\rm T}}\hat{\ttimes} A^{^{\rm T}}.
		\end{align}
		\item[(ii)]
		\begin{align}\label{2.3.101}
			(A\hat{+}B)^{^{\rm T}}=A^{^{\rm T}}\hat{+} B^{^{\rm T}}.
		\end{align}
	\end{itemize}
\end{prp}

\noindent{\it Proof.}

We prove (\ref{2.3.10}) only. The proof of (\ref{2.3.101}) is similar.

Assume $A\in {\cal M}_{m\times n}$, $B\in {\cal M}_{p\times q}$, $\xi=\lcm(m,p)$, $\eta=\lcm(n,q)$,
and $\d=\lcm(\xi,\eta)$. Then
$$
\begin{array}{l}
	A\hat{\ttimes} B=(A\otimes \E_{\xi/m\times \eta/n})\ttimes
	(B\otimes \E_{\xi/p\times \eta/q})\\
	~=(A\otimes \E_{\xi/m\times \eta/n}\otimes \E^{^{\rm T}}_{\d/\eta})
	(B\otimes \E_{\xi/p\times \eta/q}\otimes \E_{\d/\xi})\\
	~=(A\otimes \E_{\xi/m\times \d/n})
	(B\otimes \E_{\d/p\times \eta/q}).\\
\end{array}
$$
It follows that
$$
%\begin{array}{l}
(A\hat{\ttimes} B)^{^{\rm T}}= (B^{^{\rm T}}\otimes \E_{\eta/q\times \d/p}).\\
(A^{^{\rm T}}\otimes \E_{\d/n\times \xi/m})
$$
On the other hand,
$$
\begin{array}{l}
	B^{^{\rm T}}\hat{\ttimes} A^{^{\rm T}}=(B^{^{\rm T}}\otimes \E_{\eta/q\times \xi/p})\ttimes
	(A^{^{\rm T}}\otimes \E_{\eta/n\times \xi/m})\\
	~=(B^{^{\rm T}}\otimes \E_{\eta/q\times \d/p})
	(A^{^{\rm T}}\otimes \E_{\d/n\times \xi/m}).
\end{array}
$$
The conclusion follows.

\hfill $\Box$

\begin{prp}\label{p2.3.10} $({\cal  M}, \hat{\ttimes})$ is a semi-group. That is, assume $A,~B,~C\in {\cal M}$, then
	\begin{align}\label{2.3.11}
		(A\hat{\ttimes}B)\hat{\ttimes} C=A\hat{\ttimes} (B\hat{\ttimes} C).
	\end{align}
\end{prp}

\noindent{\it Proof.}
Assume
$$
A\in {\cal M}_{m\times n},\;B\in {\cal M}_{p\times q},\; C\in {\cal M}_{r\times s},
$$
and
$$
\begin{array}{ll}
	m\vee p=\xi,~n\vee q=\eta,&\xi\vee \eta=\d,~\xi\vee r=\a,\\
	\eta\vee s=\b,~\a\vee \b=\theta,&p\vee r=a,~q\vee s=b,\\
	a\vee b=e,&m\vee a=m\vee p\vee r=\a,\\
	n\vee b=n\vee q\vee s=\b,&n\vee p=\mu,\\
	q\vee r=\lambda,&~\\
\end{array}
$$
Then
$$
\begin{array}{l}
	A\hat{\ttimes} B=(A\otimes \E_{\xi/m\times \eta/n}  )\ttimes
	(B\otimes \E_{\xi/p\times \eta/q})\\
	~=(A\otimes E_{\xi/m\times \eta/n} \otimes \E^{^{\rm T}}_{\d/\eta} )
	(B\otimes \E_{\xi/p\times \eta/q}\otimes \E_{\d/\xi}  )\\
	~=(A\otimes \E_{\xi/m\times \d/n} )
	(B\otimes \E_{\d/p\times \eta/q} )\in {\cal M}_{\xi\times \eta},\\
\end{array}
$$
it follows that
$$
\begin{array}{l}
	(A\hat{\ttimes} B)\hat{\ttimes} C
	=\left\{\left[(A\otimes E_{\xi/m\times \d/n} )
	(B\otimes \E_{\d/p\times \eta/q} )\right]\right.\\
	~\left.\otimes \E_{\a/\xi\times \b/\eta}\right\}\ttimes(C\otimes \E_{\a/r\times \b/s})\\
	=\left\{\left[(A\otimes \E_{\xi/m\times \d/n} )
	(B\otimes \E_{\d/p\times \eta/q} )\right]\right.\\
	~\left.\otimes (\E_{\a/\xi}\E^{^{\rm T}}_{\b/\eta})\right\}\ttimes(C\otimes \E_{\a/r\times \b/s})\\
	=\left\{\left[(A\otimes \E_{\a/m\times \d/n} )
	(B\otimes \E_{\d/p\times \b/q} )\right]\right\}\\
	~\ttimes(C\otimes \E_{\a/r\times \b/s})\\
	=(\left\{\left[(A\otimes \E_{\a/m\times \d/n} )
	(B\otimes \E_{\d/p\times \b/q} )\right]\right\}\otimes \E^{^{\rm T}}_{\theta/\b})\\
	((C\otimes \E_{\a/r\times \b/s})\otimes \E_{\theta/\a})\\
	=(A\otimes \E_{\a/m\times \d/n} )
	(B\otimes \E_{\d/p\times \theta/q} )
	(C\otimes \E_{\theta/r\times \b/s})\\
	=(A\otimes \E_{\a/m\times \mu/n} )
	(B\otimes \E_{\mu/p\times \lambda/q} )
	(C\otimes \E_{\lambda/r\times \b/s}).\\
\end{array}
$$
On the other hand,
$$
\begin{array}{l}
	B\hat{\ttimes} C=(B\otimes \E_{a/p\times b/q}   )\ttimes
	(C\otimes \E_{a/r\times b/s})\\
	~=(B\otimes \E_{a/p\times e/q})
	(C\otimes \E_{e/r\times b/s})\in {\cal M}_{a\times b},\\
\end{array}
$$
and then
$$
\begin{array}{l}
	A\hat{\ttimes} (B\hat{\ttimes} C)
	=(A\otimes \E_{\a/m\times \b/n} )\ttimes\\
	~\left\{\left[(B\otimes \E_{a/p\times e/q})
	(C\otimes \E_{e/r\times b/s})\right]\otimes \E_{\a/a\times \b/b}\right\}\\
	=(A\otimes \E_{\a/m\times \b/n} )\ttimes
	\left[(B\otimes \E_{\a/p\times e/q})
	(C\otimes \E_{e/r\times \b/s})\right]\\
	=\left[(A\otimes \E_{\a/m\times \b/n} )\otimes \E^{^{\rm T}}_{\theta/\b}\right]\\
	~\left\{\left[(B\otimes \E_{\a/p\times e/q})
	(C\otimes \E_{e/r\times \b/s})\right] \otimes \E_{\theta/\a}\right\}\\
	=(A\otimes \E_{\a/m\times \theta/n} )
	(B\otimes \E_{\theta/p\times e/q} )
	(C\otimes \E_{e/r\times \b/s})\\
	=(A\otimes \E_{\a/m\times \mu/n} )
	(B\otimes \E_{\mu/p\times \lambda/q} )
	(C\otimes \E_{\lambda/r\times \b/s}).\\
\end{array}
$$
The conclusion follows.
\hfill $\Box$

Before end of this section, we would like to summarize the new algebraic structures over MVMDs by introducing STPs and STAs.

\begin{rem}\label{r2.3.11} There are some new algebraic structures, which are caused by STPs/STAs.
	\begin{itemize}
		\item[(i)] New Algebraic Structures:
		\begin{itemize}
			\item $({\cal M},\ltimes)$ is a monoid with identity $1$.
			\item $({\cal M},\circ)$ is a semi-group.
			\item $({\cal M}_{m\times n}, +, \ttimes)$ is a ring.
			\item $({\cal M}, \hat{\ttimes})$ is a semi-group.
			\item $({\cal M}_{\mu},\hat{\ttimes})<({\cal M}, \hat{\ttimes})$ is a sub-semi-group.
		\end{itemize}
		
		\item[(ii)] Observing the aforementioned STAs and STPs over MVMDs, in many cases it seems difficult to introduce some algebraic structures on them, because the MVMDs always have ``multiple identities".  Recall the pseudo-metric space introduced in the Introduction, we may introduce a new algebraic structure, which allows multiple identities. The ``hyper group", which will be proposed in next section is one of such objects.
	\end{itemize}
\end{rem}

\section{Hyper Group}

\subsection{Lattice}

As a preparation, we briefly review the basic concepts of lattice.

We refer to \cite{bur81} or \cite{fan14} for some basic concepts and properties of lattice.

\begin{dfn}\label{d3.0.1} A partial order set $\Lambda$ is  a lattice, if for any two elements $\lambda,~\mu \in \Lambda$ there exist a least common upper bound, denoted by $\sup(\lambda,\mu)$ (or $\lambda\vee \mu$), and a greatest common lower bound, denoted by $\inf(\lambda,\mu)$ (or $\lambda\wedge \mu$).
\end{dfn}

\begin{figure}
	\centering
	\setlength{\unitlength}{6mm}
	\begin{picture}(12,7)\thicklines
		\put(1,3){\line(1,1){2}}
		\put(1,3){\line(1,-1){2}}
		\put(3,3){\line(1,1){1}}
		\put(3,3){\line(1,-1){1}}
		\put(3,3){\line(-1,1){1}}
		\put(3,3){\line(-1,-1){1}}
		\put(5,3){\line(-1,1){2}}
		\put(5,3){\line(-1,-1){2}}
		\put(1,3){\circle*{0.3}}
		\put(2,2){\circle*{0.3}}
		\put(2,4){\circle*{0.3}}
		\put(3,1){\circle*{0.3}}
		\put(3,3){\circle*{0.3}}
		\put(3,5){\circle*{0.3}}
		\put(4,2){\circle*{0.3}}
		\put(4,4){\circle*{0.3}}
		\put(5,3){\circle*{0.3}}
		\put(2.2,6.2){$Lattice$}
		\put(8,6.2){$Not Lattice$}
		\put(2.8,4.5){$e$}
		\put(4.2,4.0){$f$}
		\put(5.2,3.0){$g$}
		\put(2.8,3.3){$h$}
		\put(4.2,4.0){$f$}
		\put(1.5,1.8){$i$}
			\put(0.5,2.8){$a$}
			\put(1.5,3.8){$c$}
			\put(4.2,1.5){$b$}
			\put(2.8,0.4){$d$}
		\put(8,4){\line(1,1){1}}
		\put(7,3){\line(1,-1){2}}
		\put(9,3){\line(1,1){1}}
		\put(9,3){\line(1,-1){1}}
		\put(9,3){\line(-1,1){1}}
		\put(9,3){\line(-1,-1){1}}
		\put(11,3){\line(-1,1){2}}
		\put(10,2){\line(-1,-1){1}}
		\put(7,3){\circle*{0.3}}
		\put(8,2){\circle*{0.3}}
		\put(8,4){\circle*{0.3}}
		\put(9,1){\circle*{0.3}}
		\put(9,3){\circle*{0.3}}
		\put(9,5){\circle*{0.3}}
		\put(10,2){\circle*{0.3}}
		\put(10,4){\circle*{0.3}}
		\put(11,3){\circle*{0.3}}
		\put(2.2,6.2){$Lattice$}
		\put(8,6.2){$Not  Lattice$}
			\put(6.5,2.8){$a$}
			\put(11.2,2.8){$b$}
	\end{picture}
	
	\caption{Hasse Diagram \label{Fig3.1.1}}
\end{figure}
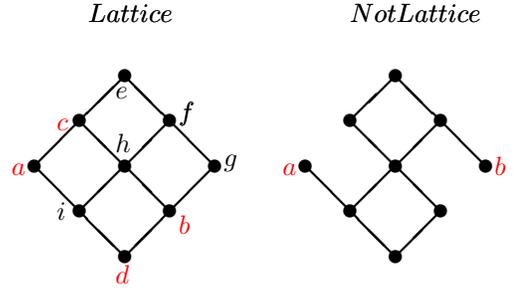

A lattice can be described by a graph called the Hasse diagram.  A Hasse diagram is a graph with no horizontal sides. The order of nodes is decided by oblique sides. For instance, in Figure \ref{Fig3.1.1} the left is a lattice, say, $a\vee b=c$ and $a\wedge b=d$. The right is not a littice, say neither the $a\vee b$ no the $a\wedge b$ exists.

In the following example, we give two lattices which will be used in the sequel.

\begin{exa}\label{e3.0.2}
	\begin{itemize}
		\item[(i)] Consider $\Z_+$. Set  $a\prec b$  if  $a\neq b$ and $a|b$, i.e., $a$ is a proper factor of $b$. It follows that
		$$
		a\vee b=\lcm(a,b);\quad a\wedge b=\gcd(a,b).
		$$
		Then $(\Z_+,\prec)$ is a lattice. We call it the multiply-division lattice of dimension 1, briefly, MD-1 lattice.
		\item[(ii)] Consider $\Z_+\times \Z_+$. Set $(a,c)\prec (b,d)$ if both $a\prec b$ and $c\prec d$ (defined as in (i)). It follows that
		$$
		\begin{array}{l}
			(a,c)\vee (b,d)=(\lcm(a,b),\lcm(c,d));\\
			(a,c)\wedge (b,d)=(\gcd(a,b),\gcd(c,d)).
		\end{array}
		$$
		Then $(\Z_+\times \Z_+,\prec)$ is a lattice, called the MD-2 lattice. Similarly, we can defined the MD-k lattice for $k>2$.
	\end{itemize}
\end{exa}

\begin{dfn}\label{d3.0.3} Let $L$ and $H$ be two lattices.
	\begin{itemize}
		\item[(i)] A mapping $\pi: L\ra H$ is  a lattice homomorphism if
		\begin{align}\label{3.0.1}
			\pi(\lambda\vee \mu)=\pi(\lambda)\vee \pi(\mu),\quad \lambda,~\mu\in L;
		\end{align}
		and
		\begin{align}\label{3.0.2}
			\pi(\lambda\wedge \mu)=\pi(\lambda)\wedge \pi(\mu).
		\end{align}
		\item[(ii)] A homomorphism $\pi: L\ra H$ is  a lattice isomorphism if
		it is  bijective and its inverse is also a homomorphism.
		\item[(iii)] A mapping $\pi: L\ra H$ is an order-preserving mapping, if
		\begin{align}\label{3.0.3}
			a\prec_L b \Rightarrow \pi(a)\prec_H \pi(b), \quad a,b\in L.
		\end{align}
	\end{itemize}
\end{dfn}

The following result comes from definition immediately.

\begin{prp}\label{p3.0.4} Let $L$ and $H$ be two lattices, and $\pi: L\ra H$ be a bijective mapping.
	If both $\pi$ and $\pi^{-1}$ are order-preserving mappings, then $\pi$ is a lattice isomorphism.
\end{prp}

\begin{dfn}\label{d3.0.5} Let $(L,\sup,\inf)$ be a lattice.
	\begin{itemize}
		\item[(i)] $H\subset L$, $H$ is  a sublattice of $L$, if $(H,\sup,\inf)$ is also a lattice.
		\item[(ii)] $H$ is a sublattice of $L$. If for any $a\in L$, there exists an $h\in H$ such that $a\prec h$ implies $a\in H$, then $H$
		is called an ideal of $L$.
		\item[(iii)] $H$ is a sublattice of $L$. If for any $a\in L$, there exists an $h\in H$ such that $a\succ h$ implies $a\in H$, then $H$
		is called a filter of $L$.
	\end{itemize}
\end{dfn}

\subsection{Hyper Group}

This subsection proposes a ``group" with a set of identities, called the identity set. To begin with, we define an identity set.

\begin{dfn}\label{d3.1.3}  Consider a semi-group ${\cal G}=(G,*)$ .
	${\bf e}\subset G$ is called an identity set of ${\cal G}$, if
	\begin{itemize}
		\item[(i)]
		\begin{align}\label{3.1.1}
			g*e=e*g,\quad \forall g\in G, ~\forall e\in {\bf e}.
		\end{align}
		\item[(ii)] There is a lattice $\Lambda$ as the index set of ${\bf e}$,  i.e.,
		$$
		{\bf e}=\{e_{\lambda}\;|\; \lambda \in \Lambda\}.
		$$
		Moreover,
		\begin{align}\label{3.1.2}
			e_{\lambda}*e_{\mu}=e_{\lambda\vee \mu},\quad \lambda,\mu\in \Lambda.
		\end{align}
		
		\item[(iii)] For each $g\in G$ there exists a unique $e_g=e_{\lambda_g}\in {\bf e}$ such that
		\begin{align}\label{3.1.3}
			e_{\lambda} *g=g,
		\end{align}
		if and only if, $\lambda\prec \lambda_g$ (including $\lambda=\lambda_g$).
	\end{itemize}
\end{dfn}

\begin{dfn}\label{d3.1.5}  Consider a semi-group ${\cal G}=(G,*)$ .
	\begin{itemize}
		\item[(i)]  ${\cal G}=(G,*)$ is  a hyper-monoid, if
		there exists an identity set ${\bf e}\subset G$.
		
		\item[(ii)] A hyper-monoid ${\cal G}=(G,*)$ is  a hyper group, if
		for each $g\in G$, there exists a $g^{-1}\in G$ such that
		\begin{align}\label{3.1.4}
			g*g^{-1}=g^{-1}*g=e_{g}.
		\end{align}
	\end{itemize}
	In addition, if $a*b=b*a,\quad \forall a,b\in G$, then it is called an abelian hyper-monoid/hyper group.
\end{dfn}

\subsection{Hyper Groups Over MVMDs}

Recall the additions over MVMDs, defined in Definitions \ref{d2.2.5} and \ref{d2.2.7}, we have the following result.

\begin{prp}\label{p3.2.1}
	\begin{itemize}
		\item[(i)] $(\R^{\infty}, \bar{+})$  with $\bar{+}$  defined by (\ref{2.2.6}) and
		$(\R^{\infty}, \hat{+})$  with $\hat{+}$  defined by (\ref{2.2.11}) are hyper groups,
		where the identity set is
		\begin{align}\label{3.2.1}
			{\bf e}=\{0_n\in \R^n\;|\; n\in \Z_+\},
		\end{align}
		with respect to MD-1 lattice.
		
		\item[(ii)]  $({\cal M}_1,\bar{+})$ is a hyper group with $\bar{+}$ defined by (\ref{2.2.7}), where the identity set is
		\begin{align}\label{3.2.101}
			{\bf e}=\{0_{n\times n} \in {\cal M}_{n\times n} \;|\; n\in \Z_+\},
		\end{align}
		with respect to MD-1 lattice. %with lattice structure $(\Z_+,\prec)$ defined in (i) of Example \ref{e3.0.2}.
		
		Also,   $({\cal M}_1,\hat{+})$ is a hyper group with $\hat{+}$ defined by (\ref{2.2.8}), where the identity set is
		the same as (\ref{3.2.101}).
		
		\item[(iii)] $({\cal M}_{\mu}, \bar{+})$ (or $({\cal M}_{\mu}, \hat{+})$) is a hyper group  with $\bar{+}$ defined by (\ref{2.2.7}) (or $\hat{+}$ defined by (\ref{2.2.8})), where the identity set is
		\begin{align}\label{3.2.102}
			{\bf e}=\{0_{n\mu_y\times n\mu_x} \in {\cal M}_{n\mu_y\times n\mu_x} \;|\; n\in \Z_+\},
		\end{align}
		with respect to MD-1 lattice. %with lattice structure $(\Z_+,\prec)$ defined in (i) of Example \ref{e3.0.2}.
		
		\item[(iv)]
		$({\cal M}, \hat{+})$ is a hyper group with $\hat{+}$  defined by (\ref{2.2.11}), where the identity set is
		\begin{align}\label{3.2.2}
			{\bf e}=\{0_{m\times n}\in {\cal M}_{m\times n}\;|\; (m, n)\in \Z_+\times \Z_+\},
		\end{align}
		with respect to MD-2 lattice. %with lattice structure $(\Z_+\times \Z_+,\prec)$ defined in (ii) of Example \ref{e3.0.2}.
		
	\end{itemize}
\end{prp}

\noindent{\it Proof.}
We prove (iv) only. The proofs for others are similar.

First, we show that $({\cal M},\hat{+})$ is a semi-group. Assume $A\in {\cal M}_{m
	\times n}$, $A\in {\cal M}_{m \times n}$, $B\in {\cal M}_{p \times q}$, $C\in {\cal M}_{r \times s}$, and
$$
\begin{array}{l}
	m\vee p=a,~p\vee r=c,~m\vee p\vee r=\xi,\\
	n\vee q=b,~q\vee s=d,~n\vee q\vee s=\eta.
\end{array}
$$
Then
$$
\begin{array}{l}
	(A\hat{+} B)\hat{+} C=\left[(A\otimes \E_{a/m\times b/n})+(B\otimes \E_{a/p\times b/q})\right]\hat{+} C\\
	~=\left[(A\otimes \E_{a/m\times b/n})+(B\otimes \E_{a/p\times b/q})\right]\otimes \E_{\xi/a\times \eta/b}\\
	~+C\otimes \E_{\xi/r\times \eta/s}\\
	~=(A\otimes \E_{\xi/m\times \eta/n})+(B\otimes \E_{\xi/p\times \eta/q})+(C\otimes \E_{\xi/r\times \eta/s})\\
\end{array}
$$
A similar calculation shows that
$$
\begin{array}{l}
	A\hat{+} (B\hat{+} C)=(A\otimes \E_{\xi/m\times \eta/n})\\
	~+(B\otimes \E_{\xi/p\times \eta/q})+(C\otimes \E_{\xi/r\times \eta/s})\\
\end{array}
$$
Hence $({\cal M},\hat{+})$ is a semi-group, which is obviously abelian.

Then we consider the identity set ${\bf e}$.

Let $A\in {\cal M}_{m\times n}$, $0_{p\times q}\in {\bf e}$, $m\vee p=a$, and $n\vee q=b$. Then
a straightforward computation shows that
$$
A\hat{+} 0_{p\times q}=0_{p\times q}\hat{+} A=A\otimes \E_{a/m\times b/n}.
$$
(\ref{3.1.1}) holds.

Next, it is ready to verify that the lattice structure of (\ref{3.2.2}) satisfies (\ref{3.1.2}).

Finally, assume $A$ as before, then  it is obvious that $e_A=0_{m\times n}\in {\bf e}$, which leads to
$$
A+(-A)=e_A.
$$
(\ref{3.1.3}) is satisfied.

\hfill $\Box$

\begin{rem}\label{r3.2.2}
	\ All the additive hyper groups are abelian.
\end{rem}

Next, we consider the product hyper groups.

Set
$$
{\cal T}:=\bigcup_{i=1}^{\infty}{\cal T}_n.
$$
Then we have the product hyper group of matrices as follows.

\begin{prp}\label{p3.2.3}
	$({\cal T}, \ltimes)$ is a hyper group with $\ltimes$  defined by (\ref{2.2.1}), where the identity set is
	\begin{align}\label{3.2.3}
		{\bf e}=\{I_n\in {\cal M}_{n\times n}\;|\; n\in \Z_+\},
	\end{align}
	with respect to MD-1 lattice.
\end{prp}

\noindent{\it Proof.}

It is well known that $({\cal M}.\ltimes)$ is a semi-group. Let $A,B\in {\cal T}$. It is easy to verify that if $A\ltimes B\in {\cal T}$. Hence, $({\cal T},\ltimes)$ is a sub-semi-group.

Let $A\in {\cal T}_m$, $I_n\in {\bf e}$, and $t=m\vee n$. Then
$$
A\ltimes I_n=I_n\ltimes A=A\otimes I_{t/m}.
$$
(\ref{3.1.1}) holds.

It is obvious that $e_A=I_m$ and
$$
A\ltimes A^{-1}=AA^{-1}=I_m.
$$
(\ref{3.1.2}) holds.

Finally, the lattice structure of the identity set is obvious. Moreover,
$$
I_m\ltimes I_n=I_t.
$$
(\ref{3.1.3}) is also true.

\hfill $\Box$

To emphasize the identity set of a hyper group, we may also denote a hyper group by a triple $(G,*,{\bf e})$.

\subsection{Component Groups and Equivalence  Group}

This subsection considers the component groups and equivalence group  of a hyper group.

\begin{prp}\label{p3.3.1} Assume ${\cal G}=(G,*, {\bf e})$ is a hyper group, where ${\bf e}$ is determined by a lattice $\Lambda$.
	Assign
	for each $e_{\lambda}\in {\bf e}$, $\lambda\in \Lambda$ a set $G_{\lambda}\subset G$ as
	\begin{align}\label{3.3.1}
		G_{\lambda}:=\{x\in G\;|\; e_x=e_{\lambda}\}.
	\end{align}
	Then
	\begin{itemize}
		\item[(i)] For each $\lambda\in \Lambda$,  ${\cal G}_{\lambda}=(G_{\lambda},*)$ is a group, called the component group of ${\cal G}$, and its identity is $e_{\lambda}$.
		\item[(ii)]
		\begin{align}\label{3.3.2}
			G=\bigcup_{\lambda\in \Lambda}G_{\lambda}
		\end{align}
		is a partition.
		\item[(iii)] If $x\in G_{\lambda}$ and $y\in G_{\mu}$, then
		$x*y\in G_{\lambda \vee \mu}$.
	\end{itemize}
\end{prp}

\noindent {\it Proof.}
\begin{itemize}
	\item[(i)] To see $G_{\lambda}$ is a group, we have to show that if $x,y\in G_{\lambda}$, then $x*y\in G_{\lambda}$ . It obvious that
	$$x*y*e_{\lambda}=e_{\lambda}*x*y.
	$$
	We also have
	$$
	x*y*y^{-1}x^{-1}=x*e_{\lambda}x^{-1}=e_{\lambda},
	$$
	and similarly, we have
	$$
	y^{-1}x^{-1}*x*y=e_{\lambda}.
	$$
	Hence, we have that both $(x*y)^{-1}=y^{-1}*x^{-1}$ and $x*y$ are $\in G_{\lambda}$.
	
	It is clear from definition that if $x\in G_{\lambda}$ then $x^{-1}\in G_{\lambda}$.
	
	Then the associativity and the uniqueness of the inverse are inherited from the properties of hyper group.
	
	We conclude that $G_{\lambda}$ is a group.
	
	\item[(ii)] The partition (\ref{3.3.2}) comes from definition immediately.
	
	\item[(iii)] Since
	$$
	\begin{array}{l}
		x*y*y^{-1}*x^{-1}=x*e_{\mu}*x^{-1}=e_{\mu}*x*x^{-1}\\
		=e_{mu}*e_{\lambda}=e_{\lambda\vee \mu}.
	\end{array}
	$$
	By definition $x*y\in G_{\lambda\vee \mu}$.
\end{itemize}
\hfill $\Box$

\begin{dfn}\label{d3.3.2} Consider a hyper group ${\cal G}=(G,*, {\bf e})$ .
	Let $x,~y\in G$.  $x$ and $y$ are said to be equivalent, denoted by $x\sim y$, if
	there exist $e_{\a}$ and $e_{\b}$ such that
	\begin{align}\label{3.3.3}%{2.1.7}
		x*e_{\a}=y*e_{\b}.
	\end{align}
	
	The equivalence class of $x$, denoted by $\bar{x}$, is
	$$
	\bar{x}:=\{y\in G\;|\; y\sim x\}.
	$$
\end{dfn}

Note that we have to show that (\ref{3.3.3}) is an equivalence relation.  A straightforward verification ensures this.

\begin{dfn}\label{d3.3.3} Consider a hyper group ${\cal G}=(G,*)$. The equivalence $\sim$ is said  to be consistent with the group operator $*$, if the followings are satisfied.
	\begin{itemize}
		\item[(i)] If $x_1\sim x_2$ and $y_1\sim y_2$, then
		\begin{align}\label{3.3.4}%{2.1.8}
			x_1*y_1\sim x_2*y_2.
		\end{align}
		\item[(ii)] If there exists $e\in {\bf e}$ such that $x\sim e$, then $x\in {\bf e}$.
	\end{itemize}
\end{dfn}

Then we have the identity set as one element.

\begin{prp}\label{p3.3.4} Consider a hyper group ${\cal G}=(G,*, {\bf e})$. If the equivalence is consistent with respect to the operator, then
	\begin{align}\label{3.3.5}
		\bar{e}={\bf e},\quad e\in {\bf e}.
	\end{align}
\end{prp}

\noindent{\it Proof.} By definition we have
$$
e_{\lambda}\sim e_{\mu},\quad \lambda,\mu\in \Lambda,
$$
According to (ii) of Definition \ref{d3.3.3}, we also know that if $a\not\in {\bf e}$ then for any $e\in {\bf e}$, $a\not\sim e$.
(\ref{3.3.5}) follows.
\hfill $\Box$

Using Proposition \ref{p3.3.4}, we have the following result immediately.

\begin{prp}\label{p3.3.5} Consider a hyper group ${\cal G}=(G,*,{\bf e})$. Assume the equivalence $\sim$ defined by (\ref{3.3.3}) is consistent with $*$. Then
	\begin{align}\label{3.3.6}%{2.1.9}
		\bar{x}*\bar{y}:=\overline{x*y},
	\end{align}
	is properly defined. Moreover, $\bar{\cal G}=(\bar{G},*)$ is a group, called the equivalence group of ${\cal G}$, denoted by
	$$
	\bar{\cal G}:={\cal G}/\sim=(\bar{G},*),
	$$
	where $\bar{G}=\{\bar{x}\;|\; x\in G\}$.
\end{prp}

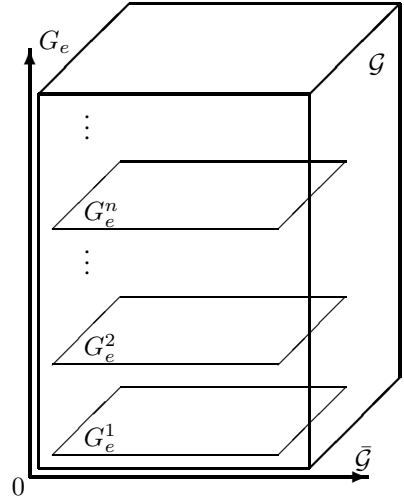
\begin{figure}
	\centering
	\setlength{\unitlength}{6mm}
	\begin{picture}(10,12)(-0.5,-0.5)
		\thicklines
		\put(6.2,0.2){\line(1,1){2}}
		\put(8.2,2.2){\line(0,1){8.3}}
		\put(0.2,0.2){\line(1,0){6}}
		\put(0.2,0.2){\line(0,1){8.3}}
		\put(0.2,8.5){\line(1,0){6}}
		\put(0.2,8.5){\line(1,1){2}}
		\put(2.2,10.5){\line(1,0){6}}
		\put(6.2,8.5){\line(1,1){2}}
		\put(6.2,8.5){\line(0,-1){8.3}}
		\put(0,0){\vector(0,1){9.5}}
		\put(0,0){\vector(1,0){7.5}}
		\put(0.2,9.5){$G_e$}
		\put(7.2,0.2){$\bar{{\cal G}}$}
		\put(1.2,0.7){$G^1_e$}
		\put(1.2,2.7){$G^2_e$}
		\put(1.2,5.7){$G^n_e$}
		\put(1.2,4.5){$\vdots$}
		\put(1.2,7.5){$\vdots$}
		\put(7.5,9){${\cal G}$}
		\put(-0.4,-0.4){$0$}
		\thinlines
		\put(0.5,0.5){\line(1,0){5}}
		\put(0.5,0.5){\line(1,1){1.5}}
		\put(7,2){\line(-1,0){5}}
		\put(7,2){\line(-1,-1){1.5}}
		\put(0.5,2.5){\line(1,0){5}}
		\put(0.5,2.5){\line(1,1){1.5}}
		\put(7,4){\line(-1,0){5}}
		\put(7,4){\line(-1,-1){1.5}}
		\put(0.5,5.5){\line(1,0){5}}
		\put(0.5,5.5){\line(1,1){1.5}}
		\put(7,7){\line(-1,0){5}}
		\put(7,7){\line(-1,-1){1.5}}
	\end{picture}
	\caption{Group Decomposition of a Hyper Group ${\cal G}$\label{Fig3.3.1}}
\end{figure}

Observing Propositions \ref{p3.3.1} and \ref{p3.3.5}, a diagram of a hyper group with its  group decomposition  is depicted by Figure \ref{Fig3.3.1}.

Since the existence of the equivalence hyper group depends on the consistence of the equivalence with respect to the group operator, the consistence is the key issue for the structure analysis of hyper groups. It is worth showing the consistence of  the hyper groups over MVMDs. To this end we need a preparation.
%%%%%%%%%%%%%%%%%%%%%%%%%%%%%%%%%%%%%%%%%%%%%%%%%%%%%%%
\begin{lem}\label{l3.3.6}%\label{l7.2.2}
	\begin{itemize}
		\item[(i)] Consider $x,y\in \R^{\infty}$, which are specified  by $x\in \R^m$, $y\in \R^n$, $m\vee n=t$, and $m\wedge n=s$.
		Assume there are $\J_p$ and $\J_q$ (or $\E_p$ and $\E_q$) such that
		\begin{align}\label{3.3.7}%{7.2.1}
			x\otimes \J_p=y\otimes \J_q,~(\mbox{or}~x\otimes \E_p=y\otimes \E_q),
		\end{align}
		then there exists a $z\in \R^s$ such that
		\begin{align}\label{3.3.8}%{7.2.2}
			x=z\otimes \J_{m/s};\quad y=z\otimes \J_{n/s}.
		\end{align}
		\item[(ii)] Consider $A,B\in {\cal M}_{\mu}$, which are specified  by $A\in {\cal M}_{m\mu_y\times m\mu_x}$, $B\in {\cal M}_{n\mu_y\times n\mu_x}$, $m\vee n=a$,  and $m\wedge n=b$.
		Assume there are $H_r=I_r$ and  $H_s=I_s$  (or  $H_r=J_r$ and  $H_s=J_s$ )  such that
		\begin{align}\label{3.3.9}%{7.2.3}
			A\otimes H_r=B\otimes H_s,
		\end{align}
		then there exists a $C\in {\cal M}_{b\mu_y\times b\mu_x}$ such that
		\begin{align}\label{3.3.10}%{7.2.4}
			A=C\otimes H_{m/b};\quad B=C\otimes H_{n/b}.
		\end{align}
		\item[(iii)] Consider $A,B\in {\cal M}$, which are specified  by $A\in {\cal M}_{m\times n}$, $B\in {\cal M}_{p\times q}$, $m\wedge p=r$,  and $n\wedge q=s$.
		Assume there are $\E_{a\times b}$ and $\E_{c\times d}$   such that
		\begin{align}\label{3.3.11}%{7.2.5}
			A\otimes \E_{a\times b}=B\otimes \E_{c\times d},
		\end{align}
		then there exists a $C\in {\cal M}_{r\times s}$ such that
		\begin{align}\label{3.3.12}%{7.2.6}
			A=C\otimes \E_{m/r\times n/s};\quad B=C\otimes \E_{p/r\times q/s}.
		\end{align}
	\end{itemize}
\end{lem}

\noindent{\it Proof.} In fact, (i) and (ii) have been proved in \cite{che19b}. The proof of (iii) is similar. We sketch it.
First, note that $ma=pc$ and $nb=qd$. By deleting common factors we can assume, without loss of generality, that
$$
a=\frac{m\vee p}{m},\quad b=\frac{n\vee q}{n},\quad c=\frac{m\vee p}{p},\quad d=\frac{n\vee q}{q}.
$$
Using the above $a,b,c,d$ and comparing both sides of (\ref{3.3.11}) carefully, the conclusion follows. The detailed arguments are similar to those for (i) and (ii) in \cite{che19b}.
\hfill $\Box$

\begin{exa}\label{e3.3.7} Consider the hyper group ${\cal G}=({\cal T},\ltimes, {\bf e})$, where the identity set is defined by
	(\ref{3.2.3}).
	Then we have
	\begin{itemize}
		\item[(i)] The component groups are
		\begin{align}\label{3.3.13}%{3.3.7}
			{\cal T}_n=\GL(n,\R),\quad  n\in \Z_+,
		\end{align}
		where $\GL(n,\R)$ are the standard general linear groups, and
		$$
		{\cal T}:=\bigcup_{i=1}^{\infty}\GL(n,\R).
		$$
		\item[(ii)]
		
		Using Lemma \ref{l3.3.6}, it is easy to verify that the equivalence is consistent with the STP, hence the equivalence quotient group exists.
		$$
		\overline{\cal T}:=\bigcup_{n=1}^{\infty} \{\bar{A}_n\;|\;A_n\in {\cal M}_{n\times n} ~\mbox{is invertible.}\}.
		$$
		If we define
		$\pi: {\cal T}\ra \overline{\cal T}$ by
		$$
		\pi(A)=\bar{A},\quad A\in {\cal T},
		$$
		then it is easy to verify that $\pi$ is a hyper group homomorphism.
	\end{itemize}
\end{exa}

In the following we consider the additive hyper groups. Recall the additive hyper groups discussed in Proposition \ref{p3.2.1}.

\begin{exa}\label{e3.3.8} In the following we assume
	the consistence of the equivalence with respect to addition in each case holds.
	
	\begin{itemize}
		\item[(i)] Consider $(\R^{\infty}, \bar{+})$  with $\bar{+}$  defined by (\ref{2.2.6}).
		It is easy to verify that the component groups are $(\R^n,+)$, $n\in \Z_+$.
		
		Next, consider the equivalence hyper group.
		Let $x\in \R^m$, $y\in \R^n$, and $t=\lcm(m.n)$. Assume $x$ and $y$ are equivalent, denoted by $x\lra y$. By definition,
		there exist $0_p$ and $0_q$, such that
		\begin{align}\label{3.3.14}%{3.3.8}
			x \bar{+} 0_p=y\bar{+} 0_q.
		\end{align}
		Comparing the dimensions and deleting common factor $\J_s$, where $s=p\wedge q$, we have
		\begin{align}\label{3.3.15}%{3.3.9}
			x\otimes \J_{t/m}=y\otimes \J_{t/n}.
		\end{align}
		In fact, (\ref{3.3.15}) is the definition of equivalence of two vectors in \cite{che19}. Note that the consistence of the equivalence with respect to $\bar{+}$ has already been proved in \cite{che19}, then the equivalence group exists, which is denoted by
		\begin{align}\label{3.3.16}%{3.3.10}
			\Omega=(\R^{\infty}/\lra)=(\{\bar{x}\;|\; x\in \R^{\infty}\}, \bar{+}).
		\end{align}
		
		Similarly, we also have $x \leftrightharpoons y$, if
		\begin{align}\label{3.3.17}%{3.3.11}
			x \hat{+} 0_p=y\hat{+} 0_q,
		\end{align}
		which leads to
		\begin{align}\label{3.3.18}%{3.3.12}
			x\otimes \E_{t/m}=y\otimes \E_{t/n}.
		\end{align}
		Finally, we have the equivalence group as
		\begin{align}\label{3.3.19}%{3.3.13}
			\Omega'=(\R^{\infty}/\leftrightharpoons )=(\{\bar{x}\;|\; x\in \R^{\infty}\}, \hat{+}).
		\end{align}
		
		\item[(ii)]  Consider $({\cal M}_1,\bar{+})$. Let $A\in {\cal M}_{m\times m}$,  $B\in {\cal M}_{n\times n}$, and $t=\lcm(m,n)$.
		Then $A\sim B$ if there exist $0_{p\times p}$ and $0_{q\times q}$ such that
		\begin{align}\label{3.3.20}%{3.3.14}
			A \bar{+} 0_{p\times p}=B\bar{+} 0_{q\times q},
		\end{align}
		which leads to
		\begin{align}\label{3.3.21}%{3.3.15}
			A\otimes I_{t/m}=B\otimes I_{t/n}.
		\end{align}
		Finally, we have the equivalence group as
		\begin{align}\label{3.3.22}%{3.3.16}
			\Sigma=({\cal M}_1/\sim) =(\{\bar{A}\;|\; A\in {\cal M}_1\}, \bar{+}).
		\end{align}
		
		Similarly, consider $({\cal M}_1,\hat{+})$. We have $A\approx B$ if
		\begin{align}\label{3.3.23}%{3.3.17}
			A \hat{+} 0_{p\times p}=B\hat{+} 0_{q\times q},
		\end{align}
		which leads to
		\begin{align}\label{3.3.24}%{3.3.18}
			A\otimes J_{t/m}=B\otimes J_{t/n}.
		\end{align}
		Finally, we have the equivalence group as
		\begin{align}\label{3.3.25}%{3.3.19}
			\Sigma'=({\cal M}_1/\approx) =(\{\bar{A}\;|\; A\in {\cal M}_1\}, \hat{+}).
		\end{align}
		
		\item[(iii)] Consider $({\cal M}_{\mu}. \hat{+})$.  Let $A\in {\cal M}_{m\mu_y\times m\mu_x}$, $B\in {\cal M}_{n\mu_y\times n\mu_x}$,
		and $t=\lcm(m,n)$. $A\approx B$ if
		\begin{align}\label{3.3.26}%{3.3.20}
			A \hat{+} 0_{p\mu_y\times p\mu_x}=B\hat{+} 0_{q\nu_y\times q\mu_x},
		\end{align}
		which leads to
		\begin{align}\label{3.3.27}%{3.3.21}
			A\otimes J_{(t/m)\mu_y\times (t/m)\mu_x}=B\otimes J_{(t/n)\mu_y\times (t/n)\mu_x}.
		\end{align}
		Finally, we have
		\begin{align}\label{3.3.28}%{3.3.22}
			\Sigma_{\mu}=({\cal M}_{\mu}/\approx) =(\{\bar{A}\;|\; A\in {\cal M}_{\mu}\}, \hat{+}).
		\end{align}
		
		\item[(iv)]
		Consider $({\cal M}, \hat{+})$. Let $A\in {\cal M}_{m\times n}$, $B\in {\cal M}_{p\times q}$, $A\sim B$ if there exist $0_{r\times s}$ and $0_{a\times b}$ such that
		\begin{align}\label{3.3.29}%{3.3.23}
			A \hat{+} 0_{r\times s}=B\hat{+} 0_{a\times b}.
		\end{align}
		Assume $\a=m\vee p$ and $\b=n\vee q$, then (\ref{3.3.23}) leads to
		\begin{align}\label{3.3.30}%{3.3.24}
			A\otimes \E_{(\a/m)\times (\b/n)}=B\otimes \E_{(\a/p)\times (\b/q) n)}.
		\end{align}
		Finally, we have
		\begin{align}\label{3.3.31}
			\Sigma=({\cal M}/\approx) =(\{\bar{A}\;|\; A\in {\cal M}\}, \hat{+}).
		\end{align}
		
	\end{itemize}
\end{exa}

In the above example we assume the consistence always exist. Using Lemma \ref{l3.3.6}, this can be proved one by one:

\begin{prp}\label{p3.3.9}%7.2.3}
\begin{itemize}
	\item[(i)] Consider $(\R^{\infty}, \bar{+})$ (or $(\R^{\infty}, \hat{+})$), its equivalence is consistent with its operator $\bar{+}$ ( or $\hat{+}$).
	
	\item[(ii)]  Consider $({\cal M}_1,\bar{+})$ (or
	$({\cal M}_1,\hat{+})$ ),  its equivalence is consistent with its operator $\bar{+}$ (or $\hat{+}$).
	
	\item[(iii)] Consider $({\cal M}_{\mu}, \hat{+})$,  its equivalence is consistent with its operator $\hat{+}$.
	
	\item[(iv)]
	Consider $({\cal M}, \hat{+})$,  its equivalence is consistent with its operator $\hat{+}$.
\end{itemize}
\end{prp}

\noindent{\it Proof.} We prove (i) and (iv). The proofs for (ii) and (iii) are similar.
\begin{itemize}
\item[(i)]
First, to use Lemma \ref{l3.3.6}, we show that if $x,y\in \R^{\infty}$ and $x\lra y$, then there exist $\J_p$ and $\J_q$ such that
$$
x\otimes \J_p=y\otimes \J_q.
$$
Assume $x\in \R^m$ and $y\in \R^n$. By definition, there exist $0_p$ and $0_q$, such that
$$
\begin{array}{l}
	x\bar{+} 0_p=(x\otimes \J_{(m\vee p)/m})+(0_p\otimes \J_{(m\vee p)/p})\\
	(x\otimes \J_{(m\vee p)/m})\\
	=y\bar{+} 0_q=(y\otimes \J_{(n\vee q)/n}).\\
\end{array}
$$
Set $t=m\vee n$. Then the above equality can be simplified to
$$
\begin{array}{l}
	x\otimes \J_{t/m} =y\otimes \J_{t/n}.
\end{array}
$$

Next, assume $x_1\sim x_2$ and $y_1\sim y_2$, we have to show that
\begin{align}\label{3.3.32}%{7.2.7}
	x_1 \bar{+} y_1\sim x_2 \bar{+} y_2.
\end{align}

Suppose $x_1\in \R^m$, $x_2\in \R^n$, $y_1\in \R^p$, $y_2\in \R^q$, and
$$
\begin{array}{ccc}
	m\vee p=\a,&n\vee q=\b,& m\wedge n=r,\\
	p\wedge q=s,&r\vee s=\xi.&\a=\xi a,\\
	\b=\xi b,&~&\\
\end{array}
$$
According to Lemma \ref{l3.3.6} there exist $x_0\in \R^r$ and $y_0\in \R^s$ such that
$$
\begin{array}{cc}
	x_1=x_0\otimes \J_{m/r},&x_2=x_0\otimes \J_{n/r},\\
	y_1=y_0\otimes \J_{p/s},&y_2=y_0\otimes \J_{q/s}.\\
\end{array}
$$
Then
$$
\begin{array}{l}
	x_1\bar{+}y_1=x_0\otimes \J_{m/r}\otimes \J_{\a/m}+y_0\otimes \J_{p/s}\otimes \J_{\a/p}\\
	~=x_0\otimes \J_{\a/r}+y_0\otimes \J_{\a/s}\\
	~=(x_0\otimes \J_{\xi/r}+y_0\otimes \J_{\xi/s})\otimes \J_{a}.\\
\end{array}
$$
Similarly,
$$
\begin{array}{l}
	x_2\bar{+}y_2=x_0\otimes \J_{\b/r}+y_0\otimes \J_{\b/s}\\
	~=(x_0\otimes \J_{\xi/r}+y_0\otimes \J_{\xi/s})\otimes \J_{b}.\\
\end{array}
$$
We conclude that
$$
x_1\bar{+}y_1\sim (x_0\otimes \J_{\xi/r}+y_0\otimes \J_{\xi/s})\sim x_2\hat{+}y_2.
$$

Next, we check condition (ii) of Definition \ref{d3.3.3}. Assume $x\in \R^m$ and $x\sim 0_a\in {\bf e}$, then there exist $\J_s$ and $\J_t$, such that
$$
x\otimes \J_s=0_a\otimes \J_t=0_{at}.
$$
It is easy to see that $x=0_m\in {\bf e}$.

\item[(iv)] To prove  the consistence of  the addition and the product with the equivalence we need to show that
for $A_1\sim A_2$ and $B_1\sim B_2$ we have
\begin{align}\label{3.3.33}%{6.4.3}
	A_1\hat{+} B_1\sim A_2\hat{+} B_2.
\end{align}

Assume $A_1\in {\cal M}_{m\times n}$, $A_2\in {\cal M}_{p\times q}$, $B_1\in {\cal M}_{\a\times \b}$, and $B_2\in {\cal M}_{\gamma\times \d}$,  and
$$
\begin{array}{ccc}
	m\vee \a=a,&n\vee \b=b,&p\vee \gamma=c,\\
	q\vee \d =d,&m\wedge p=s,&n\wedge q=t,\\
	\a\wedge \gamma =\xi,&\b\wedge \d=\eta,&~\\
	s\vee \xi=u,&t\vee \eta=v.&~\\
\end{array}
$$
Using Lemma \ref{l3.3.6}, we know that there exist
$A_0\in {\cal M}_{s\times t}$ and $B_0\in {\cal M}_{\xi\times \eta}$ such that
$$
\begin{array}{l}
	A_1=A_0\otimes \E_{m/s\times n/t},\quad A_2=A_0\otimes \E_{p/s\times q/t};\\
	B_1=B_0\otimes  \E_{\a/\xi\times \b/\eta},\quad B_2=B_0\otimes \E_{\gamma/\xi\times \d/\eta}.\\
\end{array}
$$
Now
$$
\begin{array}{l}
	A_1\hat{+} B_1=A_1\otimes \E_{a/m\times b/n}+ B_1\otimes \E_{a/\a\times b/\b}\\
	~=A_0\otimes \E_{a/s\times b/t}+B_0\otimes \E_{a/\xi\times b/\eta}\\
	~=( A_0\otimes \E_{u/s\times v/t}+B_0\otimes \E_{u/\xi\times v/\eta})\otimes \E_{a/u\times b/v}.
\end{array}
$$
$$
\begin{array}{l}
	A_2\hat{+} B_2=A_2\otimes \E_{c/p\times d/q}+ B_2\otimes \E_{c/\gamma\times d/\d}\\
	~=A_0\otimes \E_{c/s\times d/t}+B_0\otimes \E_{c/\xi\times d/\eta}\\
	~=( A_0\otimes \E_{u/s\times v/t}+B_0\otimes \E_{u/\xi\times v/\eta})\otimes \E_{c/u\times d/v}.
\end{array}
$$
It follows immediately that
$$
A_1\hat{+}B_1\sim A_2\hat{+} B_2.
$$
\end{itemize}

The condition (ii) of Definition \ref{d3.3.3} is straightforward verifiable.
\hfill $\Box$

We conclude that  for these additive hyper-vector spaces the equivalence vector spaces always exist.

\subsection{Sub Hyper Group}

\begin{dfn}\label{d3.4.1} Let ${\cal G}=(G,*,{\bf e})$ be a hyper group. ${\cal H}=(H,*,{\bf e}_0)$ is called a sub hyper group, if
$H\subset G$ is a subset, ${\bf e}_0\subset {\bf e}$ is a sub-lattice, and $(H,*,{\bf e}_0)$ is a hyper group. Denote the sub hyper group by ${\cal H}<_h {\cal G}$. If
${\cal H}<_h {\cal G}$ and $gH_{e_g}=H_{e_g}g$, $\forall g\in G$ and $e_g\in {\bf e}_0$, then $(H,*,{\bf e}_0)$ is a normal sub hyper group, and denoted by
${\cal H}\lhd_h {\cal G}$.
\end{dfn}

The following result comes from the definition immediately.

\begin{prp}\label{p3.4.2} Let ${\cal G}=(G,*,{\bf e})$ be a hyper group.
\begin{itemize}
	\item[(i)] Assume ${\bf e}_0\prec {\bf e}$ is a sub-lattice. Set
	\begin{align}\label{3.4.1}
		G_0:=\bigcup_{e\in {\bf e}_0}G_e,
	\end{align}
	then
	\begin{align}\label{3.4.2}
		{\cal G}_0:=(G_0,*,{\bf e}_0)<_h (G,*,{\bf e}).
	\end{align}
	\item[(ii)]  Consider ${\cal G}_0$ as in (i). Assume
	$H_e<G_e$, $\forall e\in {\bf e}$, and
	\begin{align}\label{3.4.3}
		H^0:=\bigcup_{e\in {\bf e}_0}H_e.
	\end{align}
	Then
	\begin{align}\label{3.4.4}
		{\cal H}_0:=(H^0, *, {\bf e}_0)<_h {\cal G}_0.
	\end{align}
	%Moreover, if $H_e\lhd G_e$, $\forall e\in {\bf e}$, then
	%\begin{align}\label{3.4.5}
	%{\cal H}:=(H, *, {\bf e}_0)\lhd_h {\cal G}_0.
	%\end{align}
	\item[(iii)] Conversely, if ${\bf e}_0\subset {\bf e}$ is a sublattice, $\emptyset\neq H_e\subset G_e$, $e\in {\bf e}_0$. $H^0$ is defined by (\ref{3.4.3}), and
	${\cal H}_0:=(H^0, *, {\bf e}_0)<_h {\cal G}$, then
	\begin{align}\label{3.4.6}
		{\cal H}_e:=(H_e, *) <{\cal G}_e,\quad \forall e\in {\bf e}_0.
	\end{align}
	\item[(iv)]
	If ${\cal H}_0:=(H^0, *, {\bf e}_0)\lhd_h {\cal G}$, then
	\begin{align}\label{3.4.7}
		{\cal H}_e:=(H_e, *) \lhd {\cal G}_e,\quad \forall e\in {\bf e}_0.
	\end{align}
\end{itemize}
\end{prp}

\begin{figure}
\centering
\setlength{\unitlength}{6mm}
\begin{picture}(18,5)\thicklines
	\put(0,0){\framebox(5,1.5){${\cal H}_0=\{H_e\;|\;e\in {\bf e}_0\}$}}
	\put(7,0){\framebox(5,1.5){${\cal H}=\{H_e\;|\;e\in {\bf e}\}$}}
	\put(0,2.5){\framebox(5,1.5){${\cal G}_0=\{G_e\;|\;e\in {\bf e}_0\}$}}
	\put(7,2.5){\framebox(5,1.5){${\cal G}=\{G_e\;|\;e\in {\bf e}\}$}}
	\put(5.2,3.45){${\bf e}_0\prec {\bf e}$}
	\put(5.7,2.6){$<$}
	\put(5.7,0.6){$<$}
	\put(2.7,1.8){$\vee$}
	\put(9.7,1.8){$\vee$}
	\put(12.5,1.8){$H_e<G_e$}
	\put(7,3.25){\vector(-1,0){2}}
	\put(12.5,3.5){\vector(0,-1){3}}
	\put(2.5,2.5){\vector(0,-1){1}}
	\put(9.5,2.5){\vector(0,-1){1}}
\end{picture}

\caption{Sub Hyper Groups \label{Fig3.4.1}}
\end{figure}
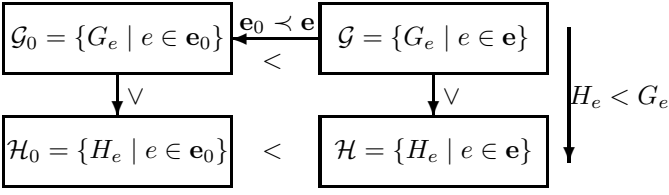

Fig. \ref{Fig3.4.1} depicts these relations.

\begin{exa}\label{e3.4.3}
\begin{itemize}
	\item[(i)] A lattice ${\cal H}$ is constructed by
	$$
	H=\{1,2,3,4,6,9,12,18,36\}\subset Z_+
	$$
	with sublattice structure of MD-1. Then its Hasse diagram is depicted by the left graph of Figure \ref{Fig3.1.1}, where
	$$
	\begin{array}{ccc}
		e=36,&c=12,&f=18,\\
		a=4,&h=6,&g=3,\\
		i=2,&b=3,&d=1.\\
	\end{array}
	$$
	Set
	$$
	{\cal T}_H:=\{{\cal T}_n\;|\;n\in H\},
	$$
	and
	$$
	{\bf e}_h=(\{I_n\;|\;n\in H\}, \vee,\wedge\}).
	$$
	Then
	${\cal H}:=({\cal T}_H, \ltimes, {\bf e}_h)<_h {\cal G}$ is a sub-hyper group, with respect to MD-1 lattice.
	${\cal G}=({\cal T},\ltimes, {\bf e})$ is defined in Example \ref{e3.3.7}.
	
	\item[(ii)] Consider $M_1=\bigcup_{n=1}^{\infty}{\cal M}_{n\times n}$,
	$\hat{+}$ is defined by (\ref{2.2.8}), and
	$$
	{\bf e}=\{0_{n\times n}\in {\cal M}_{n\times n}\;|\; n\in \Z_+\},
	$$
	with lattice structure $(\Z_+,\prec)$ defined in (i) of Example \ref{e3.0.2}.
	Then ${\cal M}_1=(M_1,\hat{+}, {\bf e})$ is a hyper group.
	\begin{itemize}
		\item Set
		$$
		S_{n}:=\{A\in {\cal M}_{n\times n}\;|\; A^{^{\rm T}}=A\},
		$$
		and
		$$
		S=\bigcup_{n=1}^{\infty}S_n.
		$$
		Then ${\cal S}:=(S,\hat{+},{\bf e})$ is a hyper group, and
		$$
		{\cal S}\lhd_h {\cal M}_1.
		$$
		\item Set
		$$
		K_{n}:=\{A\in {\cal M}_{n\times n}\;|\; A^{^{\rm T}}=-A\},
		$$
		and
		$$
		K=\bigcup_{n=1}^{\infty}K_n.
		$$
		Then ${\cal K}:=(K,\hat{+},{\bf e})$ is a hyper group, and
		$$
		{\cal K}\lhd_h {\cal M}_1.
		$$
		\item Set
		$$
		D_{n}:=\{A\in {\cal M}_{n\times n}\;|\; A ~\mbox{is diagonal.}\},
		$$
		and
		$$
		D=\bigcup_{n=1}^{\infty}D_n.
		$$
		Then ${\cal D}:=(D,\hat{+},{\bf e})$ is a hyper group, and
		$$
		{\cal D}\lhd_h {\cal M}_1.
		$$
	\end{itemize}
	
	\item[(iii)] Consider
	$T_n$, $n\in \Z_+$ as the set of $n\times n$ non-singular matrices, and set
	$$
	T:=\bigcup_{i=1}^{\infty}T_n,
	$$
	with $\ltimes$  defined by (\ref{2.2.1}), and ${\bf e}$ as in (\ref{3.2.3}). Then
	${\cal T}:=(T, \ltimes, {\bf e})$ is a  hyper group.
	\begin{itemize}
		\item Set
		$$O_{n}:=\{A\in {\cal T}_{n}\;|\; A^{^{\rm T}}=A^{-1} \},
		$$
		$$
		O=\bigcup_{n=1}^{\infty}O_n.
		$$
		Then ${\cal O}:=({O}_{n},\ltimes,{\bf e})$ is a hyper group, and
		$$
		{\cal O}<_h {\cal T}.
		$$
		\item Set
		$$
		E_{n}:=\{A\in {\cal M}_{n\times n}\;|\; \det(A)=1\},
		$$
		and
		$$
		E=\bigcup_{n=1}^{\infty}E_n.
		$$
		Then ${\bf E}=(E,\ltimes, {\bf e})$ is a hyper group, and
		$$
		{\bf E}<_h {\cal T}.
		$$
	\end{itemize}
\end{itemize}
\end{exa}

\begin{prp}\label{p3.4.4} Let ${\cal G}=(G.*,{\bf e})$ be a hyper group, and $H\subset G$. ${\cal H}=(H.*,{\bf e}_H)$ is a sub-hyper group, if and only if,
\begin{itemize}
	\item[(i)]
	\begin{align}\label{3.4.8}
		{\bf e}_H=\{e_x\;|\;x\in H\}
	\end{align}
	 is a sub-lattice;
	\item[(ii)]
	\begin{align}\label{3.4.9}
		x*y^{-1}\in H,\quad \forall x,y\in H.
	\end{align}
\end{itemize}
\end{prp}

\noindent{\it Proof.} The necessity is obvious. To see the sufficiency, assume $x,y\in H_e$, then $x*y^{-1}\in H_e$. That is,
$H_e<G_e$ is a sub-group. Hence, $H$ is a sub-hyper group of $G$.
\hfill $\Box$

\begin{rem}\label{r3.4.5} The hyper group is not a group, because it has multiple identities. Roughly speaking,  it is a cluster formed by a class of groups of the same type.
The operator of the hyper group covers all operators of  groups in the cluster as its special case, and all the identities as its identity set.
\end{rem}

\section{Permutation Hyper Group}

\subsection{Left Permutation Hyper Group}

Permutation group is one of the  most important finite groups, because any finite group can be isomorphic to a subgroup of certain permutation group \cite{dix96}. In this section we construct a hyper group structure for  permutations with different orders.

Let ${\bf S}_n$ be the permutation group of order $n$ with the group operator $\circ$ as the compound of two permutations. Set
${\bf S}:=\bigcup_{n=1}^{\infty}{\bf S}_n$.

To build a cross-order product of permutations with different orders, we need to construct a cross-dimensional mapping.

\begin{dfn}\label{d4.1.1} Assume $m|n$, and $km=n$, then we define a mapping, called the left embedding, as $\varphi^m_n:{\bf S}_m\ra {\bf S}_n$ as follows:
	Assume $\sigma\in {\bf S}_m$, then
	\begin{align}\label{4.1.1}%{2.6}
			\varphi^m_n(\sigma)((i-1)k+s):=(\sigma(i)-1)k+s,
			\quad \sigma\in{\bf S}_m, \; s\in [1,k],\; i\in [1,m].
	\end{align}
\end{dfn}

Note that when $m=n$ $\varphi^m_m$ is the identity mapping. When $m<n$, to understand this mapping, we recall permutation matrices. Recall the definition of permutation matrices in Example \ref{e2.1.2}. A permutation matrix $P\in {\cal M}_{n\times n}$ is a non-singular logical matrix. A logical matrix $A\in {\cal L}_{n\times n}$ satisfies
$\Col(A)\subset {\cal D}_n$, hence it has the form that $A=[\d_n^{i_1},\d_n^{i_2},\cdots,\d_n^{i_n}]$ and can be briefly denoted by
$A=\d_n[i_1,i_2,\cdots,i_n]$.

Recall that the set of $n\times n$ permutation matrices form a subgroup
${\cal P}_n<{\cal T}_n$.

For each $\sigma\in {\bf S}_n$, its permutation matrix is expressed as
\begin{align}\label{4.1.2}
	P_{\sigma}:=\d_n[\sigma(1),\sigma(2),\cdots,\sigma(n)].
\end{align}

According to the construction of a  permutation matrix, we have the following result.

\begin{prp}\label{p4.1.2}
	\begin{itemize}
		\item[(i)] Assume a permutation $\sigma\in {\bf S}_n$ and its permutation matrix is $P_{\sigma}$. Then
		\begin{align}\label{4.1.3}
			P_{\sigma}\d_n^i =\d_n^{\sigma(i)},\quad i\in [1,n].
		\end{align}
		\item[(ii)] Let $\sigma,~\mu\in {\bf S}_n$. Then
		\begin{align}\label{4.1.4}
			M_{\sigma\circ \mu}=M_{\sigma}M_{\mu}.
		\end{align}
		\item[(iii)] Assume $m|n$, set $km=n$, and let $\sigma\in {\bf S}_m$. Then
		\begin{align}\label{4.1.4}
			P_{\varphi(\sigma)}=P_{\sigma}\otimes I_k.
		\end{align}
	\end{itemize}
\end{prp}

\begin{rem}\label{r4.1.201}
	\begin{itemize}
		\item[(i)] Define $\pi:{\bf S}_n\ra {\cal P}_n$ by $\pi:\sigma\mapsto P_{\sigma}$, which is as shown in (\ref{4.1.2}). Then Proposition \ref{p4.1.2} shows that $\pi$ is a group isomorphism.
		\item[(ii)] Assume $m|n$, the $\varphi^m_n$ is an one-to-one homomorphism from ${\bf S}_m$ to ${\bf S}_n$. In other words,
		$\varphi^m_n$ is an isomorphism from ${\bf S}_m$ to its image $\varphi({\bf S}_m)< {\bf S}_n$. It is an embedding.
	\end{itemize}
\end{rem}

Now we are ready to define a cross-order product of permutations as follows.

\begin{dfn}\label{d4.1.3} Let $\sigma,~\mu\in {\bf S}$, say, $\sigma\in {\bf S}_m$, $\mu\in {\bf S}_n$, and $t=\lcm(m,n)$. Then the product, called the left product, of  $\sigma$ and $\mu$, denoted by $\sigma\odot \mu$, is defined by
	\begin{align}\label{4.1.5}
		\sigma\odot \mu=\varphi^m_t(\sigma)\circ \varphi^n_t(\mu)\in {\bf S}_t.
	\end{align}
\end{dfn}

The following result can be verified by the construction easily.

\begin{prp}\label{p4.1.3} Let $\sigma,~\mu \in {\bf S}$ with their permutation matrices $M_{
		\sigma}$ and $M_{\mu}$ respectively. Then
	\begin{align}\label{4.1.6}
		M_{\sigma\odot \mu}=M_{\sigma}\ltimes M_{\mu}.
	\end{align}
\end{prp}
$\varphi^m_n$ is called the left embedding because it is determined by this left STP of their permutation matrices.

Define the set of cross-order permutation matrices as
$$
{\cal P}:=\bigcup_{n=1}^{\infty}{\cal P}_n.
$$
Then it is easy to verify the following.

\begin{prp}\label{p4.1.4}
	\begin{itemize}
		\item[(i)] $({\cal P},\ltimes,{\bf e})$ is a hyper group, where
		${\bf e}$ is defined by (\ref{3.2.3}).
		\item[(ii)]
		\begin{align}\label{4.1.7}
			({\cal P},\ltimes,{\bf e})<_h {\cal T},
		\end{align}
		where ${\cal T} =(T,\ltimes,{\bf e})$ is defined in Proposition \ref{p3.2.3}.
	\end{itemize}
\end{prp}

Using Propositions \ref{p4.1.3} and \ref{p4.1.4}, the following result is obvious.

\begin{thm}\label{t4.1.8}
	${\bf S}^L:=({\bf S},\odot, {\bf e})$ is a hyper group, where the identity set is
	\begin{align}\label{4.1.8}
		{\bf e}=\{e_n={\rm Id}_n\in {\bf S}_n\;|\; n\in \Z_+\},
	\end{align}
	with respect to MD-1 lattice $(\Z_+,\prec)$.
\end{thm}

Next, we consider some basic properties of hyper group ${\bf S}^L$.
%
%The first result as follows is obvious.
%
%\begin{prp}\label{p4.1.9} The permutation hyper group is isomorphic to the hyper group of permutation matrices ${\cal P}$.
%That is,
%\begin{align}\label{4.1.9}
%{\bf S}\cong {\cal P}.
%\end{align}
%\end{prp}

\begin{prp}\label{p4.1.10}  Set the alternative set of permutations by
	\begin{align}\label{4.1.10}
		{\bf A}:=\bigcup_{n=1}^{\infty}{\bf A}_n.
	\end{align}
	\begin{itemize}
		\item[(i)] ${\bf A}^L$ is a sub-hyper group of ${\bf S}^L$. Precisely speaking,
		\begin{align}\label{4.1.11}
			({\bf A},\odot, {\bf e}) <_h ({\bf S},\odot, {\bf e}) .
		\end{align}
		\item[(ii)] ${\bf A}^L$ is a normal sub-hyper group of ${\bf S}^L$. That is,
		\begin{align}\label{4.1.12}
			({\bf A},\odot, {\bf e}) \lhd_h ({\bf S},\odot, {\bf e}) .
		\end{align}
		
	\end{itemize}
\end{prp}

\begin{rem}\label{r4.1.11} The definition of normal hyper subgroup can not be defined in the same way as for normal subgroup. It is defined in each component groups. For instance,  note that ${\bf A}
	\lhd_h {\bf S}$ does not mean
	$$
	\sigma\circ {\bf A}={\bf A}\circ \sigma,\forall \sigma\in {\bf S}.
	$$
	It means
	$$
	\sigma\circ {\bf A}_n={\bf A}_n\circ \sigma,\forall \sigma\in {\bf S}_n,\forall n\in \Z_+.
	$$
	In fact, ${\bf A}$ and $g\circ {\bf A}$, where $g\in {\bf S}\backslash {\bf A}$, is not a partition of
	${\bf S}$.
\end{rem}

Finally, we consider the consistence of the equivalence with respect to $\odot$.

\begin{itemize}
	\item[(i)]
	Assume $\sigma_1\sim \sigma_2$ and $\mu_1\sim \mu_2$,
	we verify if
	\begin{align}\label{4.1.13}
		\sigma_1\odot \mu_1 \sim \sigma_2\odot \mu_2.
	\end{align}
	Using the isomorphism between permutations and their structure matrices, (\ref{4.1.13}) is equivalent to
	\begin{align}\label{4.1.14}
		M_{\sigma_1}M_{\mu_1}\sim M_{\sigma_2}M_{\mu_2}.
	\end{align}
	Assume $\sigma_1\in {\bf S}_m$,  $\sigma_2\in {\bf S}_n$, $\mu_1\in {\bf S}_p$,  $\mu_2\in {\bf S}_q$, $m\wedge n=\a$, $p\wedge q=\b$, $m\vee p=\xi$, $n\vee q=\eta$, and $\a\vee \b=\lambda$. According to Lemma \ref{l3.3.6}, there exist $M_{\a}\in {\cal P}_{\a}$ and $M_{\b}\in {\cal P}_{\b}$ such that
	$$
	\begin{array}{cc}
		M_{\sigma_1}=M_{\a}\otimes I_{m/\a},&M_{\sigma_2}=M_{\a}\otimes I_{n/\a},\\
		M_{\mu_1}=M_{\b}\otimes I_{p/\b},&M_{\mu_2}=M_{\b}\otimes I_{q/\b}.\\
	\end{array}
	$$
	Then
	$$
	\begin{array}{l}
		M_{\sigma_1}M_{\mu_1}=(M_{\a}\otimes I_{m/\a}\otimes I_{\xi/m})
		(M_{\b}\otimes I_{p/\b}\otimes I_{\xi/p})\\
		~=(M_{\a}\otimes I_{xi/\a})
		(M_{\b}\otimes I_{\xi/\b})\\
		~=\left[(M_{\a}\otimes I_{\lambda/\a})
		(M_{\b}\otimes I_{\lambda/\b})\right]\otimes I_{\xi/\lambda}.\\
	\end{array}
	$$
	Similarly, we have
	$$
	\begin{array}{l}
		M_{\sigma_1}M_{\mu_1}\\
		~=\left[(M_{\a}\otimes I_{\lambda/\a}) (M_{\b}\otimes I_{\lambda/\b})\right]\otimes I_{\eta/\lambda}.\\
	\end{array}
	$$
	(\ref{4.1.14}) follows.
	\item[(ii)] Verify the (ii) of Definition \ref{d3.3.3}. Assume $x\sim I_r$. Then there exist $I_{\a}$ and $I_{\b}$ such that
	$$
	e\otimes I_{\a}=I_r\otimes I_{\b}.
	$$
	Then it is clear that $e=I_s\in {\bf e}$
\end{itemize}
\hfill $\Box$

Hence we have the following result.

\begin{prp}\label{p4.1.12} There exists the equivalence hyper group
	\begin{align}\label{4.1.15}
		\overline{\bf S}^L:={\bf S}^L/\sim.
	\end{align}
\end{prp}

\begin{rem}\label{r4.1.13} It is interesting to find  what are the elements in $\overline{\bf S}^L$. Let $\sigma\in {\bf S}$. $\sigma$ is said to be reducible, if $\sigma\in {\bf S}_n$ and there are $m<n$, $m|n$ and a $\sigma_0\in {\bf S}_m$ such that
	$$
	\sigma=\varphi^m_n(\sigma_0),
	$$
otherwise, it is irreducible.
	Let $\sigma \in {\bf S}_m$ be irreducible and define
	$$
	\overline{\sigma}=\{\varphi^m_{km}(\sigma)\;|\;k\in \Z_+\}.
	$$
	Then it is easy to verify that
	\begin{align}\label{4.1.16}
		\overline{\bf S}^L:=\{\overline{\sigma}\;|\;\sigma\in {\bf S} ~\mbox{is irreducible}\}.
	\end{align}
\end{rem}

\subsection{Right Permutation Hyper Group}

Alternatively, we can also define a right permutation hyper group.
We need to construct another cross-dimensional mapping.

\begin{dfn}\label{d4.e.1} Assume $m|n$, and $km=n$, then we define a mapping, called the right embedding, as $\psi^m_n:{\bf S}_m\ra {\bf S}_n$ as follows:
	Assume $\sigma\in {\bf S}_m$, then
	\begin{align}\label{4.2.1}%{2.6}
		\begin{array}{l}
			\psi^m_n(\sigma)((s-1)m+i):=(k-s)m +\sigma(i),\\
			\quad \sigma\in{\bf S}_m, \; s\in [1,k],\; i\in [1,m].\\
		\end{array}
	\end{align}
\end{dfn}

Similarly to left embedding case, we have the following result.

\begin{prp}\label{p4.2.2}
	Assume $m|n$, $km=n$, and $\sigma\in {\bf S}_m$. Then
	\begin{align}\label{4.2.4}
		P_{\psi(\sigma)}=I_k\otimes P_{\sigma}.
	\end{align}
\end{prp}

\begin{rem}\label{r4.2.201}
	Assume $m|n$, similarly to $\varphi^m_n$, the $\psi^m_n$ is an one-to-one homomorphism from ${\bf S}_m$ to ${\bf S}_n$. In other words,
	$\psi^m_n$ is an isomorphism from ${\bf S}_m$ to its image $\psi({\bf S}_m)< {\bf S}_n$. It is an  embedding.
\end{rem}

Next, we define another cross-order product of permutations as follows.

\begin{dfn}\label{d4.2.3} Let $\sigma,~\mu\in {\bf S}$, say, $\sigma\in {\bf S}_m$, $\mu\in {\bf S}_n$, and $t=\lcm(m,n)$. Then the right product of  $\sigma$ and $\mu$, denoted by $\sigma\odot_r \mu$, is defined by
	\begin{align}\label{4.2.5}
		\sigma\odot_r \mu=\psi^m_t(\sigma)\circ \psi^n_t(\mu)\in {\bf S}_t.
	\end{align}
\end{dfn}

The following result is obvious.

\begin{prp}\label{p4.2.3} Let $\sigma,~\mu \in {\bf S}$ with their permutation matrices $M_{
		\sigma}$ and $M_{\mu}$ respectively. Then
	\begin{align}\label{4.2.6}
		M_{\sigma\odot_r \mu}=M_{\sigma}\rtimes M_{\mu}.
	\end{align}
\end{prp}
$\psi^m_n$ is called the right embedding because it is determined by the right STP of their permutation matrices.

Using Propositions \ref{p4.2.3} and \ref{p4.1.4}, the following result is obvious.

\begin{thm}\label{t4.2.8}
	${\bf S}^R=({\bf S},\odot_r, {\bf e})$ is a hyper group, where the identity set is
	\begin{align}\label{4.2.8}
		{\bf e}=\{e_n={\rm Id}_n\in {\bf S}_n\;|\; n\in \Z_+\},
	\end{align}
	with respect to MD-1 lattice $(\Z_+,\prec)$.
\end{thm}

Correspondingly, the right hyper group ${\bf S}^R$ has the following properties.

%\begin{prp}\label{p4.2.9} The right permutation hyper group is isomorphic to the hyper group of permutation matrices ${\cal P}$.
%That is,
%\begin{align}\label{4.2.9}
%{\bf S}\cong {\cal P}.
%\end{align}
%\end{prp}

\begin{prp}\label{p4.2.10}
	${\bf A}^R$ is a normal sub-hyper group of ${\bf S}^R$. That is,
		\begin{align}\label{4.2.12}
			({\bf A},\odot_r, {\bf e}) \lhd_h ({\bf S},\odot_r, {\bf e}) .
		\end{align}
		
\end{prp}

$\sigma$ and $\mu$ are said to be right-equivalent, denoted by $\sigma\sim_r \mu$, if
there exsit $I_r$ and $I_s$ such that
\begin{align}\label{4.2.2.14}
	I_r\otimes M_{\sigma}=I_s\otimes M_{\mu}.
\end{align}

Then we can also prove that the equivalence is consistent with respect to the $\odot_r$, hence we have
\begin{prp}\label{p4.2.12} There exists the equivalence hyper group
	\begin{align}\label{4.2.15}
		\overline{\bf S}^R:={\bf S}^R/\sim_r.
	\end{align}
\end{prp}

Denote by
$$
\tilde{\sigma}=\{\mu\;|\;\mu\sim_r \sigma\},
$$
then we have
\begin{prp}\label{p4.2.13}
	\begin{align}\label{4.2.16}
		\overline{\bf S}^R:=\{\tilde{\sigma}\;|\;\sigma\in {\bf S} ~\mbox{is irreducible}\}.
	\end{align}
	
\end{prp}

\begin{rem}\label{r4.2.14} My conjecture is: the left permutation hyper group and the right permutation hyper group are isomorphic, i.e.,
$$
{\bf S}^L\cong_h {\bf S}^R.
$$
But it is open so far.
\end{rem}

\begin{exa}\label{e4.2.15}
\begin{itemize}
\item[(i)]
Consider ${\bf S}_2\subset {\bf S}^L$ and ${\bf S}_3\subset {\bf S}^L$.
The product $\odot$ is defined by
merging ($\varphi$) each element in ${\bf S}_2$ or ${\bf S}_3$ into ${\bf S}_6$ . Then
$$
a\odot b:=\varphi(a)\circ \varphi(b),
$$
where $\circ$ is the product of ${\bf S}_6$.

Let $e^2_1={\rm Id}_2$ (the identity),  $e^2_2=(1,2)$ .
We  use $e^2_2$ to demonstrate the merging mapping ($\varphi$).
The permutation matrix of $e^2_2$ is
$$
M_{e^2_2}=\begin{bmatrix} 0&1//1&0\end{bmatrix}=\d_2[2,1].
$$
Then
$$
M_{\varphi(e^2_2)}=M_{e^2_2}\otimes I_3=\d_6[4,5,6,1,2,3].
$$
That is
$$
\varphi(e^2_2)=(1,4)(2,5)(3,6)\in {\bf S}_6.
$$

Similarly,
$$
M_{\varphi(e^3_i)}=M_{e^3_i}\otimes I_2,\quad i\in[1,6].
$$

Finally, we have
$$
\varphi(e^2_1)={\rm Id}_6,\quad \varphi(e^2_2)=(1,4)(2,5)(3,6).
$$
Denote by $e^3_1={\rm Id}_3$, $e^3_2=(1,2)$, $e^3_3=(1,3)$, $e^3_4=(2,3)$, $e^3_5=(1,2,3)$, and $e^3_6=(1,3,2)$. Then
it is easy to calculate that
$$
\begin{array}{lll}
\varphi(e^3_1)={\rm Id}_6,&\varphi(e^3_2)=(1,3)(2,4)&\varphi(e^3_3)=(1,5)(2,6)\\
\varphi(e^3_4)=(3,5)(4,6),&\varphi(e^3_5)=(1,3,5)(2,4,6)&\varphi(e^3_6)=(1,5,3)(2,6,4).\\
\end{array}
$$
Finally, it is easy to see that the sub-hyper group generated by ${\bf S}_2$ and ${\bf S}_2$ is
$$
({\bf S}_2\bigcup {\bf S}_3\bigcup {\bf S}_6,\odot)<{\bf S}^L.
$$
which has identity set $\{{\rm Id}_2,~{\rm Id}_3,~{\rm Id}_6\}$.

\item[(ii)]
Consider ${\bf S}_2\subset {\bf S}^R$ and ${\bf S}_3\subset {\bf S}^R$.
The product $\odot_r$ is defined by
merging ($\varphi$) each element in ${\bf S}_2$ or ${\bf S}_3$ into ${\bf S}_6$ . Then
$$
a\odot_r b:=\varphi(a)\circ \varphi(b).
$$
Now the merging is defined by their corresponding permutation matrices as, for example,
$$
M_{\varphi(e^3_i)}=I_2\otimes M_{e^3_i},\quad i\in[1,6].
$$
Similar process leads to
$$
\varphi(e^2_1)=id_6,\quad \varphi(e^2_2)=(1,2)(3,4)(5,6).
$$
In addition,
$$
\begin{array}{lll}
\varphi(e^3_1)=id_6,&\varphi(e^3_2)=(1,2)(4,5)&\varphi(e^3_3)=(1,3)(4,6)\\
\varphi(e^3_4)=(2,3)(5,6),&\varphi(e^3_5)=(1,2,3)(4,5,6)&\varphi(e^3_6)=(1,3,2)(4,6,5).\\
\end{array}
$$
We also have  the sub-hyper group generated by ${\bf S}_2$ and ${\bf S}_2$ is
$$
({\bf S}_2\bigcup {\bf S}_3\bigcup {\bf S}_6,\odot_r)<{\bf S}^R.
$$

\end{itemize}
\end{exa}

\section{Homomorphism}

\subsection{Hyper Group Homomorphism/Isomorphism}

\begin{dfn}\label{d5.1.1} Let ${\cal G}_i=(G^i,\times^i,{\bf e}^i)$, $i=1,2$ be two hyper groups. $\pi:G^1\ra G^2$ is called a hyper group homomorphism, and ${\cal G}_1$ and ${\cal G}_2$ are said to be homomorphic, denoted by  ${\cal G}_1 \simeq_h {\cal G}_2$, if the followings are satisfied:
	\begin{itemize}
		\item[(i)]
		\begin{align}\label{5.1.1}
			\pi(x\times^1 y)=\pi(x)\times^2\pi(y),\quad x,y\in G_1.
		\end{align}
		\item[(ii)]
		\begin{align}\label{5.1.2}
			\pi(e^1)\in  {\bf e}^2, \quad e^1\in {\bf e}^1,
		\end{align}
		and $\pi: {\bf e}^1\ra {\bf e}^2$ is a lattice order-preserving mapping.
		\item[(iii)] If a hyper group homomorphism $\pi:G_1\ra G_2$ is bijective, and $\pi^{-1}:G_2\ra G_1$ is also a homomorphism, then
		$\pi$ is called a  hyper group isomorphism, and ${\cal G}_1$ and ${\cal G}_2$ are said to be isomorphic, denoted by ${\cal G}_1\cong_h {\cal G}_2$.
	\end{itemize}
\end{dfn}

\begin{exa}\label{e5.1.2}
	Recall ${\cal T}_H$ defined in Example \ref{e3.4.3}.
	\begin{itemize}
		\item[(i)] Let ${\rm Id}: {\cal T}_H \hookrightarrow {\cal T}$ be the including mapping.  Then ${\rm Id}$ is a hyper group homomorphism, and
		${\cal T}_H\simeq_h {\cal T}$.
		\item[(ii)] Let $\pi_k:{\cal T}_H \ra {\cal T}$ be defined by $A\mapsto A\otimes I_k$.
		Then $(\pi_k({\cal T}_H),\ltimes) \simeq_h {\cal T}$.
		\item[(iii)]  Since ${\rm Id}$ is one-to-one. Then it is easy to verify that
		$$
		{\cal T}_H \cong_h (\pi_k({\cal T}_H),\ltimes).
		$$
	\end{itemize}
\end{exa}

\begin{exa}\label{e5.1.3} Consider $({\cal M}_{\mu}, \bar{+})$  (or $({\cal M}_{\mu}, \hat{+})$).
	Let $\pi: A\mapsto A^{\rm T}$. Then $\pi$ is a hyper group isomorphism, i.e.,
	\begin{align}\label{5.1.3}
		\begin{array}{l}
			({\cal M}_{\mu}, \bar{+}) \cong_h  ({\cal M}_{1/\mu}, \bar{+});\\
			(\mbox{or}~ ({\cal M}_{\mu}, \hat{+})\cong_h  ({\cal M}_{1/\mu}, \hat{+})).\\
		\end{array}
	\end{align}
\end{exa}

\begin{exa}\label{e5.1.4} Consider  $({\bf S}^L, \odot)$ and $({\cal P},\ltimes)$ (or $({\bf S}^R, \odot_r)$ and $({\cal P},\rtimes)$).
	Then $\pi:\sigma\mapsto M_{\sigma}$ is a hyper group isomorphism, i.e.,
	\begin{align}\label{5.1.3}
		\begin{array}{l}
			({\bf S}^L, \odot) \cong_h ({\cal P},\ltimes);\\
			(\mbox{or}~ ({\bf S}^R, \odot_r) \cong_h  ({\cal P},\rtimes)).\\
		\end{array}
	\end{align}
\end{exa}

\subsection{Quotient Hyper Group Homomorphism/Isomorphism}

Assume $\varphi: H\ra G$ is a hyper group homomorphism. Set
\begin{align}\label{5.2.1}
	\begin{array}{l}
		\Im(\varphi):=\{g\in G\;|\; \exists h\in H~\mbox{such that}~g=\varphi(h)\},\\
		{\rm Ker}(\varphi):=\{h\in H\;|\; \varphi(h)\in {\bf e}_G\}.
	\end{array}
\end{align}
Then we have the following result.
\begin{prp}\label{p5.2.1}
	\begin{itemize}
		\item[(i)] ${\rm Im}(\varphi)<_h G$.
		\item[(ii)] ${\rm Ker}(\varphi) \lhd_h H$.
	\end{itemize}
\end{prp}

\noindent{\it Proof.}
\begin{itemize}
	\item[(i)]
	
	Since $\varphi:{\bf e}_H\ra {\bf e}_G$ is a lattice homomorphism, it is easy to see that $\varphi({\bf e}_H)$ is a sub-lattice of ${\bf e}_G$.  Let $x,y\in {\rm Im}(\varphi)$. Then there exist $a,b\in H$ such that $\varphi(a)=x$, $\varphi(b)=y$. Since $H$ is a hyper group, $y^{-1}\in H$. Then $\varphi(ab^{-1})=xy^{-1}\in {\rm Im}(\varphi)$. We conclude that ${\rm Im}(\varphi)<_h G$.
	
	\item[(ii)] Since ${\bf e}_H\subset  K:={\rm Ker}(\varphi)$, ${\bf e}_{K}={\bf e}_H$. It is obvious that $K <_h H$. To see $K$ is a normal hyper-sub-group, we need to show that for any $a\in H$
	$aKa^{-1}=K$. This is true because
	$$
	\varphi(aka^{-1}) =\varphi(a)e(\varphi(a))^{-1}=e\varphi(a)(\varphi(a))^{-1}\in {\bf e}_G.
	$$
\end{itemize}
\hfill $\Box$

\begin{prp}\label{p5.2.2} Let $\varphi:H\ra G$ be a hyper group homomorphism, and for each $e\in {\bf e}_H$, $\varphi_e:H_e\ra G$ be the restriction of $\varphi$ on $H_e$.
	Then $\varphi_e:H_e\ra G_{\varphi{e}}$ is a group homomorphism. Moreover, for each $e\in {\bf e}_H$
	\begin{itemize}
		\item[(i)] ${\rm Im}(\varphi_e)< G_e$;
		\item[(ii)] ${\rm Ker}(\varphi_e) \lhd H$.
	\end{itemize}
\end{prp}

\noindent{\it Proof.} We claim that $\varphi_e(H_e)\subset G_{\varphi(e)}$. If this claim is true,  all the remaining statements are obvious. Let $h\in H_e$, then
$$
\varphi(hh^{-1})=\varphi(h)\varphi(h^{-1})=\varphi(e)\in {\bf e}_G.
$$
Hence $\varphi(h)\in G_{\varphi(e)}$, which proves the claim.
\hfill $\Box$

The three Isomorphism Theorems for groups can be modified to suit the hyper groups.

\begin{thm}\label{t5.2.2} (First Isomorphism Theorem)
	\begin{itemize}
		\item[(i)] Assume $\varphi:H\ra G$ is a hyper group homomorphism. For each $e\in {\bf e}_H$, set $K_e={\rm Ker}(\varphi_e)$ and ${\rm Im}_e={\rm Im}(\varphi_e)$. Define a mapping
		$\theta_e:H_e/K_e\ra {\rm Im}_e$ as $h_eK_e\mapsto \varphi_e(h_e)$, $h_e\in H_e$. Then $\theta_e$ is an isomorphism. That is.
		\begin{align}\label{5.2.2}
			H_e/K_e\cong {\rm Im}_e,\quad e\in {\bf e}_H.
		\end{align}
		\item[(ii)] Let $H\lhd_h G$ be a normal hyper-sub-group. For each $e\in {\bf e}_H$ define $\pi_e:G_{\varphi(e)}\ra G_{\varphi(e)}/H_e$ by $g_e\mapsto g_eH_e$. Then $\pi_e$ is a surjective
		homomorphism with ${\rm Ker}(\pi_e)=H_e$.
	\end{itemize}
\end{thm}

\begin{thm}\label{t5.2.3} (Second Isomorphism Theorem)
	Assume $H,~N,~G$ are hyper groups.
	\begin{itemize}
		\item[(i)] If $H<_h G$ and $N\lhd_h G$,   then for each $e\in {\bf e}_H\cap {\bf e}_N$
		\begin{align}\label{5.2.3}
			N_e\bigcap H_e \lhd_h H_e.
		\end{align}
		\item[(ii)] Set $NH=\{nh\;|\; n\in N,~h\in H\}$, where $N,H$ are as in (i). Then for each $e\in {\bf e}_H\bigcap {\bf e}_N$
		\begin{align}\label{5.2.4}
			N_eH_e <_h G_e.
		\end{align}
		\item[(iii)] Define $\pi_e:H_e/(N_e\cap H_e)\ra N_eH_e/N_e$ by
		$$
		h_e(N_e\bigcap H_e)\mapsto h_eN_e,\quad h_e\in H_e,
		$$
		where $e=e_he_h$.
		Then $\pi_e$ is an isomorphism, i.e.,
		\begin{align}\label{5.2.5}
			H_e/(H_e\bigcap N_e)\cong N_eH_e/N_e.
		\end{align}
	\end{itemize}
\end{thm}

\begin{thm}\label{t5.2.4} (Third Isomorphism Theorem)
	Let $M,N,G$ be hyper groups, $M\lhd_h G$, $N\lhd_h G$, and $N<_h M$. Then for eacg $e\in {\bf e}_N$
	\begin{itemize}
		\item[(i)]
		\begin{align}\label{5.2.6}
			M_e/N_e\lhd_h G_e/N_e.
		\end{align}
		\item[(ii)]
		\begin{align}\label{5.2.7}
			(G_e/N_e)/(M_e/N_e)\cong G_e/M_e.
		\end{align}
	\end{itemize}
\end{thm}

\begin{exa}\label{e5.2.5} Consider $\Lambda=\{(2,2),(2,4),(4,2),(4,4)\}$ as a sub-lattice of $\Z_+\times \Z_+$,
	$H=\bigcup_{\lambda\in \Lambda}{\cal M}_{\lambda}$, and
	$$
	{\cal H}:=(H,\hat{+})<_h {\cal G}=({\cal M},\hat{+}).
	$$
	Define a hyper group homomorphism: $\varphi:{\cal H}\ra {\cal G}$ as
	$$
	\varphi(A)=A\hat{-} A^{^{\rm T}}.
	$$
	A straightforward computation shows that
	$$
	K:={\rm Ker}(\varphi)=K_{4\times 4}\bigcup K_{4\times 2}\bigcup K_{2\times 4}\bigcup K_{2\times 2},
	$$
	where
	$$
	\begin{array}{l}
		K_{4\times 4}=\{A\in {\cal M}_{4\times 4}\;|\;A^{^{\rm T}}=A\};\\
		K_{4\times 2}=\{A\otimes \J_2\;|\; A\in {\cal M}_{2\times 2}, ~\mbox{and}~A^{^{\rm T}}=A\};\\
		K_{2\times 4}=\{A\otimes \J^{^{\rm T}}_2\;|\; A\in {\cal M}_{2\times 2}, ~\mbox{and}~A^{^{\rm T}}=A\};\\
		K_{2\times 2}=\{A\in {\cal M}_{2\times 2}\;|\;A^{^{\rm T}}=A\}.\\
	\end{array}
	$$
	$$
	{\rm Im}:={\rm Im}(\varphi)={\rm Im}_{4\times 4}\bigcup {\rm Im}(2\times 2),
	$$
	where
	$$
	\begin{array}{l}
		{\rm Im}_{4\times 4}=\{A\in {\cal M}_{4\times 4}\;|\;A^{^{\rm T}}=-A\};\\
		{\rm Im}_{2\times 2}=\{A\in {\cal M}_{2\times 2}\;|\;A^{^{\rm T}}=-A\}.\\
	\end{array}
	$$
	
	Since
	$$
	K_{\a\times \b}\hat{+}K_{\gamma\times \d}=
	\begin{cases}
		K_{4\times 4},\quad \max(\a,\b,\gamma,\d)=4,\\
		K_{2\times 2},\quad \max(\a,\b,\gamma,\d)=2,\\
	\end{cases}
	$$
	$K\lhd_h {\cal H}$ is a hyper-sub-group. It is normal because ${\cal H}$ is abelian.
	Since
	$$
	{\rm Im}_{\a\times \a}\hat{+}{\rm Im}_{\b\times \b}=
	\begin{cases}
		{\rm Im}_{4\times 4},\quad \max(\a,\b)=4,\\
		{\rm Im}_{2\times 2},\quad \max(\a,\b)=2,\\
	\end{cases}
	$$
	${\rm Im}\lhd_h {\cal G}$ is a hyper-sub-group.  Note that the lattice order-preserving mapping is:
	$$
	\pi(e_{\lambda})=
	\begin{cases}
		(4,4),\quad \lambda\in\{(2,4),(4,2),(4,4)\},\\
		(2,2),\quad \lambda=(2,2).
	\end{cases}
	$$
	
	Next, we verify First Isomorphism Theorem.  Consider each case when $\lambda\in {\bf e}_H$.
	\begin{itemize}
		\item[(i)] $\lambda=(4,4)$.
		$$
		K_{\lambda}=\{A\in {\cal M}_{4\times 4}\;|\;A^{^{\rm T}}=A\}.
		$$
		$$
		{\rm Im}_{\lambda}=\{A\in {\cal M}_{4\times 4}\;|\;A^{^{\rm T}}=-A\}.
		$$
		Let $A\in H_{\lambda}={\cal M}_{4\times 4}$. Since
		$$
		A=\frac{1}{2}(A+A^{^{\rm T}})+\frac{1}{2}(A-A^{^{\rm T}}),
		$$
		$$
		AK_{\lambda}=\frac{1}{2}(A+A^{^{\rm T}})+ \frac{1}{2}(A-A^{^{\rm T}})+K_{\lambda} = \frac{1}{2}(A-A^{^{\rm T}})+K_{\lambda}.
		$$
		Define $\pi:H_{\lambda}/K_{\lambda}\ra {\rm Im}_{\lambda}$ as
		$$
		\pi(AK_{\lambda}):= \frac{1}{2}(A-A^{^{\rm T}}).
		$$
		Then $\pi$ is an isomorphism, i.e.,
		$$
		H_{\lambda}/K_{\lambda}\cong {\rm Im}_{\lambda}.
		$$
		
		\item[(ii)] $\lambda=(4,2)$.
		A straightforward calculation shows that
		$$
		K_{\lambda}=\left\{\left.
		\begin{bmatrix}
			a&a&b&b\\
			b&b&c&c
		\end{bmatrix}\;\right|\;a,b,c\in \R
		\right\}.
		$$
		Set
		$$
		\begin{array}{l}
			A=\begin{bmatrix}
				a_{1,1}&a_{1,2}&a_{1,3}&a_{1,4}\\
				a_{2,1}&a_{2,2}&a_{2,3}&a_{2,4}\\
			\end{bmatrix}\\
			~=\begin{bmatrix}
				a_{1,1}&a_{1,1}&a_{1,3}&a_{1,3}\\
				a_{1,3}&a_{1,3}&a_{2,3}&a_{2,3}\\
			\end{bmatrix}\\
			~+
			\begin{bmatrix}
				0&a_{1,2}-a_{1,1}&0&a_{1,4}-a_{1,3}\\
				a_{2,1}-a_{1,3}&a_{2,2}-a_{1,3}&0&a_{2,4}-a_{1,3}\\
			\end{bmatrix}\\
			~:=\begin{bmatrix}
				a_{1,1}&a_{1,1}&a_{1,3}&a_{1,3}\\
				a_{1,3}&a_{1,3}&a_{2,3}&a_{2,3}\\
			\end{bmatrix}+
			\begin{bmatrix}
				0&d&0&e\\
				f&g&0&h\\
			\end{bmatrix}\\
		\end{array}
		$$
		Hence
		$$
		AK_{\lambda}=
		\begin{bmatrix}
			0&d&0&e\\
			f&g&0&h\\
		\end{bmatrix}+K_{\lambda}.
		$$
		The ${\rm Im}_{\lambda}$ is calculated as
		$$
		\begin{array}{ccl}
			{\rm Im}_{\lambda}&=&\left\{\left.
			\begin{bmatrix}
				0&d&-f&e-f\\
				-d&0&-g&e-g\\
				f&g&0&h\\
				f-e&g-e&-h&0
			\end{bmatrix}\;\right|\right.\\
			~&~&\left.\;d,e,f,g,h\in \R\right\}
			\lhd {\cal K}_4.
		\end{array}
		$$
		That is, ${\rm Im}_{\lambda}$ is a (normal) sub-group of ${\cal K}_4$.
		
		Define $\pi:H_{\lambda}/K_{\lambda}\ra {\rm Im}_{\lambda}$ as
		$$
		%\begin{array}{l}
		\begin{bmatrix}
			0&d&0&e\\
			f&g&0&h\\
		\end{bmatrix}+K_{\lambda}%\\
		\mapsto%\\
		\begin{bmatrix}
			0&d&-f&e-f\\
			-d&0&-g&e-g\\
			f&g&0&h\\
			f-e&g-e&-h&0
		\end{bmatrix}
		%\end{array}
		$$
		Then it is easy to verify that $\pi$ is an isomorphism.That is,
		$$
		H_{\lambda}/K_{\lambda}\cong {\rm Im}_{\lambda}.
		$$
		
		\item[(iii)] $\lambda=(2,4)$.
		Similarly to (ii), we have
		$$
		K_{\lambda}=\left\{\left.
		\begin{bmatrix}
			a&b\\
			a&b\\
			b&c\\
			b&c
		\end{bmatrix}\;\right|\;a,b,c\in \R
		\right\}.
		$$
		$$
		AK_{\lambda}=
		\begin{bmatrix}
			0&f\\
			d&g\\
			0&0\\
			e&h\\
		\end{bmatrix}+K_{\lambda}.
		$$
		$$
		\begin{array}{ccl}
			{\rm Im}_{\lambda}&=&\left\{\left.
			\begin{bmatrix}
				0&-d&f&f-e\\
				d&0&g&g-e\\
				-f&-g&0&-h\\
				e-f&e-g&h&0
			\end{bmatrix}\;\right|\right.\\
			~&~&\left.\;d,e,f,g,h\in \R\right\}
			<{\cal K}_4.
		\end{array}
		$$
		That is, ${\rm Im}_{\lambda}$ is a sub-group of ${\cal K}_4$.
		
		Define $\pi:H_{\lambda}/K_{\lambda}\ra {\rm Im}_{\lambda}$ as
		$$
		%\begin{array}{l}
		\begin{bmatrix}
			0&f\\
			d&g\\
			0&0\\
			e&h\\
		\end{bmatrix}+K_{\lambda}%\\
		\mapsto%\\
		\begin{bmatrix}
			0&-d&f&f-e\\
			d&0&g&g-e\\
			-f&-g&0&-h\\
			e-f&e-g&h&0
		\end{bmatrix}
		%\end{array}
		$$
		Then $\pi$ is an isomorphism, i.e.,
		$$
		H_{\lambda}/K_{\lambda}\cong {\rm Im}_{\lambda}.
		$$
		\item[(iv)] $\lambda=(2,2)$.
		$$
		K_{\lambda}=\{A\in {\cal M}_{2\times 2}\;|\;A^{\rm T}=A\}.
		$$
		$$
		{\rm Im}_{\lambda}=\{A\in {\cal M}_{2\times 2}\;|\;A^{\rm T}=-A\}.
		$$
		Let $A\in H_{\lambda}={\cal M}_{2\times 2}$. Since
		$$
		A=\frac{1}{2}(A+A^{\rm T})+\frac{1}{2}(A-A^{\rm T}),
		$$
		$$
		AK_{\lambda}=\frac{1}{2}(A+A^{\rm T})+ \frac{1}{2}(A-A^{\rm T})+K_{\lambda} = \frac{1}{2}(A-A^{\rm T})+K_{\lambda}.
		$$
		Define $\pi:H_{\lambda}/K_{\lambda}\ra {\rm Im}_{\lambda}$ as
		$$
		\pi(AK_{\lambda}):= \frac{1}{2}(A-A^{\rm T}).
		$$
		Then $\pi$ is an isomorphism, i.e.,
		$$
		H_{\lambda}/K_{\lambda}\cong {\rm Im}_{\lambda}.
		$$
	\end{itemize}
\end{exa}

\subsection{Symmetric and Skew symmetric Matrices}

\begin{dfn}\label{d5.3.1} Consider the hyper group ${\cal G}_M:=({\cal M},\hat{+})$.
	\begin{itemize}
		\item[(i)] A mapping $\pi_s: {\cal G}_M\ra {\cal G}_M$, called the symmetrizing mapping, is defined by
		\begin{align}\label{5.3.1}
			\pi_s: A\mapsto \frac{1}{2}(A\hat{+} A^{\rm }T),\quad A\in {\cal M}.
		\end{align}
		\item[(ii)] A mapping $\pi_a: {\cal G}_M\ra {\cal G}_M$, called the alternating mapping, is defined by
		\begin{align}\label{5.3.2}
			\pi_a: A\mapsto \frac{1}{2}(A\hat{-} A^{\rm T}),\quad A\in {\cal M}.
		\end{align}
	\end{itemize}
\end{dfn}

The following result makes the symmetry and skew-symmetry of matrices meaningful for non-square matrices.

\begin{thm}\label{t5.3.2}  Consider the hyper group ${\cal G}_M:=({\cal M},\hat{+})$.
	\begin{itemize}
		\item[(i)] Both $\pi_s$ and $\pi_a$ are hyper group homomorphisms (precisely speaking, endomorphisms). According to Proposition \ref{p5.2.1}, the kernels of
		$\pi_s$ and $\pi_a$ are hyper groups.
		
		\item[(ii)] The kernel of $\pi_s$, denoted by
$$
K^s:=\left\{A\in {\cal M}\;|\;\pi_s(A)=0\right\},
$$
called the set of skew-symmetric matrices, is $\pi_a$ invariant in the sense that
		\begin{align}\label{5.3.3}
			\pi_a(A)=A\otimes \E_{\a/m\times \a/n},\quad A\in K^s.
		\end{align}

		\item[(iii)] Set
$$
K^s_{m\times n}:=\left\{A\in {\cal M}_{m\times n}\;|\;\pi_s(A)=0\right\},
$$
then $K^s_{n\times n}={\cal K}_n$.
That is, when $m=n$  the kernel set $K^s_{m\times n}$ coincides the classical set of skew-symmetric matrices.

\item[(iv)] The kernel of $\pi_a$, denoted by
$$
K^a:=\left\{A\in {\cal M}\;|\;\pi_a(A)=0\right\},
$$
called the set of symmetric matrices, is $\pi_s$ invariant in the sense that
		\begin{align}\label{5.3.4}
			\pi_s(A)=A\otimes \E_{\a/m\times \a/n},\quad A\in {\cal S}.
		\end{align}

		\item[(v)] Set
$$
K^a_{m\times n}:=\left\{A\in {\cal M}_{m\times n}\;|\;\pi_a(A)=0\right\},
$$
then $K^a_{n\times n}={\cal S}_n$,
i.e., when $m=n$  the kernel set $K^a_{m\times n}$ coincides the classical set of symmetric matrices.
	\end{itemize}
\end{thm}

\noindent{\it Proof.}

\begin{itemize}
	\item[(i)] We show that $\pi_s$ is a hyper group homomorphism.
	Let $A\in {\cal M}_{m\times n}$, and $t=\lcm(m,n)$. Note that
	$$
	A\hat{+} A^{^{\rm T}}=A\otimes \E_{t/m\times t/n}+ A^{^{\rm T}}\otimes \E_{t/n\times t/m}\in {\cal M}_{t\times t}.
	$$
	It is clear that $\pi_s: 0_{m\times n}\ra 0_{m\vee n}$, which is  order-keeping.
	
	Next, we assume $A\in {\cal M}_{m\times n}$, $B\in {\cal M}_{p\times q}$, $s=m\vee p$, $t=n\vee q$,  $\xi=s\vee t$,
	$\a=m\vee n$, and  $\b=p\vee q$.  Then
	
	$$
	\begin{array}{l}
		\pi_s(A\hat{+} B)=\pi_s[ (A\otimes \E_{s/m\times t/n})+(B\otimes \E_{s/p\times t/q})]\\
		~=\frac{1}{2}[(A\otimes \E_{s/m\times t/n})+(B\otimes \E_{s/p\times t/q})]\\
		~\hat{+} \frac{1}{2}[(A^{^{\rm T}}\otimes \E_{t/n\times s/m})+(B^{^{\rm T}}\otimes \E_{t/q\times s/p})]\\
		~=\frac{1}{2}[(A\otimes \E_{s/m\times t/n})+(B\otimes \E_{s/p\times t/q})]\otimes \E_{\xi/s\times \xi/t}\\
		~+ \frac{1}{2}[(A^{^{\rm T}}\otimes \E_{t/n\times s/m})+(B^{^{\rm T}}\otimes \E_{t/q\times s/p})]
		\otimes \E_{\xi/t\times \xi/s}\\
		~=\frac{1}{2}[(A\otimes \E_{\xi/m\times \xi/n})+(A^{^{\rm T}}\otimes \E_{\xi/n\times \xi/m})]\\
		~+\frac{1}{2}[(B\otimes \E_{\xi/p\times \xi/q})+(B^{^{\rm T}}\otimes \E_{\xi/q\times \xi/p})]\\
		~=\frac{1}{2}[(A\otimes \E_{\a/m\times \a/n})+(A^{^{\rm T}}\otimes \E_{\a/n\times \a/m})]\\
		~~\otimes \E_{\xi/\a\times \xi/\a}+\frac{1}{2}[(B\otimes \E_{\b/p\times \b/q})\\
		~+(B^{^{\rm T}}\otimes \E_{\b/q\times \b/p})]\otimes \E_{\xi/\b\times \xi/\b}\\
		~=\frac{1}{2}[(A\otimes \E_{\a/m\times \a/n})+(A^{^{\rm T}}\otimes \E_{\a/n\times \a/m})]\\
		~\hat{+}\frac{1}{2}[(B\otimes \E_{\b/p\times \b/q})+(B^{^{\rm T}}\otimes \E_{\b/q\times \b/p})]\\
		~=\pi_s(A) \hat{+} \pi_s(B).
	\end{array}
	$$
	We conclude that $\pi_s$ is a hyper group homomorphism. Similarly, we can prove  $\pi_a$ is also a
	hyper group homomorphism.
	
	\item[(ii)] Assume $A\in K^s={\rm Ker} (\pi_s)$. Then
	$$
	A\otimes \E_{\a/m\times \a/n}=-A^{^{\rm T}}\otimes \E_{\a/n\times \a/m}.
	$$
	(\ref{5.3.3}) follows immediately.
\end{itemize}

All the rest parts can be proved similarly.
\hfill $\Box$

\begin{cor}\label{c5.3.3} Let $A\in {\cal M}_{m\times n}$ and $t=\lcm(m,n)$. Then there exist unique symmetric matrix $A_s\in {\cal S}_t$ and a skew-symmetric matrix $A_a\in {\cal K}_t$ such that
	\begin{align}\label{5.3.5}
		A\otimes E_{t/m\times t/n}=A_s+A_a,
	\end{align}
	where
	\begin{align}\label{5.3.6}
		A_s=\frac{1}{2}(A\hat{+} A^{^{\rm T}}),\quad
		A_a=\frac{1}{2}(A\hat{-} A^{^{\rm T}}).
	\end{align}
\end{cor}

\subsection{Squaring Matrix}

\begin{dfn}\label{d5.4.1} Let $A\in {\cal M}_{m\times n}\subset {\cal M}$ and $t=\lcm(m,n)$.
	A mapping $\Box:{\cal M}\ra {\cal M}_1$, called a squaring mapping, is defined by
	\begin{align}\label{5.4.1}
		\Box(A):=A\hat{\ttimes} I_t\in {\cal M}_{t\times t}.
	\end{align}
\end{dfn}

\begin{rem}\label{r5.4.2}
	\begin{itemize}
		\item[(i)] It is obvious that if $A$ is a square matrix, then
		\begin{align}\label{5.4.2}
			\Box(A)=A,\quad A\in {\cal M}_1.
		\end{align}
		\item[(ii)]   Let $A\in {\cal M}_{m\times n}$ and $t=\lcm(m,n)$. Then
		\begin{align}\label{5.4.3}
			\Box(A)=A\otimes \E_{t/m\times t/n}.
		\end{align}
	\end{itemize}
\end{rem}

Using squaring mapping,  (\ref{5.3.3})-(\ref{5.3.5}) can be rewritten as  (\ref{5.4.4})-(\ref{5.4.6}) respectively as
\begin{align}\label{5.4.4}
	\pi_a(A)=\Box(A) ,\quad A\in {\cal K}.
\end{align}
\begin{align}\label{5.4.5}
	\pi_s(A)=\Box(A),\quad A\in {\cal S}.
\end{align}
\begin{align}\label{5.4.6}
	\Box(A)=A_s+A_a.
\end{align}

Hence, we have the following result.

\begin{prp}\label{p5.4.3} A non-square matrix $A$ is symmetric (skew-symmetric), if and only if, its squaring matrix $\Box(A)$ is symmetric (skew-symmetric).
\end{prp}

\begin{prp}\label{p5.4.4} $\Box: ({\cal M},\hat{+})\ra ({\cal M}_1,\hat{+})$ is a hyper group homomorphism. Hence,
	\begin{align}\label{5.4.601}
		\Box(A\hat{+} B)=\Box(A) \hat{+} \Box(B),\quad A,B\in {\cal M}.
	\end{align}

\end{prp}

The proof is similar to the proof of Theorem \ref{t5.3.2}.

\begin{dfn}\label{d5.4.5} Let $A\in {\cal M}_{m\times n}$ and $t=\lcm(m.n)$.
	\begin{itemize}
		\item[(i)]
		\begin{align}\label{5.4.7}
			p(\lambda):=\det(\lambda I_t\hat{-}A)
		\end{align}
		is called the s-characteristic function of $A$.
		
		\item[(ii)] $\sigma\in \mathbb{C}$ is called an s-eigenvalue of $A$, if $
		p(\sigma)=0$.
		
		\item[(iii)] $X\in \mathbb{C}^{rm T}$ is called an s-eigenvector of $A$ with respect to $\sigma$, if
		\begin{align}\label{5.4.8}
			\Box(A)X=\sigma X.
		\end{align}
	\end{itemize}
\end{dfn}

We have the s-Cayley-Hamilton theorem as follows.

\begin{thm}\label{t5.4.6} Let $p(\lambda)$ be the s-characteristic function of $A$, then
	\begin{align}\label{5.4.9}
		p(\Box(A))=0.
	\end{align}
\end{thm}

\begin{rem}\label{r5.4.7}
	The squaring mapping poses a non-square matrix all the properties of square matrix. For instance, we can define
	\begin{align}\label{5.4.10}
		\Tr^{\Box}(A)=\Tr(\Box(A)),\quad A\in {\cal M}.
	\end{align}
	Then it is clear that if $\Tr^{\Box}(A)=0$, then $\det(\exp(\Box(A)))=1$.
\end{rem}

\subsection{Homomorphisms Among Component Groups}

\begin{prp}\label{p5.5.1} Let $(G,*)$ be a hyper group with identity set ${\bf e}=\{e_{\lambda}\;|\; \lambda\in \Lambda\}$, where $\Lambda$ is the lattice of identity set. Assume $H<G_{\lambda}$ (or $H\lhd G_{\lambda}$) and $e_{\mu}\in {\bf e}$, then
	\begin{align}\label{5.5.1}
		H*e_{\mu}<G_{\lambda\vee \mu}, (H*e_{\mu}\lhd G_{\lambda\vee \mu}).
	\end{align}
\end{prp}

\begin{exa}\label{e5.5.2}
	\begin{itemize}
		\item[(i)] Since ${\bf A}_n\lhd {\bf S}_n$,
		$$
		{\bf A}\odot {\rm Id}_s \lhd {\bf S}_{n\vee s}.
		$$
		\item[(ii)] Since ${\cal P}_n < {\cal T}_n$.
		$$
		\begin{array}{l}
			{\cal P}_n\ltimes I_s < {\cal T}_{n\vee s};\\
			{\cal P}_n\rtimes I_t < {\cal T}_{n\vee t}.\\
		\end{array}
		$$
	\end{itemize}
\end{exa}

\begin{prp}\label{p5.5.3} Let $(G,*)$ be a hyper group with identity set ${\bf e}=\{e_{\lambda}\;|\; \lambda\in \Lambda\}$, where $\Lambda$ is the lattice of identity set. Assume $H<G_{\lambda}$, $N\lhd G_{\mu}$,  then
	\begin{align}\label{5.5.2}
		HN<G_{\lambda\vee \mu}.
	\end{align}
\end{prp}

\section{Ring and Hyper Ring of Matrices}

\subsection{Matrix Ring}

It has been shown  in Example \ref{e2.1.4} that
$({\cal M}_{n\times n},+,\times)$ is a unitary ring, and it has some sub-rings as
(i) diagonal sub-ring ${\cal D}_n$;
(ii) upper triangular sub-ring ${\cal U}_n$;
(iii) lower triangular sub-ring ${\cal B}_n$; etc.

In addition to these rings with square matrices, are there any rings with non-square matrices?

Recall the DK-STP defined by Definition \ref{d2.3.1}. It leads to a perfect ring structure.

\begin{thm}\label{t6.1.1} \cite{chepr} Consider ${\cal R}:=({\cal M}_{m\times n}, +, \ttimes)$, where $+$ is the classical matrix addition/subtraction. Then ${\cal R}$ is a ring.
\end{thm}

Based on these matrix rings we construct hyper rings.

\subsection{Hyper Ring}

\begin{dfn}\label{d6.2.1} A triple ${\cal R}:=(R,+,\times)$ is called a hyper ring, if the following conditions are satisfied.
	\begin{itemize}
		\item[(i)] $(R,+)$ is an Abelien hyper group.
		\item[(ii)] $(R,\times)$ is a semi-group.
		\item[(iii)]
		\begin{align}\label{6.2.1}
			\begin{array}{l}
				(a+b)\times c=a\times c+ b\times c;\\
				c\times (a+b)=c\times a + c\times b,\quad a,b,c\in R.
			\end{array}
		\end{align}
	\end{itemize}
	If $(R,\times)$ is a hyper-monoid, ${\cal R}$ is called a unitary hyper ring. If $(R,\times)$ is commutative, ${\cal R}$ is called a commutative hyper ring.
\end{dfn}

In the following we construct the most general matrix ring, which covers all matrices.

\begin{thm}\label{t6.2.2} The triple $({\cal M},\hat{+},\hat{\ttimes})$ is a hyper ring, where the addition $\hat{+}$ is defined by (\ref{2.2.11}) and the product $\hat{\ttimes}$ is defined by (\ref{2.3.8}).
\end{thm}

\noindent{\it Proof.} It was shown in Proposition \ref{p2.3.10} that $({\cal M},\hat{\ttimes})$ is a semi-group, and in Proposition \ref{p3.2.1} (iv) that $({\cal M},\hat{+})$ is an abelian hyper group. Hence what remaining to be proved is the distributivity, i.e.,
\begin{align}\label{6.2.2}
	\begin{array}{l}
		(A\hat{+}B)\hat{\ttimes} C=A\hat{\ttimes} C \hat{+} B\hat{\ttimes} C;\\
		C\hat{\ttimes}(A\hat{+} B)=C\hat{\ttimes} A \hat{+} C\hat{\ttimes} B.\\
	\end{array}
\end{align}

We prove the first equality only. Using the first equality and Proposition \ref{p2.3.8}, the second equality is obvious.

Assume $A\in {\cal M}_{m\times n}$, $B\in {\cal M}_{p\times q}$, $C\in {\cal M}_{r\times s}$, and
$$
\begin{array}{c}
	\begin{array}{cc}
		\xi=m\vee p,&\eta=n\vee q,\\
		\a=\xi\vee r=m\vee p\vee r,&\b=\eta\vee s=n\vee q\vee s,\\
	\end{array}\\
	\theta=\a\vee \b=m\vee p\vee r\vee n\vee q\vee s,\\
	\begin{array}{cc}
		u=m\vee r,&v=n\vee s,\\
	\end{array}\\
	w=u\vee v=m\vee n\vee r\vee s,\\
	\begin{array}{cc}
		a=p\vee r,&b=q\vee s,\\
	\end{array}\\
	e=a\vee b=p\vee  r\vee q\vee s,\\
	u\vee a=m\vee p\vee r=\a,\\
	v\vee b=n\vee q\vee s =\b,\\
	\begin{array}{cc}
		\mu=r\vee q,&\lambda=r\vee n.\\
	\end{array}\\
\end{array}
$$
First we have
$$
A\hat{+} B=(A\otimes \E_{\xi/m\times \eta/n})+(B\otimes \E_{\xi/p\times \eta/q}).
$$
Using formula (\ref{2.2.13}), we can calculate that
$$
\begin{array}{l}
	(A\hat{+} B)\hat{\ttimes} C=\{[(A\otimes \E_{\xi/m\times \eta/n})+(B\otimes \E_{\xi/p\times \eta/q})]
	\otimes \E_{\a/\xi\times \b/\eta}\}\ttimes (C\otimes \E_{\a/r\times \b/s})\\
	=[(A\otimes \E_{\a/m\times \b/n})+(B\otimes \E_{\a/p\times \b/q})]\ttimes (C\otimes \E_{\a/r\times \b/s})\\
	=\{[(A\otimes \E_{\a/m\times \b/n})+(B\otimes \E_{\a/p\times \b/q})]\otimes \E^{\rm T}_{\theta/\b}\}[(C\otimes \E_{\a/r\times \b/s})\otimes \E_{\theta/\a}]\\
	=[(A\otimes \E_{\a/m\times \theta/n})+(B\otimes \E_{\a/p\times \theta/q})](C\otimes \E_{\theta/r\times \b/s})\\
	=(A\otimes \E_{\a/m\times \theta/n})(C\otimes \E_{\theta/r\times \b/s})
	+(B\otimes \E_{\a/p\times \theta/q})(C\otimes \E_{\theta/r\times \b/s})\\
	=(A\otimes \E_{\a/m\times \lambda/n})(C\otimes \E_{\lambda/r\times \b/s})
	+(B\otimes \E_{\a/p\times \mu/q})(C\otimes \E_{\mu/r\times \b/s}).\\
\end{array}
$$
On the other hand,
$$
\begin{array}{l}
	A\hat{\ttimes} C= (A\otimes \E_{u/m\times v/n})\ttimes
	(C\otimes \E_{u/r\times v/s}).\\
	=[(A\otimes \E_{u/m\times v/n})\otimes \E^{^{\rm T}}_{w/v}]
	[(C\otimes \E_{u/r\times v/s})\otimes \E_{w/u}]\\
	=(A\otimes \E_{u/m\times w/n})(C\otimes \E_{w/r\times v/s})\\
	=(A\otimes \E_{u/m\times \lambda/n})(C\otimes \E_{\lambda/r\times v/s}).\\
\end{array}
$$
$$
\begin{array}{l}
	B\hat{\ttimes} C= (B\otimes \E_{a/p\times b/q})\ttimes (C\otimes \E_{a/r\times b/s})\\
	=[(B\otimes \E_{a/p\times b/q})\otimes \E^{^{\rm T}}_{e/b}]\left[(C\otimes \E_{a/r\times b/s})\otimes \E_{e/a}\right]\\
	=(B\otimes \E_{a/p\times e/q})(C\otimes \E_{e/r\times b/s})\\
	=(B\otimes \E_{a/p\times \mu/q})(C\otimes \E_{\mu/r\times b/s}).\\
\end{array}
$$
Then
$$
\begin{array}{l}
	(A\hat{\ttimes} C)\hat{+}(A\hat{\ttimes} C)\\
	=[(A\otimes \E_{u/m\times \lambda/n})(C\otimes \E_{\lambda/r\times v/s})]\otimes E_{\a/u\times \b/v}+
	[(B\otimes \E_{a/p\times \mu/q})(C\otimes \E_{\mu/r\times b/s})]
	\otimes \E_{\a/a\times \b/b}\\
	=(A\otimes \E_{\a/m\times \lambda/n})(C\otimes \E_{\lambda/r\times \b/s})+(B\otimes \E_{\a/p\times \mu/q})(C\otimes \E_{\mu/r\times \b/s}).
\end{array}
$$
The first equality of (\ref{6.2.2}) follows.

\hfill $\Box$

\begin{dfn}\label{d6.2.3} Consider a hyper ring $(R,+,\times)$. Assume the identity set of its addition hyper group
	is ${\bf e}=\{e_{\lambda}\;|\;\lambda\in \Lambda\}$.
	\begin{itemize}
		\item[(i)]
		If for each $\lambda\in \Lambda$ the corresponding component group $(R_{\lambda},+)$ satisfies
		\begin{align}\label{6.2.3}
			x\times y\in R_{\lambda},\quad x,y\in R_{\lambda},
		\end{align}
		it is said that the addition $+$ and the product $\times$ of $R$ are consistent.
		\item[(ii)]
		$x,~y\in R$ are said to be equivalent, if there exist $e_{\a},e_{\b}\in {\bf e}$ such that
		\begin{align}\label{6.2.4}
			x+e_{\a}=y+e_{\b}.
		\end{align}
		The equivalence is said to be consistent with operators, if $x_1\sim x_2$ and $y_1\sim y_2$ imply
		\begin{align}\label{6.2.5}
			x_1+y_1\sim x_2+y_2;\quad x_1\times y_1\sim x_2\times y_2.
		\end{align}
	\end{itemize}
\end{dfn}

\begin{prp}\label{p6.2.4} Consider $({\cal M}, \hat{+},\hat{\ttimes})$. (i) $\hat{+}$ and $\hat{\ttimes}$ are consistent; and (ii) the
	equivalence is consistent with operators.
\end{prp}

\noindent{\it Proof.}
\begin{itemize}
	\item[(i)] Assume $A_1\sim A_2$ and $B_1\sim B_2$. Using all the notations/parameters in the proof of Proposition \ref{p3.3.9} (iv), we have
	$$
	A_1\hat{\ttimes} B_1=( A_0\otimes \E_{a/s\times b/t})\ttimes (B_0\otimes \E_{a/\xi\times b/\eta})
	$$
	Assume
	$$
	x:=a\vee b=m\vee \a\vee n\vee \b:= b\mu.
	$$
	Then
	$$
	\begin{array}{l}
		A_1\hat{\ttimes} B\\
=(A_0\otimes \E_{a/s\times b/t}\otimes \E^{^{\rm T}}_{x/b})
		(B_0\otimes \E_{a/\xi\times b/\eta}\otimes \E_{x/a})\\
		=(A_0\otimes \E_{a/s\times x/t})(B_0\otimes \E_{x/\xi\times b/\eta})\\
		=(A_0\otimes \E_{a/s\times b/t}\otimes \E^{^{\rm T}}_{\mu})
		(B_0\otimes \E_{b/\xi\times b/\eta}\otimes \E_{\mu})\\
		~=(A_0\otimes \E_{u/s\times v/t}  \otimes\E_{a/u\times b/v})
		(B_0\otimes \E_{u/\xi\times v/\eta}\otimes\E_{a/u\times b/v}\\
		~=[(A_0\otimes \E_{u/s\times v/t})(B_0\otimes \E_{u/\xi\times v/\eta})]\otimes \E_{a/u\times b/v}.\\
	\end{array}
	$$
	Similarly,
	$$
	%\begin{array}{l}
		A_2\hat{\ttimes} B_2=[(A_0\otimes \E_{u/s\times v/t})(B_0\otimes \E_{u/\xi\times v/\eta})]\otimes \E_{c/u\times d/v}.
	%\end{array}
	$$
	Hence
	$$
	A_1\hat{\ttimes}B_1\sim A_2\hat{\ttimes} B_2.
	$$
	
	\item[(ii)] The fact that the equivalence is consistent with respect to $\hat{+}$ was proved in Proposition \ref{p3.3.9}, and in (i) it has  just been proved that the equivalence is consistent with respect to $\hat{\ttimes}$.
\end{itemize}
\hfill $\Box$

Then similar to the decomposition result of hyper group, we have the following result.

\begin{prp}\label{p6.2.5}  Consider a hyper ring $(R,+,\times)$ with the  identity set of its addition hyper group
	as ${\bf e}=\{e_{\lambda}\;|\;\lambda\in \Lambda\}$.
	\begin{itemize}
		\item[(i)]
		If  the addition $+$ and the product $\times$ of $R$ are consistent, then each $(R_{\lambda},+,\times)$ is a ring, called the component ring of the hyper ring $(R,+,\times)$ .
		\item[(ii)] In addition to (i), if the equivalence is consistent with the operators, then the operators over quotient space can be defined properly by
		\begin{align}\label{6.2.6}
			\begin{array}{l}
				\bar{x}+\bar{y}:=\overline{x+y},\quad x,y\in R,\\
				\bar{x}\times \bar{y}:=\overline{x\times y}.
			\end{array}
		\end{align}
		
		Therefore, the quotient space becomes a ring, called the equivalence ring.
	\end{itemize}
\end{prp}

\begin{figure}
	\centering
	\setlength{\unitlength}{6mm}
	\begin{picture}(10,12)(-0.5,-0.5)
		\thicklines
		\put(6.2,0.2){\line(1,1){2}}
		\put(8.2,2.2){\line(0,1){8.3}}
		\put(0.2,0.2){\line(1,0){6}}
		\put(0.2,0.2){\line(0,1){8.3}}
		\put(0.2,8.5){\line(1,0){6}}
		\put(0.2,8.5){\line(1,1){2}}
		\put(2.2,10.5){\line(1,0){6}}
		\put(6.2,8.5){\line(1,1){2}}
		\put(6.2,8.5){\line(0,-1){8.3}}
		\put(0,0){\vector(0,1){9.5}}
		\put(0,0){\vector(1,0){7.5}}
		\put(0.2,9.5){$R_e$}
		\put(7.2,0.2){$\bar{{\cal R}}$}
		\put(1.2,0.7){$R^1_e$}
		\put(1.2,2.7){$R^2_e$}
		\put(1.2,5.7){$R^n_e$}
		\put(1.2,4.5){$\vdots$}
		\put(1.2,7.5){$\vdots$}
		\put(7.5,9){${\cal R}$}
		\put(-0.4,-0.4){$0$}
		\thinlines
		\put(0.5,0.5){\line(1,0){5}}
		\put(0.5,0.5){\line(1,1){1.5}}
		\put(7,2){\line(-1,0){5}}
		\put(7,2){\line(-1,-1){1.5}}
		\put(0.5,2.5){\line(1,0){5}}
		\put(0.5,2.5){\line(1,1){1.5}}
		\put(7,4){\line(-1,0){5}}
		\put(7,4){\line(-1,-1){1.5}}
		\put(0.5,5.5){\line(1,0){5}}
		\put(0.5,5.5){\line(1,1){1.5}}
		\put(7,7){\line(-1,0){5}}
		\put(7,7){\line(-1,-1){1.5}}
	\end{picture}
	\caption{Ring Decomposition of a hyper ring ${\cal R}$\label{Fig6.2.1}}
\end{figure}
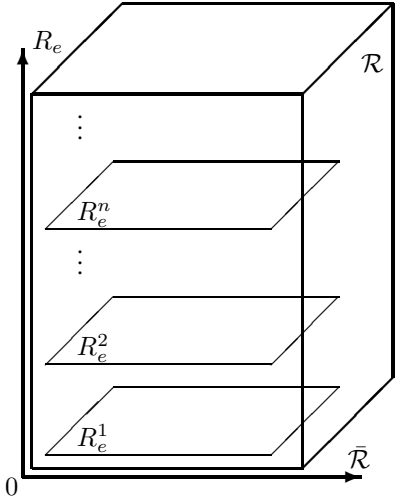

Similarly to hyper group, a diagraph a hyper ring is shown by Figure \ref{Fig6.2.1}.

\begin{exa}\label{e6.2.6}
	Consider ${\cal R}:=({\cal M},\hat{+},\hat{\ttimes})$. To see the consistence of two operators, let $A,B\in M_{e}$, where $e=0_{m\times n}$. Then
	$$
	A\hat{\ttimes}B:=A\ttimes B\in {\cal M}_{m\times n}.
	$$
	The consistence of two operators is proved. Secondly,
	it was proved in Proposition \ref{p6.2.4} that the
	equivalence is consistent with operators. Hence, its decomposition exists.
	Then we have the following decomposed rings:
	\begin{itemize}
		\item[(i)] Component Rings:
		$$
		{\bf e}=\left\{e_{m\times n}=0_{m\times n}\;|\; (m,n)\in \Z_+\times \Z_+\right\}.
		$$
		Each component ring is:
		\begin{align}\label{6.2.7}
			{\cal R}_{e_{m\times n}}=({\cal M}_{m\times n},+,\ttimes ),\quad m,n\in \Z_+.
		\end{align}
		This set of special rings have been described in  the beginning of this section.
		
		\item[(ii)] Equivalence Ring.
		
		By definition, two matrices $A$ and $B$ are equivalent, if there are two identity elements $0_{m\times n}$ and $0_{p\times q}$, such that
		$$
		A\hat{+}0_{m\times n}=B\hat{+} 0_{p\times q},
		$$
		which leads to the following equivalent statement: Let $A\in {\cal M}_{m\times n}$, $B\in {\cal M}_{p\times q}$, $s=\lcm(m,p)$ and $t=\lcm(n,q)$. Then $A\sim B$, if and only if,
		\begin{align}\label{6.2.8}
			A\otimes \E_{s/m\times t/n}= B\otimes \E_{s/p\times t/q}.
		\end{align}
		
		$A\in {\cal M}_{m\times n}$ is reducible, if there exist $(a,b)\prec (m,n)$ and $(a,b)\neq (m,n)$ and $A_0\in {\cal M}_{a\times b}$, such that
		$$
		A=A_0\otimes \E_{m/a\times n/b}.
		$$
		The equivalence class of $A$ is denoted by
		$$
		\bar{A}=\{B\;|\; B\sim A\}.
		$$
		Denote the equivalence ring by
		\begin{align}\label{6.2.9}
			\overline{\cal R}:=(\overline{\cal M},\hat{+},\hat{\ttimes}),
		\end{align}
		the equivalence ring is
		$$
		\overline{\cal M}=\left\{\bar{A}\;|\; A~\mbox{is irreducible}\right\}.
		$$
	\end{itemize}
\end{exa}

\subsection{Square Hyper Ring}

It is well known that $({\cal M}_{n\times n},+,\times)$ is a ring, where $+$ and $\times$ are classical matrix addition and product.

This subsection try to extend this ring structure to a larger set of matrices, which are of mixed dimensions.

\begin{prp}\label{p6.3.2} Assume ${\cal R}=({\cal M}_1, \bar{+},\ltimes)$, where  ${\cal M}_1$ is the set of square matrices; $\bar{+}$ is defined by (\ref{2.2.7}). Then ${\cal R}$ is a hyper ring.
\end{prp}

\noindent{\it Proof.}

It was shown in Proposition \ref{p3.2.1} (ii) that $({\cal M}_1, \bar{+})$ is a commutative hyper group. It is obvious that $({\cal M}_1, \ltimes)$ is a semi-group.

To see that  ${\cal R}$ is a hyper ring, we have only to show the distributivity.

Let $A\in {\cal M}_{m\times m}$, $B\in {\cal M}_{n\times n}$, $C\in {\cal M}_{r\times r}$, $\lcm(m,n)=t$,   $\lcm(t,r)=\lcm(m,n,r)=s$,
$\lcm(m,r)=\a$, $\lcm(n,r)=\b$, and $\lcm(\a,\b)= \lcm(m,n,r)=s$. Then we have
$$
\begin{array}{l}
	(A\bar{+} B)\ltimes C =(A\otimes I_{t/m}+B\otimes I_{t/n}) \ltimes C\\
	~=[(A\otimes I_{t/m}+B\otimes I_{t/n}) \otimes I_{s/t}](C\otimes I_{s/r})\\
	~=(A\otimes I_{s/m})(C\otimes I_{s/r})+(B\otimes I_{s/n})(C\otimes I_{s/r}).\\
\end{array}
$$
$$
\begin{array}{l}
	(A\ltimes C)\bar{+} (B\ltimes C) = (A\otimes I_{\a/m})(C\otimes I_{\a/r})\\
	~\bar{+}(B\otimes I_{\b/n})(C\otimes I_{\b/r})\\
	~=[(A\otimes I_{\a/m})(C\otimes I_{\a/r})]\otimes (I_{s/\a})\\
	~+[(B\otimes I_{\b/n})(C\otimes I_{\b/r})]\otimes I_{s/\b}\\
	~=(A\otimes I_{s/m})(C\otimes I_{s/r})+(B\otimes I_{s/n})(C\otimes I_{s/r}).\\
\end{array}
$$
i.e.,
$$
(A\bar{+} B)\ltimes C=(A\ltimes C)\bar{+} (B\ltimes C).
$$
Similarly, we can show
$$
C\ltimes (A\bar{+} B)=(C\ltimes A)\bar{+} (C\ltimes B).
$$
We conclude that ${\cal R}=({\cal M}_1, \bar{+},\ltimes)$ is a hyper ring.

Finally, we consider the decomposition of this hyper ring. The identity set is
$$
{\bf e}={\bf 0}=\{0_{n\times n}\;|\;n=1,2,\cdots\}.
$$
Then
\begin{itemize}
	\item[(i)] For each $e_n=0_{n\times n}\in {\bf 0}$, ${\cal R}_{e_n}=({\cal M}_{n\times n}, \bar{+},\ltimes)=({\cal M}_{n\times n},+,\times)$ is a component ring.
	
	\item[(ii)]
	It is easy to see that $A\sim B$, if and only if, there exist $I_{\a}$ and $I_{\b}$ such that
	$$
	A\otimes I_{\a}=B\otimes I_{\b}.
	$$
	This is the one defined in \cite{che19}. And the consistence of  the addition with the product and the their consistence with respect to  the equivalence have already been proved in \cite{che19}. We conclude that
	$\overline{{\cal R}}=(\overline{{\cal M}}_1, \bar{+},\ltimes)$ is the equivalence ring.
\end{itemize}

\hfill $\Box$

\begin{prp}\label{p6.3.3}
	
	Assume ${\cal R}=({\cal M}_1, \hat{+},\circ)$, where the product $\circ$ is defined by (\ref{2.2.2}) and the addition is defined as follows:  Assume $A\in {\cal M}_{m\times m}$,  $B\in {\cal M}_{n\times n}$, and $t=\lcm(m,n)$, then
	\begin{align}\label{6.3.2}%{5.2.3}
		A\hat{+} B:=(A\otimes J_{t/m})+(B\otimes J_{t/m}).
	\end{align}
	Then ${\cal R}$ is a hyper ring.
\end{prp}

\noindent{\it Proof.}

A straightforward verification shows that
\begin{itemize}
	\item[(i)] $({\cal M}_1, \hat{+})$ is a hyper group, where the identity set is
	$$
	{\bf e}=\{0_{n\times n}\in {\cal M}_{n\times n}\;|\;n=1,2,\cdots\}.
	$$
	\item[(ii)] $({\cal M}_1, \circ)$ is a semi-group.
\end{itemize}
To see that ${\cal R}$ is a hyper ring, it is enough to show the distributivity.

Assume $C\in {\cal M}_{p \times p}$, and $\a=\lcm(m,p)$, $\b=\lcm(n,p)$, and $s=\lcm(m,n,p)$. Then
$$
\begin{array}{l}
	(A\hat{+}B)\circ C=(A\otimes J_{t/m}+B\otimes J_{t/n})\circ C\\
	~=[(A\otimes J_{t/m}+B\otimes J_{t/n})\otimes J_{s/t}](C\otimes J_{s/p})\\
	~=[(A\otimes J_{s/m}+B\otimes J_{s/n})(C\otimes J_{s/p})\\
	~=(A\otimes J_{s/m})(C+B\otimes J_{t/n})\otimes J_{s/t}](C\otimes J_{s/p})\\
	~=(A\otimes J_{s/m})(C\otimes J_{s/p})+(B\otimes J_{s/n})(C\otimes J_{s/p}).\\
\end{array}
$$
Meanwhile,
$$
\begin{array}{l}
	(A\circ C)\hat{+} (B\circ C)=(A\otimes J_{\a/m})(C\otimes J_{\a/p}\\
	\hat{+} (B\otimes J_{\b/n})(C\otimes J_{\b/p}\\
	~=(A\otimes J_{\a/m})(C\otimes J_{\a/p})\otimes J_{s/\a} \\
	~ +(B\otimes J_{\b/n})(C\otimes J_{b/p})\otimes J_{s/\b}\\
	~=(A\otimes J_{\a/m})(C\otimes J_{\a/p})\otimes J^2_{s/\a}\\
	~  +(B\otimes J_{\b/n})(C\otimes J_{b/p})\otimes J^2_{s/\b}\\
	~=(A\otimes J_{s/m})(C\otimes J_{s/p}+(B\otimes J_{s/n})(C\otimes J_{s/p}).\\
\end{array}
$$
We conclude that
$$
(A\hat{+}B)\circ C =(A\circ C)\hat{+} (B\circ C).
$$

Similarly, we can also show that
$$
C\circ (A\hat{+}B) =(C\circ A)\hat{+} (C\circ B).
$$
We conclude that ${\cal R}=({\cal M}_1, \hat{+},\circ)$ is a hyper ring.

Finally, we consider the decomposition of this hyper ring. The identity set is
$$
{\bf e}={\bf 0}=\{0_{n\times n}\;|\;n=1,2,\cdots\}.
$$
Then
\begin{itemize}
	\item[(i)] For each $e_n=0_{n\times n}\in {\bf 0}$, ${\cal R}_{e_n}=({\cal M}_{n\times n}, \hat{+},\circ) = ({\cal M}_{n\times n}, +,\times)$ is a component ring.
	
	\item[(ii)]
	Since $A\sim B$ means there exist $J_{\a}$ and $J_{\b}$ such that
	$$
	A\otimes J_{\a}=B\otimes J_{\b}.
	$$
	To prove the consistence of  the addition and the product with the equivalence we need to show that
	for $A_1\sim A_2$ and $B_1\sim B_2$ we have
	\begin{align}\label{6.3.3}%{5.2.4}
		A_1\hat{+} B_1\sim A_2\hat{+} B_2;
	\end{align}
	and
	\begin{align}\label{6.3.4}%{5.2.5}
		A_1\circ B_1\sim A_2\circ  B_2.
	\end{align}
	We prove (\ref{6.3.4}) only. The proof of (\ref{6.3.3}) is similar.
	
	Assume $A_1\in {\cal M}_{m\times m}$, $A_2\in {\cal M}_{n\times n}$, $B_1\in {\cal M}_{r\times r}$, and $B_2\in {\cal M}_{s\times s}$ with
	$$
	A_1\sim A_2;\quad B_1\sim B_2.
	$$
	Denote $\gcd(m,n)=\a$, and $\gcd(r,s)=\b$. Using Lemma \ref{l3.3.6},  we know that there exist
	$A_0\in {\cal M}_{\a\times \a}$ and $B_0\in {\cal M}_{\b\times \b}$ such that
	$$
	\begin{array}{l}
		A_1=A_0\otimes J_{m/\a},\quad A_2=A_0\otimes J_{n/\a};\\
		B_1=B_0\otimes J_{r/\b},\quad B_2=B_0\otimes J_{s/\b}.\\
	\end{array}
	$$
	Set $\lcm(m,r)=p$, $\lcm(n,s)=q$, $\lcm(p,q)=\mu$, then we have
	$$
	\begin{array}{l}
		A_1\circ B_1=(A_0\otimes J_{m/\a}\otimes J_{p/m}) (B_0\otimes J_{r/\b}\otimes J_{p/r})\\
		~=(A_0\otimes J_{p/\a}) (B_0\otimes J_{p/\b}).\\
	\end{array}
	$$
	Similarly, we have
	$$
	\begin{array}{l}
		A_2\circ B_2=(A_0\otimes J_{n/\a}\otimes J_{q/n}) (B_0\otimes J_{s/\b}\otimes J_{q/s})\\
		~=(A_0\otimes J_{q/\a}) (B_0\otimes J_{q/\b}).\\
	\end{array}
	$$
	Note that
	$$
	\begin{array}{l}
		(A_1\circ B_1)\otimes J_{\mu/p}=(A_0\otimes J_{p/\a}) (B_0\otimes J_{p/\b})\otimes J^2_{\mu/p}\\
		~=(A_0\otimes J_{\mu/\a}) (B_0\otimes J_{\mu/\b});
	\end{array}
	$$
	and
	$$
	\begin{array}{l}
		(A_2\circ B_2)\otimes J_{\mu/q}=(A_0\otimes J_{q/\a}) (B_0\otimes J_{q/\b})\otimes J^2_{\mu/q}\\
		~=(A_0\otimes J_{\mu/\a}) (B_0\otimes J_{\mu/\b}).
	\end{array}
	$$
	(\ref{6.3.4}) follows.
	
	We conclude that
	$\overline{{\cal R}}=(\overline{{\cal M}}_1, \hat{+},\circ)$ is the equivalence ring.
\end{itemize}

\hfill $\Box$

\subsection{Sub Hyper Ring}

\begin{dfn}\label{d6.4.1} Let $(R,+,\times)$ be a hyper ring and $H\subset R$. $(H,+,\times)$ is said to be a sub-hyper ring of
	$(R,+.\times)$, if (i)  $(H,+)<_h (R,+)$ is a sub-hyper group, and $(H,\times)$ is a sub-semi-group of $(R,\times)$.
\end{dfn}

\begin{prp}\label{p6.4.2}
	\begin{itemize}
		\item[(i)] ${\cal R}_{\mu}=({\cal M}_{\mu}, \hat{+},\hat{\ttimes})$ is a hyper ring,
		where $\hat{+}$ is defined by (\ref{2.2.12}) and $\hat{\ttimes}$ is defined by (\ref{2.3.9}).
		\item[(ii)] ${\cal R}_{\mu}$ is a sub-hyper ring of
		${\cal R}=({\cal M}, \hat{+},\hat{\ttimes})$.
	\end{itemize}
\end{prp}

\noindent{\it Proof.}
\begin{itemize}
	\item[(i)] It has been proved in Proposition \ref{p3.2.1} that $({\cal M}_{\mu}, \hat{+})$ is a hyper group, and in Remark \ref{r2.3.11} that $({\cal M}_{\mu}, \hat{\ttimes})$ is a hyper-semigroup. The proof of distributivity is exactly the same as the one for $({\cal M},
	\hat{+},\hat{\ttimes})$.
	\item[(ii)] A careful observation shows that the  $\hat{+}$ defined by (\ref{2.2.12}) and $\hat{\ttimes}$ defined by (\ref{2.3.9}) respectively are exactly the same as the restrictions of the corresponding $\hat{+}$ and $\hat{\ttimes}$ for ${\cal M}$ or ${\cal M}_{\mu}$. Then the conclusion follows immediately.
\end{itemize}
\hfill $\Box$

\begin{exa}\label{e6.4.3}
	Define
	$$
	\begin{array}{l}
		U^1=\{A\in {\cal M}_1\;|\; A~\mbox{is upper triangular}\};\\
		B^1=\{A\in {\cal M}_1\;|\; A~\mbox{is lower triangular}\};\\
		D^1=\{A\in {\cal M}_1\;|\; A~\mbox{is diagonal}\}.\\
	\end{array}
	$$
	Then (i) $(U^1,\bar{+}, \ltimes)$ is a sub-hyper ring of $({\cal M}_1,\bar{+}, \ltimes)$. (ii) $(B^1,\bar{+}, \ltimes)$ is a sub-hyper ring of $({\cal M}_1,\bar{+}, \ltimes)$. (iii) $(D^1,\bar{+}, \ltimes)$ is a sub-hyper ring of $({\cal M}_1,\bar{+}, \ltimes)$.
\end{exa}

\subsection{Hyper Ring Homomorphism}

\begin{dfn}\label{d6.5.1} Let ${\cal N}=(N,+,\times)$ and ${\cal R}=(R,+,\times)$ be two hyper rings. A mapping $\varphi:N\ra R$ is called a hyper ring homomorphism, if $\varphi: (N,+)\ra (R,+)$ is a hyper group homomorphism, and
	\begin{align}\label{6.5.1}
		\varphi(x\times y)=\varphi(x)\times \varphi(y),\quad x,y\in N.
	\end{align}
\end{dfn}

\begin{prp}\label{p6.5.2} Consider ${\cal R}=({\cal M},\hat{+},\hat{\ttimes})$ and ${\cal N}:= ({\cal M}_1,\hat{+},\circ)$. Let   $\Box: {\cal R}\ra {\cal N}$. Then $\Box$ is a hyper ring homomorphism.
\end{prp}

\noindent{\it Proof.}
It was proved in Proposition \ref{p5.4.4}  that $\Box: ({\cal M},\hat{+})\ra ({\cal M}_1,\hat{+})$ is a hyper group homomorphism.
So we have only to prove (\ref{6.5.1}).

Assume $A\in {\cal M}_{m\times n}$, $B\in {\cal M}_{p\times q}$, $s=m\vee p$, $t=n\vee q$,
$\a=m\vee n$, $\b=p\vee q$,
and $\xi=s\vee t=\a\vee \b$. Then
$$
\begin{array}{l}
	\Box(A\hat{\ttimes} B)=\Box\left[(A\otimes \E_{s/m\times t/n})\ttimes (B\otimes \E_{s/p\times t/q})
	\right]\\
	~=\Box\left[(A\otimes \E_{s/m\times \xi/n})(B\otimes \E_{\xi/p\times t/q})
	\right]\\
	~=(A\otimes \E_{\xi/m\times \xi/n})(B\otimes \E_{\xi/p\times \xi/q}).
\end{array}
$$
$$
\begin{array}{l}
	\Box(A)=A\otimes \E_{\a/m\times \a/n};\\
	\Box(B)=B\otimes \E_{\b/p\times \b/q}.\\
\end{array}
$$
$$
\begin{array}{l}
	\Box(A)\circ\Box(B)=\left[(A\otimes \E_{\a/m\times \a/n})\otimes J_{\xi/\a}\right]\\
	~\left[(B\otimes \E_{\b/p\times \b/q})\otimes J_{\xi/\b}\right]\\
	~=\left[(A\otimes \E_{\a/m\times \a/n})\otimes \E_{\xi/\a\times \xi/\a}\right]\\
	~\left[(B\otimes \E_{\b/p\times \b/q})\otimes \E_{\xi/\b\times \xi/\b}\right]\\
	~=(A\otimes \E_{\xi/m\times \xi/n})(B\otimes \E_{\xi/p\times \xi/q}).
\end{array}
$$
Hence,
$$
\Box(A\hat{\ttimes} B)=\Box(A)\hat{\ttimes} \Box(B).
$$
\hfill $\Box$

\section{Hyper Module and General Linear (Control) System}

\subsection{Hyper Module}

Similarly to the definition of module \cite{hun74}, a hyper module can be defined as follows.

\begin{dfn}\label{d6.1.1} Let $R$ be a hyper ring. A (left) $R$-module is an additive (abelian) hyper group $A$ with a function $\pi:R\times A\ra A$, (briefly denoted by $\pi(r,a)=ra$, $r\in R$, $a\in A$,) such that
	\begin{itemize}
		\item[(i)]
		\begin{align}\label{6.1.1}
			r(a+b)=ra+rb,\quad r\in R,\; a,~b\in A;
		\end{align}
		\item[(ii)]
		\begin{align}\label{6.1.2}
			(r+s)a=ra+sa,\quad r,s\in R,\; a,~b\in A;
		\end{align}
		\item[(iii)]
		\begin{align}\label{6.1.3}
			r(sa)=(rs)a.
		\end{align}
	\end{itemize}
	
	If $R$ has an identity set ${\bf e}_R$ for the product, such that
	\begin{align}\label{6.1.4}
		e_r x\sim x\quad e_r\in {\bf e}_R, \forall x\in A,
	\end{align}
	and there exists at least one $e_x\in {\bf e}_R$, such that
	\begin{align}\label{6.1.5}
		e_x x=x,
	\end{align}
	then $A$ is said to be a unitary hyper R-module.
\end{dfn}

\begin{rem}\label{r6.1.2}
	\begin{itemize}
		\item[(i)] A module can be considered as a particular hyper module, where the $R$ is degenerated to  a ring.
		\item[(ii)] The equality (\ref{6.1.4}) can equivalently be written as
		\begin{align}\label{6.1.6}
			e_R x-x\in {\bf e},\quad e_R\in {\bf e}_R.
		\end{align}
		\item[(iii)] We are particularly interested in unitary (hyper-)R-module.
	\end{itemize}
\end{rem}

Consider the (unitary) module of (hyper-)matrix ring with (hyper-)vector group.

\subsubsection{Hyper Vector Group}

\begin{itemize}
	\item[(i)]
	\begin{align}\label{6.1.7}
		{\cal A}^V_1:=(\R^n,+),\quad n\in \Z_+,
	\end{align}
	where $+$ is the classical vector addition. ${\cal A}^V_1$ is a vector group with unique identity $0_n$.
	\item[(ii)]
	\begin{align}\label{6.1.8}
		{\cal A}^V_2:=(\R^{\infty},\hat{+}),
	\end{align}
	where $\hat{+}$ is defined as follows:
	Let $x\in \R^m$, $y\in \R^n$, and $t=\lcm(m,n)$. Then
	\begin{align}\label{6.1.9}
		x\hat{\pm} y:= (x\otimes \J_{t/m}) + (y\otimes \J_{t/n})\in  \R^{t}.
	\end{align}
	Then ${\cal A}^V_2$ is a hyper-vector group with identity set
	$$
	{\bf e}=\{e_n=0_n\;|\; n\in\Z_+\}.
	$$
	\item[(iii)]
	\begin{align}\label{6.1.10}
		{\cal A}^V_3:=(\R^{\infty},\bar{+}),
	\end{align}
	where $\bar{+}$ is the defined as follows:
	Let $x\in \R^m$, $y\in \R^n$, and $t=\lcm(m,n)$. Then
	\begin{align}\label{6.1.11}
		x\hat{\pm} y:= (x\otimes \E_{t/m}) +  (y\otimes \E_{t/n})\in  \R^{t}.
	\end{align}
	Then ${\cal A}^V_3$ is a hyper-vector group with identity set
	$$
	{\bf e}=\{e_n=0_n\;|\; n\in \Z_+\}.
	$$
\end{itemize}

\subsubsection{Hyper Matrix Ring}

\begin{itemize}
	\item[(i)]
	\begin{align}\label{6.1.12}
		{\cal R}^M_1:=({\cal M}^{n\times n},+,\times ),\quad n\in \Z_+,
	\end{align}
	where $+$ and $\times$ are  the classical matrix addition and production respectively. ${\cal R}^M_1$ is a matrix ring with unique identity $I_n$.
	
	\item[(ii)]
	\begin{align}\label{6.1.13}
		{\cal R}^M_2:=(\overline{\cal M}^{m\times n},+,\ttimes),\quad m,n\in \Z_+,
	\end{align}
	where $+$ is the classical matrix addition and the $\ttimes$ is the DK-STP, defined by (\ref{2.3.1}). ${\cal R}^M_2$ is a matrix ring with the unique identity $I_{m\times n}$.
	
	\item[(iii)]
	\begin{align}\label{6.1.14}
		{\cal R}^M_3:=({\cal M}_1,\bar{+},\ltimes),
	\end{align}
	where $\bar{+}$ is defined by (\ref{2.2.7})  and the $\ltimes$ is defined by (\ref{2.2.1}). ${\cal R}^M_3$ is a matrix unitary hype-ring with the identity set
	$$
	{\bf e}=\{0_{n\times n}\;|\; n\in \Z_+\};
	$$
	and
	$$
	{\bf e}_R=\{I_{n}\;|\; n\in \Z_+\}.
	$$
	\item[(iv)]
	\begin{align}\label{6.1.15}
		{\cal R}^M_4:=({\cal M}_1,\hat{+},\circ),
	\end{align}
	where $\hat{+}$ is defined by (\ref{2.2.11})  and the $\circ$ is defined by (\ref{2.2.2}). ${\cal R}^M_4$ is a matrix unitary hype-ring with the identity set
	$$
	{\bf e}=\{0_{n\times n}\;|\; n\in \Z_+\};
	$$
	and
	$$
	{\bf e}_R=\{I_{n}\;|\; n\in \Z_+\}.
	$$
	\item[(v)]
	\begin{align}\label{6.1.16}
		{\cal R}^M_5:=(\overline{\cal M}_{\mu},\hat{+},\hat{\ttimes}),
	\end{align}
	where the $\hat{\ttimes}$ is defined by (\ref{2.3.8}). ${\cal R}^M_5$ is a matrix unitary hype-ring with the identity set
	$$
	{\bf e}=\{0_{n\mu_y\times n\mu_x}\;|\; n\in \Z_+\};
	$$
	and
	$$
	{\bf e}_R=\{I_{n}\;|\; n\in \Z_+\}.
	$$
	\item[(vi)]
	\begin{align}\label{6.1.17}
		{\cal R}^M_6:=(\overline{\cal M},\hat{+},\hat{\ttimes}).
	\end{align}
	${\cal R}^M_6$ is a matrix unitary hype-ring with the identity set
	$$
	{\bf e}=\{0_{m\times n}\;|\; m,n\in \Z_+\};
	$$
	and
	$$
	{\bf e}_R=\{I_{m\times n}\;|\; m,n\in \Z_+\}.
	$$
\end{itemize}

\subsubsection{Hyper Matrix-Vector Module}

We construct the following combinations of matrix hyper rings with vector hyper groups.
\begin{itemize}
	\item[(i)]
	$\aleph_1:=({\cal R}^M_1, {\cal A}^V_1,\times)$,  where $\times$ is the classical matrix-vector product.
	
	\item[(ii)]
	$\aleph_2:=({\cal R}^M_2, {\cal A}^V_2,\ttimes)$ (or $({\cal R}^M_2, {\cal A}^V_3,\ttimes)$ ), where
	$\ttimes$ is the DK-STP, defined by (\ref{2.3.1}).
	
	\item[(iii)]
	$\aleph_3:=({\cal R}^M_3, {\cal A}^V_2,\ltimes)$ (or $({\cal R}^M_3, {\cal A}^V_3,\ltimes)$).
	
	\item[(iv)]
	$\aleph_4:=({\cal R}^M_4, {\cal A}^V_2,\circ)$ (or $({\cal R}^M_4, {\cal A}^V_3,\circ)$ .
	
	\item[(v)]
	$\aleph_5:=({\cal R}^M_5, {\cal A}^V_2,\ttimes)$ (or $({\cal R}^M_5, {\cal A}^V_3,\ttimes)$).
	
	\item[(vi)]
	$\aleph_6:=({\cal R}^M_6, {\cal A}^V_2,\ttimes)$ (or$({\cal R}^M_6, {\cal A}^V_3,\ttimes)$).

\end{itemize}

\begin{prp}\label{p6.1.3}
	\begin{itemize}
		\item[(i)] The  $\aleph_1$ and $\aleph_2$ are unitary modules.
		\item[(ii)] The $\aleph_3$, $\aleph_4$, $\aleph_5$ and $\aleph_6$ are unitary hyper modules.
	\end{itemize}
\end{prp}

\noindent {\it Proof.} We prove $\aleph_2$, the proof of others is similar.
What do we need to verify is (\ref{6.1.1})-(\ref{6.1.3}). In fact, (\ref{6.1.2}) is trivial, and (\ref{6.1.3}) is well know \cite{che23}. Hence, we prove (\ref{6.1.1}) only.

Let $A\in {\cal M}_{m\times n}$, $x\in \R^p$, $y\in \R^q$, $p\vee n=\a$, $q\vee n=\b$, $p\vee q=t$, and $n\vee p\vee q=s$. Then
$$
\begin{array}{l}
	A\ttimes x=(A\otimes \J^{\rm T}_{\a/n})(x\otimes \J_{\a/p})\in \R^m;\\
	A\ttimes y=(A\otimes \J^{\rm T}_{\b/n})(y\otimes \J_{\b/q})\in \R^m.\\
\end{array}
$$
$$
\begin{array}{l}
	A\ttimes(x\hat{+}y)=A\ttimes (x\otimes \J_{t/p}+y\otimes \J_{t/q})\\
	~=(A\otimes \J^{\rm T}_{s/n})( x\otimes \J_{s/p}+y\otimes \J_{s/q}  )\\
	~=(A\otimes \J^{\rm T}_{s/n})( x\otimes \J_{s/p})+
	(A\otimes \J^{\rm T}_{s/n})(y\otimes \J_{s/q}  )\\
	~=(A\otimes \J^{\rm T}_{\a/n})( x\otimes \J_{\a/p})+
	(A\otimes \J^{\rm T}_{\b/n})(y\otimes \J_{\b/q}  )\\
\end{array}
$$
(\ref{6.1.1}) follows.

\hfill $\Box$

\subsection{Linear (Control) System}

\begin{dfn}\label{d6.2.1} Let $\aleph=({\cal R}_0,{\cal A}_0, \pi) $ be a hyper matrix-vector (MV-) module.
	\begin{itemize}
		\item[(i)] A discrete-time linear system is defined as
		\begin{align}\label{6.2.1}
		%	\begin{array}{l}
				x(t+1)=M(t)x(t),\quad x(0)=x_0\in {\cal V}_0,%\\
				\quad M(t)\in {\cal R}_0,\; x(t)\in {\cal A}_0,
		%	\end{array}
		\end{align}
		where $M(t)x(t)=\pi(M(t),x(t))$.
		\item[(ii)] A continuous-time linear system is defined as
		\begin{align}\label{6.2.2}
		%	\begin{array}{l}
				\dot{x}=M(t)x(t),\quad x(0)=x_0\in {\cal V}_0,%\\
				\quad M(t)\in {\cal R}_0,\; x(t)\in {\cal A}_0.
		%	\end{array}
		\end{align}
		\item[(iii)] A discrete-time linear control system is defined as
		\begin{align}\label{6.2.3}
			%\begin{array}{l}
				x(t+1)=M(t)x(t) + B(t)u(t),\quad x(0)=x_0\in {\cal V}_0,%\\
				\quad M(t),B(t)\in {\cal R}_0,\; x(t),u(t)\in {\cal A}_0.
			%\end{array}
		\end{align}
		\item[(iv)] A continuous-time linear control system is defined as
		\begin{align}\label{6.2.4}
			%\begin{array}{l}
				\dot{x}=M(t)x(t)+B(t)u(t),\quad x(0)=x_0\in {\cal V}_0,%\\
				\quad M(t),B(t)\in {\cal R}_0,\; x(t)\in {\cal A}_0.
			%\end{array}
		\end{align}
	\end{itemize}
\end{dfn}

\begin{rem}\label{r6.2.2}
	\begin{itemize}
		\item[(i)] When $M(t)=M$ (and $B(t)=B$ for control systems) the systems are called the time-invariant linear (control) systems.
		\item[(ii)] We are only interested in the case when $\aleph\in \{\aleph^i\;|\;i\in [1,6]\}$, which were defined in the previous subsection.
		\item[(iii)] When $\aleph=\aleph^1$, the systems are the classical linear (control) systems, which have been investigated over 8 decades \cite{won74}.
		\item[(iv)] When $\aleph=\aleph^3$ or $\aleph=\aleph^4$, the systems have been investigated by \cite{che19}. Particularly, the trajectories of the time-invariant systems have been constructed.
	\end{itemize}
\end{rem}

\begin{exa}\label{e6.2.3} Consider a linear control system
	\begin{align}\label{6.2.5}
		\dot{x}=Ax(t)+Bu(t)+C\eta(t),
	\end{align}
	where $x(t)\in \R^n$ is the state, $u(t)\in \R^m$ is the control, $\eta(t)\in \R^r$ is the disturbance, $A\in {\cal M}_{n\times n}$,
	$B\in {\cal M}_{n\times m}$,$C\in {\cal M}_{n\times r}$. Assume $\eta(t)$ is a dimension-varying disturbance, then we may replace (\ref{6.2.5}) by
	$$
	\dot{x}=Ax(t)+Bu(t)+C\ttimes \eta(t).
	$$
	Then the system can handle disturbance with various dimensions.
	If there are some state variable dimension-varying vibrations, say,
	$$
	x(t_{\tau})=x(t_{\tau_{-}}),\quad \tau\in \{t_1,t_2,\cdots\},
	$$
	Then  (\ref{6.2.5}) can be replaced by
	$$
	\dot{x}=A\ttimes x(t)+Bu(t)+C\ttimes \eta(t).
	$$
\end{exa}

\section{Hyper Vector Space}

Up to now the algebraic structures of MVMDs have been investigated. From this section on, we will consider the geometric structure of MVMDs.

\subsection{From Hyper Group to Hyper Vector Space}

To be precise we first review the concept of vector space in linear algebra.

\begin{dfn}\label{d7.1.1} \cite{gre81} A vector space ${\cal V}$ over $\R$ is a triple
	${\cal V}=(X,+,\cdot)$ , where $X$ is a set with elements $x\in X$ are vectors;  addition $+:X\times X\ra X$, and scalar product $\cdot:\R\times X\ra X$, satisfying
	\begin{enumerate}
		\item[(1)] $(X,+)$ is an abelian group, i.e.,
		
		\begin{itemize}
			\item[(i)] Associativity:
			\begin{align}\label{7.1.1}
				(x+y)+z=x+(y+z),\quad x,y,z\in X.
			\end{align}
			\item[(ii)] Commutativity:
			\begin{align}\label{7.1.2}
				x+y=y+x,\quad x,y\in X.
			\end{align}
			\item[(iii)] Zero: There exists a unique $\vec{0}\in X$, such that
			\begin{align}\label{7.1.3}%{2.1.3}
				x+\vec{0}=\vec{0}+x=x,\quad \forall x\in X.
			\end{align}
			\item[(iv)] Inverse: For each $x\in X$ there exists an inverse $-x\in X$, such that
			\begin{align}\label{7.1.4}%{2.1.4}
				x+(-x)=0,\quad x\in X.
			\end{align}
		\end{itemize}
		\item[(2)] Scalar product satisfies
		\begin{itemize}
			\item[(i)] Associativity:
			\begin{align}\label{7.1.5}%{2.1.5}
				(r_1r_2)\cdot x=r_1\cdot(r_2\cdot x),\quad r_1,r_2\in \F,~x\in X.
			\end{align}
			\item[(ii)] Distributivity:
			\begin{align}\label{7.1.6}%{2.1.6}
				\begin{array}{l}
					(r_1+r_2)\cdot x=r_1\cdot x+r_2\cdot x,\quad r,r_1,r_2\in \R,\\
					r\cdot(x+y)=r\cdot x+r\cdot y,\quad ~x,~y\in X.
				\end{array}
			\end{align}
			\item[(iii)] Unit:
			\begin{align}\label{7.1.7}%{2.1.7}
				1\cdot x=x, \quad 1\in \R, x\in X.
			\end{align}
		\end{itemize}
	\end{enumerate}
\end{dfn}

When the MVMDs are considered, it is natural to consider ``pseudo-" vector space, which is defined by modifying the definition of vector space a little bit.

\begin{dfn}\label{d7.1.2}  A hyper vector space ${\cal V}$ over $\R$ is a triple
	${\cal V}=(X,+,\cdot)$  as in Definition \ref{d7.1.1}, satisfying
	\begin{enumerate}
		\item[(1')] $(X,+)$ ia an abelian hyper group.
		\item[(2)] Same as the (2) in Definition \ref{d7.1.1}.
	\end{enumerate}
\end{dfn}

Recall that each hyper group has its component groups and as the operator (it is addition now) is consistent with respect to equivalence then the equivalence group exists. Correspondingly, we have the following result.

\begin{prp}\label{p7.1.3} Consider  a hyper vector space ${\cal V}=(X,+,\cdot)$. Assume ${\bf e}=\{e_{\lambda}\;|\;\lambda\in \Lambda\}$. Then
	\begin{itemize}
		\item[(i)] For each $\lambda\in \Lambda$, $X_{\lambda}$ is a vector space.
		\item[(ii)] Assume the  addition is consistent with respect to equivalence, then the equivalence group is a vector space.
	\end{itemize}
\end{prp}

\noindent{\it Proof.}
\begin{itemize}
	\item[(i)] Let $x\in V_{\lambda}$. Then,
	$$
	xe_{\lambda}=e_{\lambda}x=x, \quad x^{-1}e_{\lambda}=e_{\lambda}x^{-1}=x^{-1}.
	$$
	It follows that for any $\a\in \R$
	$$
	\begin{array}{l}
		(\a x)e_{\lambda}=e_{\lambda}(\a x)=\a x,\\
		(\a x)^{-1}e_{\lambda}=e_{\lambda}(\a x)^{-1}=(\a x)^{-1}.
	\end{array}
	$$
	Hence $\a x\in X_{\lambda}$, which ensures that $V_{\lambda}$ is a vector space.
	\item[(ii)] Since
	$$
	\overline{\a x}=\a \bar{x},\quad \a\in \R,
	$$
	$(\bar{V},+)$ is also a vector space.
\end{itemize}
\hfill $\Box$

\subsection{Hyper Vector Space over Matrices}

It was mentioned in Remark \ref{r3.2.2} that all the addition hyper groups on MVMDs are abelian. In addition, for MVMDs the scaler product is naturally defined. Moreover, the properties of scalar product required for vector spaces are automatically satisfied. Then based on Proposition \ref{p3.2.1} the following results are obvious.

\begin{prp}\label{p7.2.1}
	\begin{itemize}
		\item[(i)] $(\R^{\infty}, \hat{+})$ is a hyper vector space.
		
		\item[(ii)]  $({\cal M}_1,\bar{+})$ is a hyper vector space.
		$({\cal M}_1,\hat{+})$ is also a hyper vector space.
		
		\item[(iii)] $({\cal M}_{\mu}. \hat{+})$ is a hyper vector space  with $\hat{+}$ defined by (\ref{2.2.8}), where the identity set is
		the same as (\ref{3.2.1}).
		
		\item[(iv)]
		$({\cal M}, \hat{+})$ is a hyper vector space.
	\end{itemize}
\end{prp}

The decomposition for hyper groups can be immediately extended to hyper vector space. Then we have

\begin{prp}\label{p7.2.1} Let $(X,+,{\bf e})$ be a hyper vector space.
	\begin{itemize}
		\item[(i)] For each $e\in {\bf e}$, $(X_e,+)$  is a vector space, called the component subspace of the hyper vector space.
		
		\item[(ii)] Assume the equivalence is consistent with operator $+$ and scalar product, then the quotient space is also a vector
		space called the equivalence vector space.
		
	\end{itemize}
\end{prp}

\section{Metric Space Structure on $\R^{\infty}$ }

\subsection{Hyper Inner Product Space}

Replacing vector space by hyper vector space, the hyper inner product space is defined as follows  \cite{tay80}.

\begin{dfn}\label{d8.1.1} A hyper vector space $X$ over $\F$  ($\F=\mathbb{C}$ or $\F=\R$) is called an hyper inner product space if to each pair of elements $x,~y\in X$ there is associated a number $\langle x,y\rangle \in \F$, called the hyper inner product of $x$ and $y$, with the following properties:
	\begin{itemize}
		\item[(i)]
		\begin{align}\label{8.1.1}
			\langle x+y,z\rangle=\langle x,z\rangle+\langle y,z\rangle,\quad x,y,z\in X.
		\end{align}
		\item[(ii)]
		\begin{align}\label{8.1.2}
			\langle x,y\rangle=\overline{\langle y,x\rangle}.
		\end{align}
		\item[(iii)]
		\begin{align}\label{8.1.3}
			\langle ax,y\rangle=a\langle x,y\rangle.
		\end{align}
		\item[(iv)]
		\begin{align}\label{8.1.4}
			\langle x,x\rangle\geq 0, ~\mbox{and}~\langle x,x\rangle\neq 0, ~\mbox{if}~x\not\in {\bf 0}.
		\end{align}
	\end{itemize}
\end{dfn}

Consider $\R^{\infty}$, an inner product  has been defined as follows.

\begin{dfn}\label{d8.1.2} \cite{che19} Consider $x,y\in \R^{\infty}$, say, $x\in \R^m$, $y\in \R^n$, and $\lcm(m,n)=t$
	The inner product of $x$ and $y$ is defined by
	\begin{align}\label{8.1.5}
		%\begin{array}{l}
			\langle x,y\rangle_{{\cal V}}:=\frac{1}{t}x\vec{\cdot}y%\\
			=\frac{1}{t}(x\otimes \J_{t/m})^{^{\rm T}} (y\otimes \J_{t/n}).
		%\end{array}
	\end{align}
\end{dfn}

Recall that $(\R^{\infty},\hat{+},\cdot)$ is a hyper vector space,
it is readily verifying that the inner product defined by (\ref{8.1.5}) is a hyper inner product.

Note that by Definition \ref{d8.1.2} if $x\in \R^n$, then
\begin{align}\label{8.1.501}
	\langle x,x\rangle_{{\cal V}}=\frac{1}{n}x^{^{\rm T}}x=\frac{1}{n}\langle x,x\rangle,
\end{align}
where $\langle x,x\rangle$ is the Euclidean inner product. Hence the inner product defined by (\ref{8.1.5}) is the ``normalized"  Euclidean inner product. This is very important for comparing the vectors in Euclidean spaces of various dimensions.

Alternatively, we can also define
\begin{align}\label{8.1.501}
	\langle x,y\rangle_{{\cal V}}:=(x\otimes \E_{t/m})^{^{\rm T}} (y\otimes \E_{t/n}).
\end{align}

For a hyper inner product space, the inner-product is said to be consistent with the operator (addition) if
$x_1\sim x_2$ and $y_1\sim y_2$ imply
$$
\langle x_1, y_1 \rangle= \langle x_2, y_2 \rangle.
$$

The following result is an immediate consequence of the definition.

\begin{prp}\label{p8.1.3} Let ${\cal H}:=({\cal V}, \langle\cdot,\cdot\rangle)$ be a hyper inner product space, where ${\cal V}$ is a hyper vector space and $\langle\cdot,\cdot\rangle$ is the inner product. If the inner-product is consistent with the addition. Then the equivalence space, denoted by $\overline{\cal H}:=(\overline{\cal V}, \langle\cdot,\cdot\rangle)$,
	is an inner product space, with the inner product as
	\begin{align}\label{8.1.6}
		\langle \bar{x},\bar{y}\rangle:=\langle x,y \rangle.
	\end{align}
\end{prp}

\begin{exa}\label{e8.1.4} Consider the hyper inner product space ${\cal H}:=(\R^{\infty},\hat{+},\langle \cdot,\cdot\rangle)$, where the inner product is defined by (\ref{8.1.5}). It has been proved in \cite{che19} that this inner product is consistent with the addition, then it is clear that the equivalence space $\overline{\cal H}:=(\overline{\R}^{\infty},\hat{+},\langle \cdot,\cdot\rangle)$ is an inner-product space. It has also been proved in \cite{che19} that this inner-product space is not complete, hence it is not a Hilbert space.
\end{exa}

\begin{dfn}\label{d8.1.5}
	\begin{enumerate}
		\item[(1)] A norm on a hyper vector space $X$ is a real function, whose value at $x\in X$, denoted by $\|x\|$, with the properties:
		\begin{itemize}
			\item[(i)]
			\begin{align}\label{8.1.7}
				\|x+y\|\leq \|x\|+\|y\|,\quad x,y\in X.
			\end{align}
			\item[(ii)]
			\begin{align}\label{8.1.8}
				\|ax\|=|a|\|x\|.
			\end{align}
			\item[(iii)]
			\begin{align}\label{8.1.9}
				\|x\|\geq 0,~\mbox{and}~\|x\|\neq 0,~\mbox{if}~x\not\in {\bf 0}.
			\end{align}
		\end{itemize}
		\item[(2)] Assume ${\cal V}$ is a hyper inner product space, i.e., it is a hyper vector space with an inner product $\langle\cdot, \cdot \rangle$, then a norm can be defined by
		\begin{align}\label{8.1.10}
			\|x\|:=\sqrt{\langle x,x\rangle}.
		\end{align}
		\item[(3)]   The norm is said to be consistent with the equivalence, if $x\sim y$ implies $\|x\|=\|y\|$.
	\end{enumerate}
	A hyper vector space with a norm is called a normed hyper vector space.
\end{dfn}

The following is an immediate consequence of definitions.

\begin{prp}\label{p8.1.6} Consider a normed hyper vector space. If the norm is consistent with the equivalence, its equivalence space is a normed vector space.
\end{prp}

\begin{exa}\label{e8.1.7} Consider the hyper inner product space ${\cal H}:=(\R^{\infty},\hat{+},\langle \cdot,\cdot\rangle)$.
	\begin{itemize}
		\item[(i)] Its norm has been defined in \cite{che19} as
		\begin{align}\label{8.1.11}
			\|x\|_{{\cal V}}:=\sqrt{\left<x,x\right>_{{\cal V}}}.
		\end{align}
		\item[(ii)] It has been proved in \cite{che19} that the norm defined by (\ref{8.1.10}) is consistent with the equivalence.
		Hence the equivalence space $(\overline{\R}^{\infty},\hat{+},\langle \cdot,\cdot\rangle)$ is a normed vector space.
	\end{itemize}
\end{exa}

Let $X$ be a normed hyper vector space, and $x,y\in X$. A natural way to define the distance between $x$ and $y$ is
\begin{align}\label{8.1.1101}
	d(x,y):=\|x - y\|.
\end{align}

As an example, consider $\R^{\infty}$.  Let $x,y\in \R^{\infty}$. The distance of $x$ and $y$ is defined by \cite{che19}
\begin{align}\label{8.1.12}
	d_{{\cal V}}(x,y):=\|x - y\|_{{\cal V}}.
\end{align}

It is readily verifying that the distance defined in (\ref{8.1.1101}), or particularly (\ref{8.1.12}) for $\R^{\infty}$,  is a hyper-metric, which turns the hyper vector space into a hyper-metric space.

\begin{prp}\label{p8.1.8}
	\begin{itemize}
		\item[(i)] The metric defined by (\ref{8.1.12}) is a hyper-metric, which makes $(\R^{\infty},d)$ a hyper-metric space.
		\item[(ii)]  The hyper-metric space  is a special pseudo-matric space.
		\item[(iii)] The metric defined by (\ref{8.1.12}) is consistent with equivalence, i.e., assume $x_1\sim x_2$ and $y_1\sim y_2$, then $d(x_1,y_1) = d(x_2,y_2)$.
		\item[(iv)] The equivalence space $\bar{\R}^{\infty}$ is a metric space.
	\end{itemize}
\end{prp}

The topology deduced by the hyper-metric $d_{{\cal V}}$, defined by (\ref{8.1.12}) and denoted by ${\cal T}_d$, makes $\R^{\infty}$ a topological space. This topological space has the following properties \cite{che19}:

\begin{prp}\label{p8.1.9} Consider the topological space $(\R^{\infty},{\cal T}_d)$.
	\begin{itemize}
		\item[(i)] It is separable, i.e., there exists a countable subset, which is dense in $\R^{\infty}$.
		\item[(ii)] It is not a Hausdorff space. (In fact, it is not even $T_0$).
	\end{itemize}
\end{prp}

Finally, we briefly review its equivalence space. It is easy to verify the following equivalent statements.

\begin{prp}\label{p8.1.9} Let $x,~y\in \R^{\infty}$. $x$ and $y$  are equivalent, if one of the following  equivalent statements holds.
	\begin{itemize}
		\item[(i)] $x,~y$ are equivalent, denoted by
		$x\lra y$.
		\item[(ii)]  There exist $0_p$ and $0_q$  such that
		$x\hat{+} 0_p=y\hat{+} 0_q$.
		\item[(iii)] There exist
		$\J_p$ and $\J_q$ such that
		$x\otimes \J_p=y\otimes \J_q$.
		\item[(iv)]
		$d_{{\cal V}}(x,y)=0$.
	\end{itemize}
\end{prp}
Denote the quotient space by
\begin{align}\label{8.1.13}
	\Omega:=\R^{\infty}/\lra.
\end{align}

According to the Proposition \ref{p8.1.9}, it is clear that $\Omega$ is both the equivalence space and the topological quotient space.
Moreover, this space is  a Hausdorff space.

\begin{rem}\label{r8.1.10} Note that if $x,~y\in \R^n$ and $x\lra y$, then $x=y$. So
\begin{align}\label{8.1.14}
\Omega_n:=\R^n/\lra\cong \R^n.
\end{align}
Hence the coordinates of $\R^n$ can be used as the coordinates of $\Omega_n$. Now $\Omega=\bigcup_{n=1}^{\infty}\Omega_n$. It follows that $x\in \Omega$ can have infinite ``coordinate neighborhoods", (which are global coordinate charts). For instance,
\begin{itemize}
\item[(i)] let $x=a\in \R^1$. The coordinate neighborhoods of $\bar{x}$ are
$$
\bar{x}\in \Omega_1\subset \Omega_2\subset \Omega_3\subset \cdots;
$$
\item[(ii)] let $x=(a,b)^{\mathrm{T}}\in \R^2$. The coordinate neighborhoods of $\bar{x}$ are
$$
\bar{x}\in \Omega_2\subset \Omega_4\subset \Omega_6\subset \cdots.
$$
\end{itemize}
If the cross-dimensional manifold $M$ is considered. Then locally it is homeomorphic to an open subset of $\Omega$. Hence,
 \begin{itemize}
\item[(i)] the coordinate neighborhoods of $\bar{x}\in M$ are
$$
\bar{x}\in O_1\subset O_2\subset O_3\subset \cdots,
$$
where $O_i\subset \Omega_i$, $i=1,2,3,\cdots$ are open subsets.

\item[(ii)] the coordinate neighborhoods of $\bar{x}$ are
$$
\bar{x}\in O_2\subset O_4\subset O_6\subset \cdots.
$$
\end{itemize}
\end{rem}

We refer to \cite{cho66} for related topics about topological spaces; and to  \cite{tay80,con85} for some related knowledge for inner product space.

\subsection{Geometric Structure on $\R^{\infty}$}%

First, we consider the projection from $\R^m$ to $\R^n$.

\begin{dfn} \label{d8.2.1}
	\cite{che19} Let $\xi\in \R^m$. The projection of $\xi$ on $\R^n$, denoted by $\pi^m_n(\xi)$, is defined as
	\begin{align}\label{8.2.1}
		\pi^m_n(\xi):=\argmin_{x\in \R_n}\|\xi \hat{-} x\|_{{\cal V}}.
	\end{align}
\end{dfn}

The projection can be expressed by matrix product form as follows.

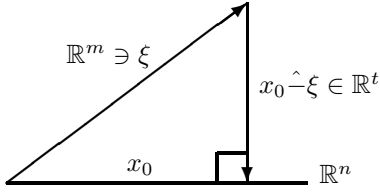
\begin{figure}
	\centering
	\setlength{\unitlength}{8mm}
	\begin{picture}(8,5)\thicklines
		\put(1,1){\line(1,0){5}}
		\put(1,1){\vector(4,3){4}}
		\put(5,4){\vector(0,-1){3}}
		\put(4.5,1.5){\line(1,0){0.5}}
		\put(4.5,1.5){\line(0,-1){0.5}}
		\put(2,3){$\R^{m}\ni\xi$}
		\put(3,1.2){$x_0$}
		\put(5.2,2.5){$x_0\hat{-}\xi\in \R^{rm T}$}
		\put(6.2,1){$\R^n$}
	\end{picture}
	\vskip -5mm
	\caption{Cross-dimensional Projection\label{Fig8.2.1}}
\end{figure}

\begin{prp} \label{p8.2.2}
	\cite{che19}
	\begin{align}\label{8.2.2}
		x_0:=\pi^m_n(\xi)=\Pi^m_n\xi \in \R^n,\quad \xi\in \R^m,
	\end{align}
	where ($t=\lcm(m,n)$)
	\begin{align}\label{8.2.3}
		\Pi^m_n=\frac{1}{t/n} ( I_n\otimes \J^{^{\rm T}}_{t/n})  ( I_m\otimes \J_{t/m}) ,
	\end{align}
	Moreover, $ x_0 \bot (\xi\hat{-} x_0 )$ (see Figure \ref{Fig8.2.1}).
\end{prp}

We give an example:

\begin{exa}\label{e8.2.3}
	Assume $m=4$, $n=6$. Let $\xi=(\xi_1,\xi_2,\xi_3,\xi_4)^{^{\rm T}}\in \R^m$. We consider its projection on $\R^6$.
	First,  $t=\lcm(m,n)=12$, then
	$$
	\begin{array}{l}
		\Pi^4_6=\frac{1}{2}(I_6\otimes \J^{^{\rm T}}_{2})(I_4\otimes \J_{3})\\
		=\begin{bmatrix}
			1&0&0&0\\
			0.5&0.5&0&0\\
			0&1&0&0\\
			0&0&1&0\\
			0&0&0.5&0.5\\
			0&0&0&1\\
		\end{bmatrix}.
	\end{array}
	$$
	Hence,
	$$
	x_0=(\xi_1,0.5\xi_1+0.5\xi_2,\xi_2,\xi_3,0.5\xi_3+0.5\xi_4,\xi_4)^{^{\rm T}}.
	$$
	$$
	\begin{array}{l}
		x_0\hat{-}\xi=(0,0,0.5\xi_2-0.5\xi_1,0.5\xi_1-0.5\xi_2,0,0\\
		~0,0,0.5\xi_4-0.5\xi_3,0.5\xi_3-0.5\xi_4,0,0)^{^{\rm T}}.
	\end{array}
	$$
	It is ready to check that
	$$
	\langle x_0\hat{-} \xi, x_0\rangle_{{\cal V}}=0.
	$$
\end{exa}

\begin{figure}
	\centering
	\setlength{\unitlength}{8mm}
	\begin{picture}(8,5)\thicklines
		\put(5,1){\vector(-1,0){4}}
		\put(1,1){\vector(4,3){4}}
		\put(5,1){\vector(0,1){3}}
		\put(4.5,1.5){\line(1,0){0.5}}
		\put(4.5,1.5){\line(0,-1){0.5}}
		\put(2.5,3){$\xi\hat{-} z$}
		\put(5,0.5){$x_0$}
		\put(5.2,2.5){$\xi\hat{-}x_0\in \R^{rm T}$}
		\put(3,1.2){$\R^n$}
		\put(0.8,0.5){$z$}
		\put(2.5,0.5){$z-x_0$}
		\put(5,4.2){$\xi\in \R^m$}
	\end{picture}
	\vskip -5mm
	\caption{Projection is orthogonal to $\R^n$\label{Fig8.2.2}}
\end{figure}
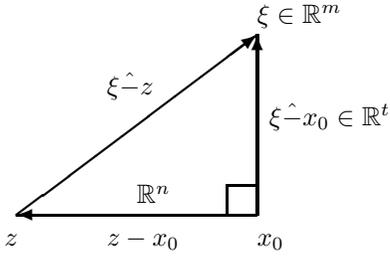

\begin{thm}\label{t8.2.4}
	(Generalized Pythagoras' Theorem)
	\begin{itemize}
		\item[(i)] Let $\xi\in \R^m$, $x_0=\pi^m_n(\xi)\in \R^n$, $x_0\hat{-}\xi\in \R^{rm T}$, where $t=\lcm(m,n)$. Then
		\begin{align}\label{8.2.4}
			\|\xi\|_{{\cal V}}^2=\|x_0\|_{{\cal V}}^2+\|\xi-x_0\|_{{\cal V}}^2.
		\end{align}
		
		Figure \ref{Fig8.2.1} depicts the relationship of above vectors.
		
		\item[(ii)]  Let $\xi\in \R^m$, and the projection of $\xi$ on $\R^n$ be $x_0=\pi^m_n(\xi)$. Then
		$\xi-x_0$ is perpendicular to $\R^n$, i.e., for any $z\in \R^n$ we have
		\begin{align}\label{8.2.5}
			\langle (\xi\hat{-}x_0),(z-x_0)\rangle_{{\cal V}} =0.
		\end{align}
		
		Moreover,
		\begin{align}\label{8.2.6}
			\|\xi\hat{-}z\|^2_{{\cal V}}=\|\xi\hat{-}x_0\|^2_{{\cal V}}+\|z\hat{-}x_0\|^2_{{\cal V}}.
		\end{align}
		
		Figure \ref{Fig8.2.2} depicts the relationship of the vectors.
		
	\end{itemize}
\end{thm}

\noindent{\it Proof.} We prove (ii) only. The proof of (i) is similar.

Embedding $\R^n$ into $\R^{rm T}$ by $x\mapsto x\otimes \J_{t/n}$ and using Proposition \ref{p8.1.3},
we have
$$
(\xi\hat{-}x_0)\perp (\R^n\otimes \J_{t/n}).
$$
Then
$$
(\xi\hat{-}x_0)\perp (z-x_0\otimes \J_{t/n}).
$$
Using  Pythagoras' Theorem on $\R^{rm T}$, we have
$$
\|\xi\hat{-}z\|^2_{{\cal V}}=\|\xi\hat{-}x_0\|^2_{{\cal V}}+\|(z-x_0)\otimes \J_{t/n}\|^2_{{\cal V}}.
$$
Then Proposition  \ref{p8.1.3} yields (\ref{8.2.6})
immediately.
\hfill $\Box$

\begin{exa}\label{e8.2.5}
	Recall Example \ref{e8.2.3},
	a straightforward computation shows that
	$$
	\|\xi\|_{{\cal V}}^2=\frac{1}{4}(\xi_1^2+\xi_2^2+\xi_3^2+\xi_4^2);
	$$
	$$
		\|x_0\|_{{\cal V}}^2=\frac{1}{6}(1.25\xi_1^2 +1.25\xi_2^2+1.25\xi_3^2+1.25 \xi_4^2
		+0.5\xi_1\xi_2+0.5\xi_3\xi_4).
	$$
	$$
		\|x_0\hat{-}\xi\|_{{\cal V}}^2=\frac{1}{12}(0.5\xi_1^2 +0.5\xi_2^2+0.5\xi_3^2
	+0.5 \xi_4^2+\xi_1\xi_2+\xi_3\xi_4).
	$$
	Using these equalities,  (\ref{8.2.4}) follows immediately.%{2.3.9}).
\end{exa}

Using  projection and Pythagoras' Theorem, more geometric structure over $\R^{\infty}$ can be explored.

\begin{exa}\label{e8.2.6}%{e2.3.7}
	\begin{itemize}
		\item[(i)] The relationship between $\R^m$ and $\R^n$ ($m\neq n$)?
		
		Assume $s=\gcd(m,n)$ is the great common divisor. Then it is easy to see that
		under the topology ${\cal T}_d$, some points of $\R^m$ and $\R^n$  are glued together, when they are of zero distance. Moreover,
		$$
		d(x,y)=0,\quad x\in \R^m, y\in \R^n,
		$$
		if and only if, there exists $z\in \R^s$ such that $x=z\otimes \J_{m/s}$ and $y=z\otimes \J_{n/s}$. Roughly speaking, $\R^m$ intersects $\R^n$ at a subspace $Z\sim\R^s$. Figure \ref{Fig8.2.3} shows this.
		
		\begin{figure}
			\centering
			\setlength{\unitlength}{8mm}
			\begin{picture}(5.5,4)(-0.5,0)\thicklines
				\put(1,2){\line(1,0){3}}
				\put(2.5,2.2){$Z$}
				\put(0.5,2.6){$\R^n$}
				\put(0.5,0.2){$\R^m$}
				%
				%{\color{red}
					\put(0,0){\line(1,2){2}}
					\put(3,0){\line(1,2){2}}
					\put(0,0){\line(1,0){3}}
					\put(2,4){\line(1,0){3}}
				%}
				%{\color{blue}
					\put(1.9,1){\line(-1,1){2}}
					\put(4.9,1){\line(-1,1){2}}
					\put(1.9,1){\line(1,0){3}}
					\put(-0.1,3){\line(1,0){3}}
				%}
			\end{picture}
			
			\caption{Intersection of $\R^m$ and $\R^n$ \label{Fig8.2.3}}
		\end{figure}
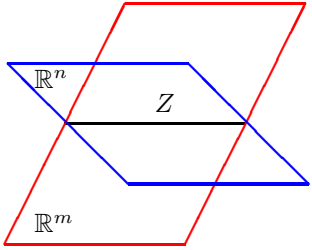

		\item[(ii)] The sphere in $\R^{\infty}$:
		
		Let $c_0\in \R^{\infty}$ and $r>0$. The sphere  in $\R^{\infty}$, centering at $c_0$ with radius $r$, is defined by
		
		\begin{align}\label{8.2.7}
			S_r(c_0):=\left\{x\in \R^{\infty}\;|\; d_{{\cal V}}(x,c_0)=r \right\}.
		\end{align}
		
		Assume $c_0\in \R^m$. Now consider $\R^n$. We project $c_0$ to $\R^n$, say, $\pi^m_n(c_0)=c_n$. Denote
		$h_n:=c_0-c_n$. Assume
		\begin{align}\label{8.2.8}
			\|h_n\|_{{\cal V}}\leq r.
		\end{align}
		Then a sphere can be constructed as
		\begin{align}\label{8.2.9}
			S_n:=\left\{x\in \R^n\;|\; d(x,c_n)=\sqrt{r^2-\|h_n\|^2_{{\cal V}}}\right\}.
		\end{align}
		
		We conclude that the sphere $S_r(c_0)$ in $\R^{\infty}$ consists of all spheres $S_n\in \R^n$, where the distance of $\R^n$ with $c_0$ satisfies (\ref{8.2.8})
		and $S_n$ is constructed by (\ref{8.2.9}).
		Figure \ref{Fig8.2.4} shows this.

		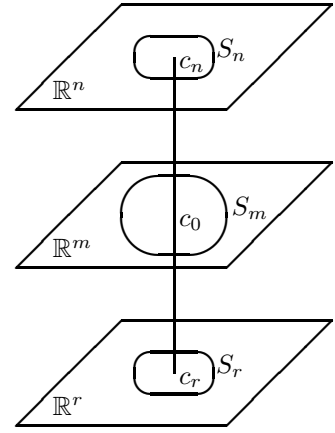
\begin{figure}
			\centering
			\setlength{\unitlength}{7mm}
			\begin{picture}(6,8)\thicklines
				\put(0,0){\line(1,0){4}}
				\put(2,2){\line(1,0){4}}
				\put(0,3){\line(1,0){4}}
				\put(2,5){\line(1,0){4}}
				\put(0,6){\line(1,0){4}}
				\put(2,8){\line(1,0){4}}
				\put(0,0){\line(1,1){2}}
				\put(4,0){\line(1,1){2}}
				\put(0,3){\line(1,1){2}}
				\put(4,3){\line(1,1){2}}
				\put(0,6){\line(1,1){2}}
				\put(4,6){\line(1,1){2}}
				\put(3,1){\line(0,1){6}}
				\put(3,1){\oval(1.5,0.8)}
				\put(3,7){\oval(1.5,0.8)}
				\put(3,4){\oval(2,1.5)}
				\put(0.7,6.2){$\R^n$}
				\put(0.7,3.2){$\R^m$}
				\put(0.7,0.2){$\R^r$}
				\put(3.1,0.8){$c_r$}
				\put(3.1,3.8){$c_0$}
				\put(3.1,6.8){$c_n$}
				\put(3.8,1){$S_r$}
				\put(3.8,7){$S_n$}
				\put(4.1,4){$S_m$}
			\end{picture}
			
			\caption{A Sphere on $\R^{\infty}$ \label{Fig8.2.4}}
		\end{figure}
		
	\end{itemize}
\end{exa}

\section{Differentiable Structure on $\R^{\infty}$}

Since each component space $\R^{n}\subset \R^{\infty}$ has its differentiable structure, which provides a differentiable structure for $\R^{\infty}$.

\begin{dfn}\label{d10.0.1}
	\begin{itemize}
		\item[(i)] A hyper vector space is said to be differentiable, if all its component subspaces are differentiable.
		
		\item[(ii)]
		Let $\Theta\subset \Z_+$, with $(\Theta,\prec)$  be a sublattice of MD-1. Then
		$$
		\R^{\Theta}:=\bigcup_{\theta\in \Theta}\R^{\theta}\subset \R^{\infty},
		$$
		as a sub hyper vector space
		with the inherited topological structure from $\R^{\infty}$ and the differentiable structure from $\R^{\theta}$, where $\theta\in \Theta$, is called an Euclidean hyper manifold.
	\end{itemize}
\end{dfn}

\begin{rem}\label{r10.0.2}
	\begin{itemize}
		\item[(i)] It follows from the definition that the differentiability is component-wise defined. Hence
		$\R^{\infty}$ is obviously differentiable. Then so is $\R^{\Theta}$.
		\item[(ii)] Unlike the differentiability, the topological structure of  an Euclidean hyper manifold is globally defined.
	\end{itemize}
\end{rem}

In the following, the major objects of differential geometry will be defined over Euclidean hyper manifold.

\subsection{Smooth Functions}

\begin{dfn}\label{d10.1.1} A continuous function $h:\R^{\Theta}\ra \R$ is said to be $C^{r}$ (or $C^{\infty}$, or $C^{\omega}$), Denoted the set of  $C^{r}$ (or $C^{\infty}$, or $C^{\omega}$, which means the set of analytic functions.) functions by  $C^{r}(\R^{\Theta})$ (or $C^{\infty}(\R^{\Theta})$, or $C^{\omega}(\R^{\Theta})$),
	if for each $\theta\in \Theta$ the restriction $f|_{\R^{\theta}}$ is $C^{r}$ (or $C^{\infty}$, or $C^{\omega}$).
\end{dfn}

\begin{exa}\label{e10.1.2}
	\begin{itemize}
		\item[(i)] Define $f:\R^{\infty}\ra \R$ by
		$$
		f|_{\R^n}(x)=\dsum_{i=1}^n|x_i|,\quad x\in \R^n,\;n\in \Z_+.
		$$
		Then it is obvious that $f\in C^0(\R^{\infty})$.
		\item[(ii)] Define  $f:\R^{\infty}\ra \R$ by
		$$
		f(x)=\|x\|_{{\cal V}}^2.
		$$
		Then
		$$
		f|_{\R^n}(x)=\frac{1}{n}\dsum_{i=1}^nx^2_i,\quad x\in \R^n,\;n\in \Z_+.
		$$
		It is clear that $f\in C^{\omega}(\R^{\infty})$.
		
	\end{itemize}
\end{exa}

\begin{rem}\label{r10.1.3} A necessary condition for $h$ to be continuous is
	\begin{align}\label{10.1.1}
		h(x)=h(y),\quad x\lra y.
	\end{align}
\end{rem}

An alternative definition for smooth functions is as follows.

\begin{dfn}\label{d10.1.4} A function $h:\R^{\Theta}\ra \R$ is said to be $C^{r}$ (or $C^{\infty}$, or $C^{\omega}$),  if
	\begin{itemize}
		\item[(i)] (\ref{10.1.1}) holds.
		\item[(ii)]
		For each $\theta\in \Theta$ the restriction $f|_{\R^{\theta}}$ is $C^{r}$ (or $C^{\infty}$, or $C^{\omega}$).
	\end{itemize}
\end{dfn}

Since (\ref{10.1.1}) ensures the continuity among different component spaces, the Definition \ref{d10.1.4} is obviously equivalent to Definition \ref{d10.1.1}. Hereafter, for statement ease, we allow $r=\infty$ or $r=\omega$.

In general, such smooth functions are hardly constructed. Particularly, in further discussions (about vector fields etc.) and from application point of view, we may weaken the requirement.

\begin{dfn}\label{d10.1.5} Assume $\Theta$ has a lattice structure and $\Xi\subset \Theta$ is a filter. A function $h:\R^{\Theta}\ra \R$ is said to be $C^{r}$ with respect to $\Xi$, denoted by $h\in C^r_{\Xi}(\Theta)$, if
	\begin{itemize}
		\item[(i)] (\ref{10.1.1}) holds for $x,~y\in \R^{\Xi}=\bigcup_{\xi\in \Xi}\R^{\xi}$.
		\item[(ii)]
		For each $\theta\in \Theta$ the restriction $f|_{\R^{\theta}}$ is $C^{r}$ .
	\end{itemize}
\end{dfn}

We need a lemma for further construction.

\begin{lem}\label{l10.1.6} Let $u,v,w\in \Z_+$, where $u|v$ and $v|w$. Then
	\begin{align}\label{10.1.2}
		\Pi^w_u=\Pi^w_v\Pi^v_u.
	\end{align}
\end{lem}

\noindent{\it Proof.} Set $s=v/u$ and $t=w/v$. Then
$$
\Pi^w_v=\frac{1}{t}(I_v\otimes \J^{^{\rm T}}_t)(I_w\otimes 1)=\frac{1}{t}(I_v\otimes \J^{^{\rm T}}_t).
$$
Similarly,
$$
\Pi^v_u=\frac{1}{s}(I_u\otimes \J^{^{\rm T}}_s).
$$
Then
$$
	\Pi^w_v\Pi^v_u=\frac{1}{st}\left\{ [(I_u\otimes \J_s^{^{\rm T}})I_r]\otimes \J^{^{\rm T}}_t\right\}
	~=\frac{1}{st}(I_u\otimes \J^{^{\rm T}}_{st})=\Pi^w_u.
$$
\hfill $\Box$

The following example shows how to construct a $C_{\Xi}^r$ function from a component space to whole hyper-manifold.

\begin{exa}\label{e10.1.7} Let $h(x)\in C^r(\R^{\theta})$. Define
	$$
	\Xi=\{n\theta\;|\;n\in \Z_+\}.
	$$
	Then it is clear that $(\Xi,\prec)$ is a sublattice of $(\Z_+,\prec)$, and it is a filter.
	Define
	\begin{align}\label{10.1.3}
		H(z):=\begin{cases}
			h(\Pi^m_{\theta}(z)),\quad z\in \R^m,\; m\in \Xi,\\
			0,\quad m\not\in \Xi.
		\end{cases}
	\end{align}
	
	We claim that $H(z)\in C_{\Xi}^r(\R^{\Theta})$.
	Condition (ii) of Definition \ref{d10.1.5} is obviously true. We have only to prove (i).
	Let $x\in \R^{r\theta}$, $y\in \R^{s\theta}$ and $x\lra y$. Denote  $r\wedge s=t$, then there exists
	$z\in \R^{rm T}$ such that $x=z\otimes \J_{r/t}$ and  $y=z\otimes \J_{s/t}$.
	
	Set $\Pi^{tn}_n(z):=z_0$. Using Lemma \ref{l10.1.6}, we have
	$$
		H(x)=h(\Pi^{rn}_nx)=h(\Pi^{tn}_{n}\Pi^{rn}_{tn}(x))=h(\Pi^{tn}_{n}(z))=h(z_0).
	$$
	Similarly,
	$$
	H(y)=h(z_0).
	$$
	 (i) follows.
	
	Note that instead of (\ref{10.1.3}) we can also define
	\begin{align}\label{10.1.4}
		H(z):=
		h(\Pi^m_{\theta}(z)).
	\end{align}
	But we can still only prove (i) for $x,y\in \R^{\Xi}$, where $\Xi=\{kn\;|\;k\in \Z_+\}$.
\end{exa}

\subsection{Vector Fields and Covector Fields}

\begin{dfn}\label{d10.2.1}  Consider $\R^{\Theta}$.
	\begin{itemize}
		\item[(i)] The tangent space of $\R^{\Theta}$ is defined by
		\begin{align}\label{10.2.1}
			T(\R^{\Theta})=\bigcup_{\theta\in \Theta}T(\R^{\theta}).
		\end{align}
		\item[(ii)] A vector field $f\in V^r(\R^{\Theta})$ if it satisfies
		\begin{itemize}
			\item[(a)]
			\begin{align}\label{10.2.2}
				f(x)\lra f(y),\quad x\lra y.
			\end{align}
			\item[(b)]
			\begin{align}\label{10.2.3}
				f|_{\R^{\theta}}\in V^r(\R^{\theta}),\quad \theta\in \Theta.
			\end{align}
		\end{itemize}
	\end{itemize}
\end{dfn}

\begin{exa}\label{e10.2.2} Let $f(x)\in V^r(\R^n)$ and $n\in \Theta$. Then the smallest filter containing $n$ is
	$$
	\Xi:=\{kn\;|\;k\in \Z_+\}.
	$$
	It is easy to extend $f(x)$  to a vector field on $V^r_{\Xi}(\R^{\Theta})$.
	To this end,  we define
	\begin{align}\label{10.2.4}
		F(z):=
		\begin{cases}
			\Pi^n_{m}f(\Pi^{m}_n(z)),\quad m=kn,\\
			0,\quad m\neq kn,\; z\in \R^{m}.\\
		\end{cases}
	\end{align}
	
	It is ready to verify that the condition (ii) of Definition \ref{d10.1.5} is satisfied. We prove condition (i) is also satisfied for $x, y\in \R^{\Xi}$.
	Let $x\in \R^{r\theta}$, $y\in \R^{s\theta}$ and $x\lra y$. Denote $r\wedge s=t$, then there exists
	$z\in \R^{rm T}$ such that $x=z\otimes \J_{r/t}$ and  $y=z\otimes \J_{s/t}$.
	Then
	$$
	\begin{array}{ccl}
		F(x)&=&\Pi^n_{rn}f(\Pi^{tn}_n\Pi^{rn}_{tn}(x))=\Pi^n_{rn}f(\Pi^{tn}_n(z))=\Pi^n_{rn}f(z_0)\\
		~&=&(I_n\otimes \J_r)f(z_0)=f(z_0)\otimes I_r.
	\end{array}
	$$
	Similarly, we have
	$$
	F(y)~=f(z_0)\otimes I_s.
	$$
	Hence
	$$
	F(x)\lra F(y).
	$$
	
\end{exa}

\begin{dfn}\label{d10.2.3}  Consider $\R^{\Theta}$.
	\begin{itemize}
		\item[(i)] The cotangent space of $\R^{\Theta}$ is defined by
		\begin{align}\label{10.2.5}
			T^*(\R^{\Theta})=\bigcup_{\theta\in \Theta}T^*(\R^{\theta}).
		\end{align}
		\item[(ii)] A co-vector field $\omega \in V^{*r}(\R^{\Theta})$ if it satisfies
		\begin{itemize}
			\item[(a)]
			\begin{align}\label{10.2.6}
				\omega(x)\lra \omega(y),\quad x\lra y.
			\end{align}
			\item[(b)]
			\begin{align}\label{10.2.7}
				\omega|_{\R^{\theta}}\in V^{*r}(\R^{\theta}),\quad \theta\in \Theta.
			\end{align}
		\end{itemize}
	\end{itemize}
\end{dfn}

Similarly to Example \ref{e10.2.2}, we are also able to extend a co-vector field $\omega(x)\in V^{*r}(\R^n)$ to $V^{*r}_{\Xi}(\R^{\Theta})$.

\subsection{Distribution and Co-Distribution}

\begin{dfn}\label{d10.3.1}  Consider $\R^{\Theta}$. A distribution (co-distribution) is a rule, which assign for each $x\in \R^{\Theta}$ a
	subspace $D_x\subset T_x(\R^{\Theta})$  ( $D^*_x\subset T^*_x(\R^{\Theta}))$ , satisfying
	\begin{itemize}
		\item[(i)] For each $n\in \Theta$,
		$D|_{\R^n}$ is a distribution over $T(\R^n)$ ($D^*|_{\R^n}$ is a co-distribution over $T^*(\R^n)$ );
		\item[(ii)] Assume $x\in \R^m$, $y\in \R^n$, $m,n\in \Theta$, and $x\lra y$, then $t=m\wedge n\in \Theta$ and there exists $z\in \R^{rm T}$ such that $x=z\otimes \J_{m/t}$ and $y=z\otimes \J_{n/t}$, and then
		\begin{align}\label{10.3.1}
				D_x=D_z\otimes \J_{m/t},\quad D_y=D_z\otimes \J_{n/t}.~
				(
				D^*_x=D^*_z\otimes \J^{^{\rm T}}_{m/t},\quad D^*_y=D^*z\otimes \J^{^{\rm T}}_{n/t}.
				)
		\end{align}

		\item[(iii)] A distribution $D$ (co-distribution $D^*$) is said to be $C^r$ if it is generated by a set of $C^r$ vector fields (co-vector fields) as
		\begin{align}\label{10.3.2}
				D(x)=\dsum_{i=1}^sc_i(x)f_i(x),
				(
				D^*(x)=\dsum_{i=1}^sc_i(x)\omega_i(x),
				)
		\end{align}
		where $f_i(x)\in V^r(\R^{\Theta})$ ($\omega_i(x)\in V^{*r}(\R^{\Theta})$) are $C^r$ vector fields (co-vector fields) and the coefficients $c_i(x)\in C^r(\R^{\Theta})$ are  $C^r$ functions.
		
	\end{itemize}
\end{dfn}

\subsection{Tensor Fields}

\begin{dfn}\label{d10.4.1}  Consider $\R^{\Theta}$. A tensor field with covariant order $r$ and contra-variant order $s$, denoted by
	${\bf t}\in {\bf T}^r_s(\R^{\Theta})$,
	is a rule, which assign for each $x\in \R^{\Theta}$ a set of multi-linear mappings, satisfying
	\begin{itemize}
		\item[(i)] For each $n\in \Theta$,
		for ${\bf t}|_{\R^n}$ is a tensor field
		$$
		{\bf t}|_{\R^n}\in {\bf T}^r_s(\R^n).
		$$
		\item[(ii)] Assume $x\in \R^m$, $y\in \R^n$, $m,n\in \Theta$ and  $x\lra y$, then $t=m\wedge n\in \Theta$ and there exists $z\in \R^{rm T}$ such that $x=z\otimes \J_{m/t}$ and $y=z\otimes \J_{n/t}$, and
		\begin{align}\label{10.4.1}
			\begin{array}{l}
				M^{\bf t}(x)=\underbrace{\E^{^{\rm T}}_{m/t}\otimes\cdots \otimes \E^{^{\rm T}}_{m/t}}_s \otimes %\\
				M^{\bf t}_z\otimes \underbrace{\E_{m/t}\otimes\cdots \otimes \E_{m/t}}_r,\\
				M^{\bf t}(y)=\underbrace{\E^{^{\rm T}}_{n/t}\otimes\cdots \otimes \E^{^{\rm T}}_{n/t}}_s \otimes %\\
				M^{\bf t}_z\otimes \underbrace{\E_{n/t}\otimes\cdots \otimes \E_{n/t}}_r,\\
			\end{array}
		\end{align}
		where $M^{\bf t}(p)$, $p=x,z,y$ are the structure matrix of ${\bf t}$ at $x$, $y$, and $z$ respectively.
		
		\item[(iii)] A tensor field ${\bf t}$ is said to be $C^r$ if the structure matrix of ${\bf t}$ is overall continuous, and within each $n\in \Theta$ , ${\bf t}|_{\R^n}$ is $C^r$.
		
%		\item[(iv)] Let $\Xi$ be a filter of $\Theta$, if the structure matrix of ${\bf t}$ is continuous within $\Xi$ (under the subspace topology), and within each $n\in \Theta$ , ${\bf t}|_{\R^n}$ is $C^r$, denoted by ${\bf t}\in {{\bf T}^r_s}_{\Xi}(\R^{\Theta})$.
	\end{itemize}
\end{dfn}

(We refer to \cite{che23} for detailed discussion about the differential structure over hyper-manifolds.)

\section{Hyper Manifold}

\begin{figure}
	\centering
	\setlength{\unitlength}{4.5mm}
	\begin{picture}(18,8.5)\thicklines
		\put(0,3){\line(1,0){4}}
		\put(0,3){\line(1,1){2}}
		\put(4,3){\line(1,1){2}}
		\put(2,5){\line(1,0){4}}
		\put(12,3){\line(1,0){4}}
		\put(12,3){\line(1,1){2}}
		\put(16,3){\line(1,1){2}}
		\put(14,5){\line(1,0){4}}
		\put(8.5,4){\oval(2.6,1.4)}
		\put(9.5,4){\oval(2.6,1.4)}
		\put(8.5,4){\vector(-1,0){6}}
		\put(9.5,4){\vector(1,0){6}}
		\put(8.5,6.5){\vector(-1,0){6}}
		\put(9.5,6.5){\vector(1,0){6}}
		%\put(2,7.2){$\R^{\Theta}$}
		%\put(16,7.2){$\R^{\Theta}$}
		\put(8.3,6.8){{\tiny $U\cap V$}}
		\put(7.3,0.5){${\bf M}$}
		\put(0.2,0.2){$\R^{\Theta}$}
		\put(12.2,0.2){$\R^{\Theta}$}
		\put(3,6.8){$\Psi$}
		\put(14,6.8){$\Phi$}
		\put(3,3.2){$\Psi_{\theta}$}
		\put(14,3.2){$\Phi_{\theta}$}
		\put(7.5,3.6){{\tiny $U_{\theta}$}}
		\put(9.8,3.6){{\tiny $V_{\theta}$}}
		\put(7.5,6.8){$U$}
		\put(9.8,6.8){$V$}
		\thinlines
		\put(0,0){\line(1,0){4}}
		\put(4,0){\line(1,1){2}}
		\put(0,0){\line(0,1){6}}
		\put(4,0){\line(0,1){6}}
		\put(6,2){\line(0,1){6}}
		\put(0,6){\line(1,0){4}}
		\put(0,6){\line(1,1){2}}
		\put(4,6){\line(1,1){2}}
		\put(2,8){\line(1,0){4}}
		\put(12,0){\line(1,0){4}}
		\put(16,0){\line(1,1){2}}
		\put(12,0){\line(0,1){6}}
		\put(16,0){\line(0,1){6}}
		\put(18,2){\line(0,1){6}}
		\put(12,6){\line(1,0){4}}
		\put(12,6){\line(1,1){2}}
		\put(16,6){\line(1,1){2}}
		\put(14,8){\line(1,0){4}}
		\put(9,1){\oval(4,2)[b]}
		\put(9,7){\oval(4,2)}
		\put(9,4){\oval(4,2)}
		\put(7,1){\line(0,1){6}}
		\put(11,1){\line(0,1){6}}
		\put(8.5,7){\oval(2.6,1.4)}
		\put(9.5,1){\oval(2.6,1.4)[b]}
		\put(8.5,1){\oval(2.6,1.4)[b]}
		\put(8.5,0.3){\line(1,0){1}}
		\put(7.2,0.7){\line(0,1){6}}
		\put(8.2,0.9){\line(0,1){6}}
		\put(9.8,0.9){\line(0,1){6}}
		\put(10.8,0.7){\line(0,1){6}}
		\put(9.5,7){\oval(2.6,1.4)}
	\end{picture}
	
	\caption{Hyper Manifold \label{Fig11.1}}
\end{figure}
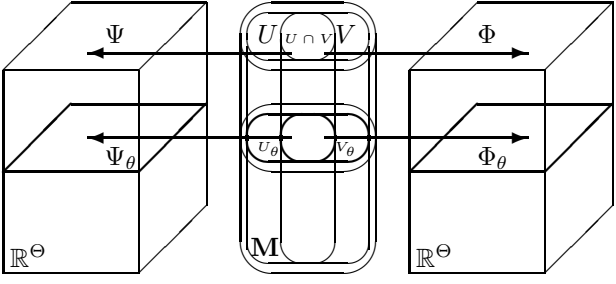

\begin{dfn}\label{d11.1} Let ${\bf M}$ be a second countable topological space with a partition
	$$
	{\bf M}=\bigcup_{\theta\in \Theta}{\bf M}_{\theta},
	$$
	where each ${\bf M}_{\theta}$ is a Hausdorff space (with respect to the inherited topology).
	
	Assume
	$$
	{\cal U}:=\left\{{\cal U}_{\lambda},\;|\;\lambda\in \Lambda \right\}
	$$
	is an open cover of ${\bf M}$, such that for each $U\in {\cal U}$  there exists a mapping $\pi: U\ra U^{\R}$,  where $U^{\R}\subset \R^{\Theta}$ is an open set, and $\pi$ is a homeomorphism.
	
	Denote by $U_{\theta}=U\cap {\bf  M}_{\theta}$, $\theta\in \Theta$. Let $\Psi:U_{\theta}\ra U^{\R}_{\theta}$ with
	$\Psi:U_{\theta}\ra U^{\R}_{\theta}\subset \R^{\theta}$,  and $\Phi:V_{\theta}\ra V^{\R}_{\theta}$ with
	$\Phi:V_{\theta}\ra V^{\R}_{\theta}\subset \R^{\theta}$, and $U_{\theta}\cap V_{\theta}\neq \emptyset$. Then
	\begin{align}\label{11.1}
		\Phi\circ \Psi^{-1}: \Psi(U_{\theta}\cap V_{\theta})\ra \Phi(U_{\theta}\cap V_{\theta})
	\end{align}
	is $C^{r}$.
	
	Then ${\bf M}$ is called a $C^r$ hyper-manifold.
	
	Each ${\cal U}_{\lambda}\in {\cal U}$ is called a coordinate chart.
	
\end{dfn}

Figure \ref{Fig11.1} depicts a hyper-manifold.

All the differentiable objects of ${\bf M}$ can be defined over each coordinate chart $U$ (or $V$) in such a way, that the image of the objects over $\Psi(U)$ are properly defined. Then it is also properly defined on each $U$, and hence on ${\bf M}$.

\begin{exa}\label{e11.2}
	\begin{itemize}
		\item[(i)]
		Consider the set of matrices, ${\cal M}$. Let $\Theta=\{ m\times n\;|\; m,n\in \Z_+\}\subset Z_+$. Define a mapping
		$\Psi: {\cal M}\ra \R^{\Theta}$ as
		\begin{align}\label{11.2}
			\Psi: A\mapsto V_c(A),
		\end{align}
		then $\R^{\Theta}$ can be considered as a global coordinate chart of ${\cal M}$. Hence ${\cal M}$ is a hyper-manifold.
		\item[(ii)]
		Consider $\overline{\cal M}$. Let $\Theta=\{ m\times n+1\;|\; m,n\in \Z_+\}\subset Z_+$. Define a mapping
		$\Psi: \overline{\cal M}\ra \R^{\Theta}$ as
		\begin{align}\label{11.3}
			\Psi: A=aI_{m\times n}+A_0\mapsto (a,V_c(A_0)),
		\end{align}
		then $\R^{\Theta}$ can be considered as a global coordinate chart of $\overline{\cal M}$. Hence $\overline{\cal M}$ is a hyper-manifold.
	\end{itemize}
\end{exa}

\begin{rem}\label{r11.3} The cross-dimensional Euclidean space has give coordinates, but the cross-dimensional manifold does not. The cross-dimensional Euclidean space provides the coordinates for cross-dimensional manifold as local coordinate charts. This relationship is similar to that  between Euclidian space and differential manifold.
\end{rem}

\section{Lie Algebra and Lie Group of Non-square Matrices}

\subsection{Lie Algebra of Non-square Matrices}

\begin{dfn}\label{dn.1.1} \cite{chepr}Consider ${\cal M}_{m\times n}$. The Lie bracket is defined by
	\begin{align}\label{n.1.1}
		[A,B]_{\ttimes}:=A\ttimes B-B\ttimes A.
	\end{align}
\end{dfn}

\begin{prp}\label{pn.1.2} \cite{chepr} $\gl(m\times n, \R):=({\cal M}_{m\times n}, [\cdot,\cdot]_{\ttimes})$ is a Lie algebra.
\end{prp}

Note that when $m=n$, $\gl(m\times n, \R)=\gl(n,\R)$. Hence $\gl(m\times n,\R)$ is a natural generalization of the classical general linear algebra  to non-square matrices.

To construct a Lie group $\GL(m\times n,\R)$ with $\gl(m\times n,\R)$ as its Lie algebra, we need some preparations.

\begin{dfn}\label{dn.1.2}\cite{var84}
	\begin{itemize}
		\item[(i)] Let ${\bf g}=(G,[\cdot,\cdot])$ be a Lie algebra, and $H\subset G$ be a subspace. If
		\begin{align}\label{n.1.1}
			[h_1,h_2]\in H,\quad h_1,h_2\in H,
		\end{align}
		then ${\bf h}=(H, [\cdot,\cdot])$ is called a Lie sub-algebra of ${\bf g}$.
		\item[(ii)] Let ${\bf h}$ be a Lie sub-algebra of ${\bf g}$. If
		\begin{align}\label{n.1.2}
			[h,g]\in H,\quad h\in H,\; g\in G,
		\end{align}
		then ${\bf h}$ is called an ideal of ${\bf g}$.
	\end{itemize}
\end{dfn}

It follows from the definition immediately that
\begin{prp}\label{pn.1.3} Consider $\gl(m\times n,\R)$.
	\begin{itemize}
		\item[(i)] Set
		\begin{align}\label{n.1.201}
			Q_{m\times n}:=\left\{A\in {\cal M}_{m\times n}\;|\; \Tr(\Pi_A)=0\right\},
		\end{align}
		then it is easy to verify that $Q_{m\times n}$ is an ideal of $\gl(m\times n,\R)$.
		\item[(ii)] Set
		\begin{align}\label{n.1.3}
			Z_{m\times n}:=\left\{A\in {\cal M}_{m\times n}\;|\; \Pi_A=rI_m,\;r\in \R\right\},
		\end{align}
		then $Z_{m\times n}$ is an ideal of $\gl(m\times n,\R)$.
	\end{itemize}
\end{prp}

\begin{exa}\label{en.1.4} Consider $\gl(2\times 3,\R)$. Solving
	$$
	A\Psi_{3\times 2}=\lambda I_2
	$$
	yields
	\begin{align}\label{n.1.4}
			Z_{m\times n}=\left\{\left.
			\sqrt{6}\begin{bmatrix}
				1.5&-2&1\\
				1&-2&1.5
			\end{bmatrix}\lambda+
			(1,-2,1)\begin{bmatrix}
				\a\\\b
			\end{bmatrix}\right| \lambda,\a,\b\in \R
			\right\}.
	\end{align}
	
	In general, when $m<n$,  $Z_{m\times n}$ is an $(n-m)m+1$ dimensional subspace.
\end{exa}

\begin{dfn}\label{dn.1.4}\cite{var84}
	\begin{itemize}
		\item[(i)] Let ${\bf h}=(H,[\cdot,\cdot])$ and ${\bf g}=(G,[\cdot,\cdot])$ be two Lie algebras. If there exists a linear mapping $\pi: H\ra G$ such that
		\begin{align}\label{n.1.5}
			\pi [h_1,h_2]=[\pi(h_1),\pi(h_2)],\quad h_1,h_2\in H,
		\end{align}
		then ${\bf h}$ and ${\bf g}$ are homomorphic and $\pi$ is a homomorphism.
		\item[(ii)] $\pi: H\ra G$ is an one-to-one and onto mapping, and both $\pi$ and $\pi^{-1}$ are homomorphisms, then
		${\bf h}$ and ${\bf g}$ are isomorphic and $\pi$ is an isomorphism.
	\end{itemize}
\end{dfn}

\begin{prp}\label{pn.1.5} Consider $\gl(m\times n,\R)$ and $\gl(m_0\times n_0,\R)$  Assume $m_0|m$, $n_0|n$, and denote $m=m_0a$ and $n=n_0b$. Define
	$\pi:\gl(m_0\times n_0,\R)\ra \gl(m\times n,\R)$ by
	\begin{align}\label{n.1.6}
		\pi(A):=A\otimes \E_{a\times b},\quad A\in \gl(m_0\times n_0).
	\end{align}
	Then $\pi$ is a homomorphism.

\end{prp}

\noindent{\it Proof.} We have to show that
$$
[\pi(A_0),\pi(B_0)]=\pi([A_0,B_0]),\quad A_0,B_0\in \gl(m_0\times n_0,\R).
$$
It suffices to show that setting
$A=A_0\otimes \E_{a\times b}$ and $B=B_0\otimes \E_{a\times b}$, then
\begin{align}\label{n.1.7}
	A\ttimes B=\pi(A_0\ttimes B_0),\quad A_0,B_0\in  \gl(m_0\times n_0,\R).
\end{align}

Let $t_0=m_0\vee n_0$ and $\xi=m\vee n$. Then
$$
\begin{array}{l}
	\pi(A_0\ttimes B_0)=\pi\left[(A\otimes \E^{^{\rm T}}_{t_0/n_0})(B\otimes \E_{t_0/m_0})\right]
	=\left[(A\otimes \E^{^{\rm T}}_{t_0/n_0})(B\otimes \E_{t_0/m_0})\right]\otimes \E_{a\times b}\\
	~~=\left[(A\otimes \E^{^{\rm T}}_{t_0/n_0})(B\otimes \E_{t_0/m_0})\right]\otimes (\E_{a}E^{^{\rm T}}_b)
=(A\otimes \E_{a\times t_0/n_0})(B\otimes \E_{t_0/m_0\times b}).\\
\end{array}
$$
In addition,
$$
\begin{array}{l}
	A\ttimes B=(A_0\otimes \E_{a\times b})\ttimes (B_0\otimes \E_{a\times b})
=(A_0\otimes \E_{a\times b}\otimes \E^{^{\rm T}}_{\xi/n})(B_0\otimes \E_{a\times b}\otimes \E_{\xi/m})\\
	~~=(A_0\otimes \E_{a\times b\xi/n})(B_0\otimes \E_{a\xi/m\times b})
	=(A_0\otimes \E_{a\times \xi/n_0})(B_0\otimes \E_{\xi/m_0\times b})
=(A\otimes \E_{a\times t_0/n_0})(B\otimes \E_{t_0/m_0\times b}).\\
\end{array}
$$
The conclusion follows.

\hfill $\Box$

\begin{cor}\label{cn.1.6}
	Assume  $\gl(m\times n,\R)$, $\gl(m_0\times n_0,\R)$ and the mapping $\pi$ are as in Proposition \ref{pn.1.5}. Then
	${\rm Im}(\pi)$ is a Lie sub-algebra of $\gl(m\times n,\R)$.
\end{cor}

\noindent{\it Proof.} Note that $\pi$ is a homomorphism. Moreover, it is obvious that $\pi$ is one to one, then the conclusion is obvious.

\hfill $\Box$

\subsection{Lie Group of Non-square Matrices}

Define $I_{m\times n}$ as a formal identity for the algebra $({\cal M}_{m\times n}\ttimes)$, which satisfies
\begin{align}\label{n.2.1}%{2.4.2}
	\begin{array}{l}
		I_{m\times n}\ttimes I_{m\times n}=I_{m\times n}\\
		I_{m\times n}\ttimes A=A\ttimes I_{m\times n}=A,\quad \forall A\in {\cal M}_{m\times n}.
	\end{array}
\end{align}

Then we have an extended matrix-ring as
\begin{align}\label{n.2.2}
	\overline{\cal M}_{m\times n}:=\{aI_{m\times n}+A_0\;|\; a\in\R, ~A_0\in {\cal M}_{m\times n}\},
\end{align}
with
\begin{align}\label{n.2.3}
	\begin{array}{l}
		(aI_{m\times n}+A_0)+(bI_{m\times n}+B):=(a+b)I_{m\times n}+(A_0+B_0),\\
		(aI_{m\times n}+A_0)\ttimes (bI_{m\times n}+B_0):=(ab)I_{m\times n}+aB_0+bA_0+A_0\ttimes B_0,
		~a,b\in \R,~A,B\in {\cal M}_{m\times n}.
	\end{array}
\end{align}

Note that where $m\neq n$ the objects $\{aI_{m\times n}+A_0\}$ of an extended matrix-ring are not matrices unless $a=0$.

Let $A\in {\cal M}_{m\times n}$. Then $\Pi_A=A\Psi_{n\times m}$, where $\Psi_{n\times m}$ defined by (\ref{2.3.3}), is called the restricted (square) form of $A$, because
$$
A\ttimes x=\Pi_Ax,\quad \forall x\in \R^m.
$$
Define
\begin{align}\label{n.2.4}%{2.4.5}
	A^{\langle k\rangle}=
	\begin{cases}
		I_{m\times n},\quad k=0,\\
		\underbrace{A\ttimes \cdots\ttimes A}_k,\quad k>0.\\
	\end{cases}
\end{align}
In addition, let $p(x)=x^n+c_{n-1}x^{n-1}+\cdots+c_1x+c_0$, we define
\begin{align}\label{n.2.5}%{2.4.6}
		p(\langle A\rangle):=A^{\langle n\rangle}+c_{n-1}A^{\langle n-1\rangle}+\cdots+c_1A
		+c_0I_{m\times n}\in \overline{\cal M}_{m\times n}.
\end{align}

Then we have the following result.

\begin{thm}\label{tn.2.3}%{t2.4.3}
	(Generalized  Cayley-Hamilton theorem)\cite{chepr}
	Let $p(x)=x^m+c_{m-1}x^{m-1}+\cdots+c_1x+c_0$ be the characteristic function of $\Pi_A$, which is also called the characteristic function of $A$. Then
	\begin{align}\label{n.2.6}%{2.4.7}
		A^{\langle m+1\rangle}+c_{m-1}A^{\langle m \rangle}+\cdots+c_1A^{\langle 2\rangle}+c_0A=0.
	\end{align}
\end{thm}

\begin{rem}\label{rn.2.4}%{r2.4.4}
	\begin{itemize}
		\item[(i)] When $n<m$ we can take the characteristic function of $\Pi_{A^{\rm T}}$ as the characteristic function of $A$. This alternative choice can reduce the degree of the characteristic function.
		
		\item[(ii)] Formally, we may rewritten (\ref{n.2.6}) into the form as
		\begin{align}\label{n.2.7}%{2.4.8}
			A^{\langle m\rangle}+c_{m-1}A^{\langle m-1\rangle}+\cdots+c_1A+c_0I^A_{m\times n}=0.
		\end{align}
		Then if $c_0\neq 0$ we can define
		\begin{align}\label{n.2.8}%{2.4.9}
			I^A_{m\times n}:=-\frac{1}{c_0}[A^{\langle m\rangle}+c_{m-1}A^{\langle m-1\rangle}+\cdots+c_1A].
		\end{align}
		According to  Generalized  Cayley-Hamilton theorem, it is clear that
		\begin{align}\label{n.2.9}%{2.4.10}
			A\ttimes I^A_{m\times n}=I^A_{m\times n}\ttimes A=A.
		\end{align}
		Moreover, for any polynomial $p\langle A\rangle$, we also have
		\begin{align}\label{n.2.10}%{2.4.11}
			p(\langle A\rangle)\ttimes I^A_{m\times n}=I^A_{m\times n}\ttimes p(\langle A\rangle)=p(\langle A\rangle).
		\end{align}
		Then $A$ is said to be relatively invertible and define its relative inverse by
		\begin{align}\label{n.2.11}%{2.4.12}
				A^{\langle -1\rangle}= -\frac{1}{c_0} \left[A^{\langle m-1\rangle}+c_{m-1}A^{\langle m-2\rangle}+\cdots+c_1I^A_{m\times n}\right].
		\end{align}
		It follows that
		\begin{align}\label{n.2.12}%{2.4.13}
			A\ttimes A^{\langle -1\rangle}=A^{\langle -1\rangle}\ttimes A=I^A_{m\times n}.
		\end{align}
		It seems that if $A$ is relatively invertible, then $I_{m\times n}|_{p(\langle A\rangle)}$ is $I^A_{m\times n}$. Unfortunately, $I^A_{m\times n}$ depends on $A$. Moreover, it is a matrix while $I_{m\times n}$ is not.
	\end{itemize}
\end{rem}

\begin{dfn}\label{dn.2.5}%{d2.4.5}
	$A\in \overline{\cal M}_{m\times n}$ is invertible, if there exists $B\in \overline{\cal M}_{m\times n}$, such that
	$$
	A\ttimes B=B\ttimes A=I_{m\times n}.
	$$
\end{dfn}

\begin{rem}\label{rn.2.6}%{r2.4.6}
	\begin{itemize}
		\item[(i)] Assume $A=aI_{m\times n}+A_0$, where $A_0\in {\cal M}_{m\times n}$.  It is obvious that a necessary condition for $A$ being invertible is $a\neq 0$.
		\item[(ii)] If $A=aI_{m\times n}+A_0$ is invertible and $A^{\langle -1\rangle}=B$. Then $A'=\frac{1}{a}A$ is also invertible and $(A')^{\langle -1\rangle}=aB$. Hence, to find the inverse of $A$, we can,without loss of generality, assume $a=1$.
		\item[(iii)] Assume $A=I_{m\times n}+A_0$, $A$ is invertible, if and only if, there exists $B=I_{m\times n}+B_0$, where $B_0\in {\cal M}_{m\times n}$ such that
		\begin{align}\label{n.2.13}%{2.4.14}
				A\ttimes B=(I_{m\times n}+A_0)\ttimes (I_{m\times n}+B_0)
				=I_{m\times n}+A_0+B_0+A_0\ttimes B_0.
		\end{align}
		\item[(iv)] Assume $m=n$. Then $\overline{\cal M}_{n\times n}={\cal M}_{m\times n}$  because now $I_{n\times n}=I_n\in {\cal M}_{n\times n}$.  Then it is easily verifiable that $A\in {\cal M}_{n\times n}$ is invertible (under the classical definition), if and only if, (\ref{n.2.13})  holds. Hence Definition \ref{dn.2.5} can be considered as an extension of the classical invertibility to non-square case.
		
		\item[(v)] The $\GL(m\times n,\R)$ has been defined by \cite{chepr} as
		\begin{align}\label{n.2.1301}
				\GL(m\times n,\R):=\left\{A=I_{m\times n}+A_0\;|\; A_0\in {\cal M}_{m\times n}~\mbox{and}~A~\mbox{is invertible} \right\}.
		\end{align}
		Then it was proved in \cite{chepr} that $\GL(m\times n,\R)$ is a Lie group, which has $\gl(m\times n,\R)$ as its Lie algebra.
	\end{itemize}
\end{rem}

Next, we consider when $A=I_{m\times n}+A_0\in \overline{\cal M}$ is invertible.
Assume the characteristic function of $A_0$ is
\begin{align}\label{n.2.14}%{2.4.15}
	p(x)=x^{k-1}+c_{k-1}x^{k-2}+\cdots+c_2x+c_1.
\end{align}
Using generalized  Cayley-Hamilton theorem, we have
$$
A_0^{\langle k\rangle}+c_{k-1}A_0^{\langle k-1\rangle}+\cdots+c_2A_0^{\langle 2\rangle}+c_1A_0=0.
$$
We assume
\begin{align}\label{n.2.15}%{2.4.16}
	B_0=x_{k-1}A_0^{\langle k-1\rangle}+x_{k-2}A_0^{\langle k-2\rangle}+\cdots+x_2A_0^{\langle 2\rangle}+x_1A_0.
\end{align}
Then
$$
\begin{array}{l}
	A_0\ttimes B_0=x_{k-1}A_0^{\langle k\rangle}+x_{k-2}A_0^{\angle k-1\rangle}+\cdots+x_2A_0^{\langle 3\rangle}+x_1A^{\langle 2\rangle}_0\\
	~~=-x_{k-1}(c_{k-1}A_0^{\langle k-1\rangle}+\cdots+c_2A_0^{\langle 2\rangle}+c_1A_0)+x_{k-2}A_0^{\langle k-1\rangle}+\cdots+x_2A_0^{\langle 3\rangle}+x_1A^{\langle 2\rangle}_0\\
\end{array}
$$
Setting
\begin{align}\label{n.2.16}%{2.4.17}
	A_0+B_0+A_0\ttimes B_0=0
\end{align}
and comparing the coefficients yield the following linear system
\begin{align}\label{n.2.17}%{2.4.18}
	\begin{cases}
		-c_{k-1}x_{k-1}+x_{k-2}+x_{k-1}=0,\\
		-c_{k-2}x_{k-1}+x_{k-3}+x_{k-2}=0,\\
		\cdots,\\
		-c_2x_{k-1}+x_1+x_2=0,\\
		-c_1x_{k-1}+x_1+1=0.
	\end{cases}
\end{align}

Solving (\ref{n.2.17}) yields
\begin{align}\label{n.2.18}%{2.4.19}
	\begin{cases}
		x_1=c_1x_{k-1}-1,\\
		x_2=(c_2-c_1)x_{k-1}+1,\\
		x_3=(c_3-c_2+c_1)x_{k-1}-1,\\
		\cdots,\\
		x_{k-2}=(c_{k-2}-c_{k-3}\pm \cdots+(-1)^{k-1}c_1)x_{k-1}+(-1)^{k-2}\\
		(c_{k-1}-c_{k-2}+c_{k-3}\pm\cdots +(-1)^{k}c_1-1)x_{k-1}+(-1)^{k-1}=0.
	\end{cases}
\end{align}

Then we have the following result.

\begin{thm}\label{tn.2.7}%{t2.4.7}
	Consider $A=I_{m\times n}+A_0$ and assume the characteristic function of $A_0$ is (\ref{n.2.14}).
	Then $A$ is invertible, if and only if,
	\begin{align}\label{n.2.19}%{2.4.20}
		c_{k-1}-c_{k-2}+c_{k-3}\pm\cdots +(-1)^{k-2}c_1-1\neq 0.
	\end{align}
\end{thm}

\noindent{\it Proof.} The previous constructing process is a constructive proof for the sufficiency. As for the necessity, it is reasonable to assume the inverse of $A$ is an analytic function of $A_0$. Taking the   generalized  Cayley-Hamilton theorem into consideration, it is obvious that the form of $B_0$, expressed by (\ref{n.2.15}) is a general expression of an analytic function of $A_0$.

\hfill $\Box$

\begin{rem}\label{pn.2.8}%{p2.4.8}
	\begin{itemize}
		\item[(i)] If (\ref{n.2.19}) holds, (\ref{n.2.18}) with (\ref{n.2.15})  provide the inverse of $A$ immediately.
		
		\item[(ii)] It is easy to see that
		\begin{align}\label{n.2.2001}%{2.4.2001}
			I_{m\times n}=\exp(0_{m\times n}).
		\end{align}
		Then by continuity, one sees easily that there exists an open neighborhood ${\cal O}$ of $0_{m\times n}$, say
		$$
		0_{m\times n}\in {\cal O}\subset \R^{mn},
		$$
		such that for each $B\in {\cal O}$ the charactoristic function
		$$
		p_B(x)=x^{k-1}+c_{k-1}x^{k-2}+\cdots+c_2x+c_1,
		$$
		satisfies
		$$
		\left|c_{k-1}-c_{k-2}+c_{k-3}\pm\cdots +(-1)^{k-2}c_1\right|<1.
		$$
		Hence
		$$
		I_{m\times n}+B\in \GL(m\times n,\R),\quad B\in {\cal O}.
		$$
		\item[(iii)] By possibly shrinking ${\cal O}$, we can assume
		\begin{align}\label{n.2.2002}%{2.4.2002}
			\ln{\langle I_{m\times n}+B\rangle}=\dsum_{i=1}^{\infty} (-1)^{i-1}B^{\langle i\rangle},\quad B\in {\cal O}
		\end{align}
		converges. Then we define
		\begin{align}\label{n.2.2003}
			%\begin{array}{l}
			\GL_0(m\times n,\R):=\left\{ I_{m\times n}+B\;|\; B\in {\cal O}\right\}.
			%\end{array}
		\end{align}
		It follows that
		\begin{align}\label{n.2.2004}%{2.4.2003}
			\GL_0(m\times n,\R)\subset \exp({\cal O}),
		\end{align}
		where $\exp({\cal O})$ is an open sub-manifold of $\GL(m\times n,\R)$,
		which shows that $\GL_0(m\times n,\R)$ itself is a sub-Lie group of $\GL(m\times n,\R)$.
	\end{itemize}
\end{rem}

\begin{exa}\label{en.2.9}%{e2.4.9}
	Consider $A=I_{2\times 3}+A_0$, where
	$A_0=(a_{i,j})\in {\cal M}_{2\times 3}$. For statement ease, we let $\Psi_{3\times 2}$ be as in (\ref{2.3.301}). Then we have
	\begin{itemize}
		\item[(i)]
		$$
		%\begin{array}{l}
			\Pi_A=A_0\Psi_{3\times 2}=A_0\begin{bmatrix}
				2&0\\1&1\\0&2\end{bmatrix}%\\
			=
			\begin{bmatrix}
				2a_{1,1}+a_{1,2}&a_{1,2}+2a_{1,3}\\2a_{2,1}+a_{2,2}&a_{2,2}+2a_{2,3}\end{bmatrix}.
		%\end{array}
		$$
		The characteristic function of $A_0$ is
		$$
		p(x)=x^2+c_2x+c1,
		$$
		where
		$$
		\begin{array}{l}
			c_2=-(2a_{1,1}+a_{1,2}+a_{2,2}+2a_{2.3}),\\
			c_1=(2a_{1,1}+a_{1,2})(a_{2,2}+2a_{2,3})
			-(2a_{2,1}+a_{2,2})(a_{1,2}+2a_{1,3}).
		\end{array}
		$$
		We conclude that as long as $c_2-c_1-1\neq 0$, $A$ is invertible. The inverse is
		$$
		B_0=x_2A_0^{\langle 2\rangle}+x_1A_0,
		$$
		where
		$$
		x_1=\frac{c_2-1}{c_1-c_2+1},\quad x_2=\frac{1}{c_1-c_2+1}.
		$$
		\item[(ii)] To verify the above result, we calculate the inverse by definition.
		Assume $B_0=(b_{i,j})\in {\cal M}_{2\times 3}$ such that
		\begin{align}\label{n.2.20}%{2.4.21}
			A_0+B_0+A_0\ttimes B_0=0.
		\end{align}
		A direct calculation shows that
		$$
		A_0\ttimes B_0:=C_0=(c_{i,j})\in {\cal M}_{2\times 3},
		$$
		where
		$$
		\begin{array}{l}
			c_{1,1}=2a_{1,1}b_{1,1}+a_{1,2}b_{1,1}+a_{1,2}b_{2,1}+2a_{1,3}b_{2,1},\\
			c_{1,2}=2a_{1,1}b_{1,2}+a_{1,2}b_{1,2}+a_{1,2}b_{2,2}+2a_{1,3}b_{2,2},\\
			c_{1,3}=2a_{1,1}b_{1,3}+a_{1,2}b_{1,3}+a_{1,2}b_{2,3}+2a_{1,3}b_{2,3},\\
			c_{2,1}=2a_{2,1}b_{1,1}+a_{2,2}b_{1,1}+a_{2,2}b_{2,1}+2a_{2,3}b_{2,1},\\
			c_{2,2}=2a_{2,1}b_{1,2}+a_{2,2}b_{1,2}+a_{2,2}b_{2,2}+2a_{2,3}b_{2,2},\\
			c_{2,3}=2a_{2,1}b_{1,3}+a_{2,2}b_{1,3}+a_{2,2}b_{2,3}+2a_{2,3}b_{2,3}.\\
		\end{array}
		$$
		
		(\ref{n.2.20}) leads to a linear system as
		\begin{align}\label{n.2.21}%{2.4.22}
			\begin{bmatrix}
				\a&0&0&\b&0&0\\
				0&\a&0&0&\b&0\\
				0&0&\a&0&0&\b\\
				\gamma&0&0&\d&0&0\\
				0&\gamma&0&0&\d&0\\
				0&0&\gamma&0&0&\d\\
			\end{bmatrix}
			\begin{bmatrix}
				b_{1,1}\\
				b_{1,2}\\
				b_{1,3}\\
				b_{2,1}\\
				b_{2,2}\\
				b_{2,3}\\
			\end{bmatrix}
			=
			\begin{bmatrix}
				-a_{1,1}\\
				-a_{1,2}\\
				-a_{1,3}\\
				-a_{2,1}\\
				-a_{2,2}\\
				-a_{2,3}\\
			\end{bmatrix},
		\end{align}
		where
		$$
		\begin{array}{l}
			\a=2a_{1,1}+a_{1,2}+1,\\
			\b=a_{1,2}+2a_{1,3},\\
			\gamma=2a_{2,1}+a_{2,2},\\
			\d=2a_{2,3}+a_{2,2}+1.\\
		\end{array}
		$$
		Then it is clear that (\ref{n.2.21}) has unique solution, if and only if,
		$$
		\a\d-\b\gamma\neq 0.
		$$
		A careful calculation shows that
		$$
		\a\d-\b\gamma=-(c_2-c1-1),
		$$
		i.e., the condition obtained in (i) is necessary and sufficient.
		
	\end{itemize}
\end{exa}

In general, the following is obvious.
\begin{prp}\label{pn.2.10}%{p2.4.10}
	Let
	\begin{align}\label{n.2.22}%{2.4.23}%{2.4.17}
		{\cal T}_{m\times n}:=(\{A\in \overline{\cal M}_{m\times n}\;|\; A~\mbox{is invertible}\},\ttimes).
	\end{align}
	Then ${\cal T}_{m\times n}$ is a group.
\end{prp}

\begin{dfn}\label{dn.2.11}%{d2.4.11}
	Let $B\in {\cal M}_{m\times n}$.
	\begin{itemize}
		\item[(i)]
		\begin{align}\label{n.2.23}%{2.4.24}
				\GL(m\times n,\R):=\left\{A=I_{m\times n}+A_0 \;|\; A_0\in {\cal M}_{m\times n} ~\mbox{and}~ A\in {\cal T}_{m\times n}\right\}.
		\end{align}
		called the $m\times n$ general linear group.
		\item[(iI)] Let $B\in {\cal M}_{m\times n}$. Define
		\begin{align}\label{n.2.24}%{2.4.25}
			\exp(\langle B\rangle):=I_{m\times n}+\dsum_{n=1}^{\infty}\frac{1}{n!}B^{\langle n\rangle},
		\end{align}
		Then
		$$
		\exp(\langle B\rangle)\in \GL(m\times n,\R).
		$$
	\end{itemize}
\end{dfn}

\begin{thm}\label{tn.2.12}%{t2.4.12}
	\cite{chepr}
	$\GL(m\times n,\R)$ is a Lie group, which has $\gl(m\times n)$ as its Lie algebra.
\end{thm}

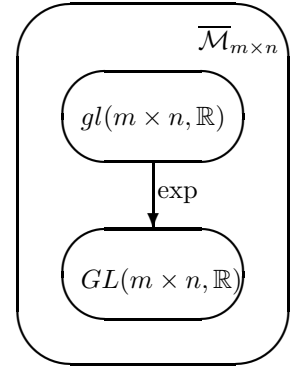
\begin{figure}
	\centering
	\setlength{\unitlength}{6mm}
	\begin{picture}(6,8)\thicklines
		\put(3,4){\oval(6,8)}
		\put(3,2){\oval(4,2)}
		\put(3,5.5){\oval(4,2)}
		\put(3,4.5){\vector(0,-1){1.5}}
		\put(1.4,5.3){$\gl(m\times n,\R)$}
		\put(1.4,1.7){$\GL(m\times n,\R)$}
		\put(3.1,3.7){$\exp$}
		\put(4,7){$\overline{\cal M}_{m\times n}$}
		%\put(1.5,1.8){$i$}
		
	\end{picture}
	
	\caption{$\gl(m\times n,\R)$ and  $\GL(m\times n,\R)$ \label{Fig2.4.1}}
\end{figure}

\begin{rem}\label{rn.2.13}%{r2.4.13}
	\begin{itemize}
		\item[(i)]
		$$
		\GL(m\times n,\R)<{\cal T}_{m\times n}.
		$$
		\item[(ii)]
		Define
		\begin{align}\label{n.2.25}%{2.4.26}
			{\cal B}:=\{\ttimes_{i=1}^k\exp(B_i)\;|\; B_i\in {\cal M}_{m\times n},\;k<\infty\}.
		\end{align}
		Then ${\cal B}$ is a group and
		$$
		{\cal B}< \GL(m\times n,\R)
		$$
		is a sub-Lie group, which is a path-connected component of $\GL(m\times n,\R)$.
		
		\item[(iii)] $A_0$ can be considered as the coordinate of $A=I_{m\times n}+A_0$. It is then clear that
		\begin{align}\label{n.2.26}%{2.4.27}
			\dim(\GL(m\times n,\R))=mn.
		\end{align}
		
		\item[(iv)] $\overline{\cal M}_{m\times n}$ is a ring, which is also an $mn+1$ dimensional vector space. It can also be naturally embedded into $\R^{mn+1}$ by
		$$
		A=aI_{m\times n}+A_0\mapsto (a,A_0^{1,1},A_0^{1,2},\cdots, A_0^{m,n})^{^{\rm T}}\in \R^{mn+1}.
		$$
		Hence $\overline{\cal M}_{m\times n}$ can be naturally considered as a topological vector space, which has the classical topology of Euclidian space $\R^{mn+1}$. Now $\gl(m\times n,\R)$ can be considered as its $mn$ dimensional subspace, and $\GL(m\times n, \R)$ is also its $m\times n$ dimensional sub-manifold. When $m\neq n$ $\gl(m\times n)\bigcap \GL(m\times n)=\emptyset$. This is depicted in Fig.\ref{Fig2.4.1}.
		
		\item[(v)] It is clear that
		$$
		I_{m\times n}=\exp(0_{m\times n}),
		$$
		which is the identity of $\GL(m\times n, \R)$. So it is not an element in $\gl(m\times n,\R)$. When we extend it to be the identity of the ring $\overline{\cal M}$, we have the relations shown in (\ref{n.2.1}).
	\end{itemize}
\end{rem}

\subsection{Relationship with General Linear Group/Algebra}

As an immediate consequence of Adi's theorem, we have the following theorem \cite{var84}.

\begin{thm}\label{tn.3.1} Let ${\bf g}$ be a Lie algebra over $\R$. Then there exists an integer $n\geq 1$ and an analytic subgroup $G$ of $\GL(n,\R)$ such that the Lie algebra of $G$ is isomorphic to ${\bf g}$.
\end{thm}

According to Theorem \ref{tn.3.1}, there is a Lie sub-algebra of ${\bf g}\subset \gl(k,\R)$, such that $\gl(m\times n,\R)$ is  isomorphic to ${\bf g}<\gl(k,\R)$.
This motivates the following result.

\begin{prp}\label{pn.3.2} Consider $\Box: {\cal M}_{m\times n}\ra {\cal M}_{t\times t}$, where $t=m\vee n$. Then
	$\Box: \gl(m\times n,\R)\ra \gl(t,\R)$ is a Lie algebra homomorphism. Moreover, ${\rm Im}(\Box)$ is a sub-algebra of $\gl(t,\R)$.
\end{prp}

\noindent{\it Proof.} Let $A,B\in {\cal M}_{m\times n}$. It suffices to show that
\begin{align}\label{n.3.1}
	\Box(A\ttimes B)=\Box(A)\Box(B),\quad A,B\in {\cal M}_{m\times n}.
\end{align}

$$
\begin{array}{ccl}
	\Box(A\ttimes B)&=&\Box((A\otimes \E^{^{\rm T}}_{t/n})(B\otimes \E_{t/m})
=((A\otimes \E^{^{\rm T}}_{t/n})(B\otimes \E_{t/m})\otimes \E_{t/m\times t/n}\\
	~&=&((A\otimes \E^{^{\rm T}}_{t/n})(B\otimes \E_{t/m})\otimes (\E_{t/m}\E^{^{\rm T}}_{t/n})
=((A\otimes \E^{^{\rm T}}_{t/n}
	\otimes \E_{t/m})(B\otimes \E_{t/m}\otimes \E^{^{\rm T}}_{t/n})\\
	~&=&((A\otimes \E_{t/m\times t/n})(B\otimes \E_{t/m\times t/n})
	~=\Box(A)\Box(B).\\
\end{array}
$$
Hence $\Box: \gl(m\times n,\R)\ra \gl(t,\R)$ is a Lie algebra homomorphism. Now it is obvious that $\Box$ is one-to-one, then
$$
{\rm Im}(\Box)=\Box({\cal M}_{m\times n})
$$
is a Lie sub-algebra of $\gl(t,\R)$.

\hfill $\Box$

Define
\begin{align}\label{n.3.2}
		\GL({\rm Im}({\cal M}_{m\times n}),\R):=\left\{
		\prod_{i=1}^k\exp(B_i)\;|\; B_i\in
		\Box({\cal M}_{m\times n}),\;k<\infty\right\}.
\end{align}
Then
\begin{align}\label{n.3.3}
	\GL({\rm Im}({\cal M}_{m\times n}),\R)< \GL(n,\R)
\end{align}
is a Lie sub-group.

Moreover, if we choose
\begin{align}\label{n.3.4}
		\GL(m\times n,\R):=\left\{
		\ttimes_{i=1}^k\exp(\langle B_i\rangle)\;|\; B_i\in\gl_{m\times n},\R),\;k<\infty
		\right\},
\end{align}
then
\begin{align}\label{n.3.5}
	\GL(m\times n,\R)\cong
	\GL({\rm Im}({\cal M}_{m\times n}),\R)
\end{align}
is a Lie group isomorphism.

\section{From Topological Hyper Group to Hyper Lie Group}

\subsection{Topological Hyper Group}

\begin{dfn}\label{d12.1.1} Assume ${\cal G}=(G,*,{\bf e})$ is both a hyper group and a topological space. If both the product and the inverse are continuous, then ${\cal G}$ is called a topological hyper group.
\end{dfn}

\begin{exa}\label{e12.1.2}
	\begin{itemize}
		\item[(i)] Consider ${\cal T}=\{A\in {\cal M}_1\;|\; A~\mbox{is invertible}\}$. Then it is easy to verify that
		$({\cal T},\ltimes,{\bf e})$ (or $({\cal T},\rtimes,{\bf e})$ ) is a topological  hyper group, where ${\bf e}=\{I_n\;|\;n\in \Z_+\}$.
		
		\item[(ii)] Consider the permutation hyper group ${\bf S}$. By identifying $\sigma\in {\bf S}_n$ with its permutation matrix
		$M_{\sigma}\in {\cal P}_n$ as $\pi_n: {\bf S}_n\ra {\cal P}_n$. Then we can define $\pi:{\bf S}\ra {\cal P}$ by
		\begin{align}\label{12.1.1}
			\pi|_{{\bf S}_n}=\pi_n,\quad n\in \Z_+.
		\end{align}
		It is easy to see that $\pi$ is a bijective mapping.
		Then we define $\circ_{\ell}:{\bf S}\times {\bf S}\ra {\bf S}$ and  $\circ_{r}:{\bf S}\times {\bf S}\ra {\bf S}$ as
		\begin{align}\label{12.1.2}
				\sigma \circ_{\ell}\mu =\pi^{-1}\left[M_{\sigma}\ltimes M_{\mu}\right],\quad
				\sigma \circ_{r}\mu =\pi^{-1}\left[M_{\sigma}\rtimes M_{\mu}\right].
		\end{align}
		
		\item[(iii)] Consider $({\cal T}, \ltimes,{\bf e})$ (or $({\cal T}, \rtimes,{\bf e})$), which is a topological hyper group. As the hyper-subgroup of $({\cal T}, \ltimes,{\bf e})$ (or $({\cal T}, \rtimes,{\bf e})$), the $({\cal P}, \ltimes,{\bf e})$ (or $({\cal P}, \rtimes,{\bf e})$) with its inherited topology is a topological hyper group.
		
		\item[(iv)] Using the topology of  $({\cal P}, \ltimes,{\bf e})$ (or $({\cal P}, \rtimes,{\bf e})$), $({\bf S},\circ_{\ell}, {\bf {\rm Id}})$
		(or $({\bf S},\circ_{r}, {\bf {\rm Id}})$) is a topological hyper group.
		
	\end{itemize}
\end{exa}

\subsection{Hyper Lie Group}

\begin{dfn}\label{d12.2.1} Assume ${\cal G}=(G,*,{\bf e})$ is a topological hyper group and is also an analytic hyper-manifold with both the product and the inverse as analytic mappings. Then ${\cal G}$ is called a hyper Lie group.
\end{dfn}

\begin{exa}\label{e12.2.2} Consider ${\cal G}:=(\R^{\infty},\hat{+},{\bf e})$, where ${\bf e}=\{0_{m\times n}\;|\; m,n\in \Z_+\}$. Equipped by its Euclidean hyper-manifold structure,  ${\cal G}$ is a hyper Lie group.
\end{exa}

\section{Hyper Linear Algebra}
\label{shgla}

\subsection{Matrix Hyper Lie Algebra}

\begin{dfn}\label{d9.1.1} Consider ${\cal M}_0\subset {\cal M}$. Assume there are two binary operators: $+:{\cal M}_0\times {\cal M}_0\ra {\cal M}_0$, $\times:{\cal M}_0\times {\cal M}_0\ra {\cal M}_0$, such that
	$({\cal M}_0,+)$ is a hyper vector space; $({\cal M}_0,+,\times)$ is a hyper ring,
	then $({\cal M}_0, +, \times)$ is called a matrix hyper-algebra.
\end{dfn}

\begin{prp}\label{p9.1.2} Assume $({\cal M}_0, +, \times )$  is a matrix hyper-algebra.
	A Lie bracket is defined by
	\begin{align}\label{9.1.1}
		[A,B]:=A\times B - B\times A.
	\end{align}
	Then it satisfies the following conditions for Lie algebra.
	\begin{itemize}
		\item[(i)] Skew Symmetry:
		\begin{align}\label{9.1.2}
			[A,B]=-[B,A], \quad A,B\in {\cal M}_0.
		\end{align}
		\item[(ii)] Bi-Linearity:
		\begin{align}\label{9.1.3}
			[aA+bB,C]=a[A,C]+b[B,C],\quad A.B.C\in {\cal M}_0.
		\end{align}
		\item[(iii)] Jacobi Identity:
		\begin{align}\label{9.1.4}
			[A,[B,C]]+[B,[C,A]]+[C[A,B]]=0,
		\end{align}
	\end{itemize}
	i.e., a matrix hyper ring with a Lie bracket defined by (\ref{9.1.1})  is called a hyper-Lie algebra.
\end{prp}

\noindent{\it Proof.} (\ref{9.1.2}) comes from definition (\ref{9.1.1}) immediately. The distributivity of hyper ring leads to (\ref{9.1.3}).
(\ref{9.1.4}) follows from a straightforward computation, using the linearity.
\hfill $\Box$

A matrix hyper-algebra over matrices with Lie bracket defined by (\ref{9.1.1}) is called a hyper linear algebra.

Moreover, the following consequence is obvious.

\begin{prp}\label{p9.1.3} Assume $({\cal M}_0, +, \times )$ is a hyper linear algebra, and the identity set of
	$({\cal M}_0, +)$ is ${\bf e}$. Then
	\begin{itemize}
		\item[(i)] each vector space ${\cal M}_0^e$, $e\in {\bf e}$,
		is a general linear algebra;
		\item[(ii)]  if $+$ and $\times$ are consistent with respect to the equivalence, the equivalence classes $({\cal M}_0/\sim, +,\times)$ is a Lie algebra, called the equivalence linear algebra.
	\end{itemize}
\end{prp}

It is well known that $({\cal M}_{n\times n}+,\times)$  with the Lie bracket defined by (\ref{9.1.1}) is the general linear  algebra  $\gl(n,\R)$. So the hyper linear algebra is a natural generalization of the classical general linear algebra.

The above argument shows that as long as a matrix hyper ring exists, it leads to a hyper linear algebra automatically. Hence, we have the following hyper matrix Lie algebras based on the matrix hyper rings\cite{che23}:

\begin{exa}\label{e9.1.4}
	\begin{itemize}
		\item[(i)] Consider $({\cal M}_1,\bar{+},\ltimes)$, defined by equation (\ref{6.1.14}). Define
		\begin{align}\label{9.1.5}
			[A,B]:=A\ltimes B\bar{-}B\ltimes A.
		\end{align}
		Then $({\cal M}_1,\bar{+},\ltimes)$ with the Lie bracket defined by (\ref{9.1.5}) is a hyper-Lie algebra.
		Denote this hyper-Lie algebra by $\gl_s^1(\R)$.
		
		\item[(ii)] Consider $({\cal M}_1,\hat{+},\circ)$, defined by equation (\ref{6.1.15}). Define
		\begin{align}\label{9.1.6}
			[A,B]:=A\circ B\hat{-}B\circ A.
		\end{align}
		Then $({\cal M}_1,\hat{+},\circ)$ with the Lie bracket defined by (\ref{9.1.6})  is a hyper-Lie algebra. Denote this hyper-Lie algebra by
		$\gl_s^2(\R)$.
		%
		%\item[(iii)] Consider $({\cal M}_{\mu},\hat{+},\hat{\ttimes})$, defined in  Remark \ref{r6.1.2} equation (\ref{6.1.16}). Define
		%\begin{align}\label{9.1.7}
		%[A,B]:=A\hat{\ttimes} B\hat{-}B\hat{\ttimes} A.
		%\end{align}
		%Then $({\cal M}_{\mu},\hat{+},\hat{\ttimes})$ with the Lie bracket defined by (\ref{9.1.7}) is a hyper-Lie algebra.
	\end{itemize}
\end{exa}

In the following we construct the most general hyper linear algebra over the set of all matrices.

\subsection{Hyper Linear Algebra $\gl(\R)$}

\begin{dfn}\label{d9.2.1} Consider
	\begin{align}\label{9.2.1}
		\gl(\R):=({\cal M},\hat{+},\hat{\ttimes}).
	\end{align}
	It is already known that $\gl(\R)$ is a hyper ring. Define
	\begin{align}\label{9.2.1}
		[A,B]_{\hat{\ttimes}}:=A\hat{\ttimes} B\hat{-} B\hat{\ttimes}A.
	\end{align}
	Then $\gl(\R)$ becomes a hyper linear algebra.
\end{dfn}

Recall that the identity set is
$$
{\bf e}=\left\{0_{m\times n}\;|\; m,n\in \Z_+\right\}.
$$

First, for each $e=0_{m\times n}\in {\bf e}$ ${\cal M}_e={\cal M}_{m\times n}$. Hence we have
\begin{prp}\label{p9.2.2}
	\begin{itemize}
		\item[(i)]
		\begin{align}\label{9.2.2}
			\gl(m\times n,\R):=( ({\cal M}_{m\times n}, +, \ttimes),[\cdot,\cdot])
		\end{align}
		is a Lie algebra.
		\item[(ii)] Particularly,
		\begin{align}\label{9.2.3}
			\gl(n,\R):=( ({\cal M}_{n\times n}, +, \times),[\cdot,\cdot])
		\end{align}
		is a Lie algebra, which is the well known $n$-th order general linear algebra.
		\item[(iii)]
		\begin{align}\label{9.2.4}
			\gl(\R)=\bigcup_{m=1}^{\infty}\bigcup_{n=1}^{\infty}\gl(m\times n,\R)
		\end{align}
		is a partition of the hyper-general linear algebra. And each Lie-algebra $\gl(m\times n,\R)$ is called its component Lie-algebra.
	\end{itemize}
\end{prp}

Then $A\sim B$ means there exist $0_{a\times b}$ and $0_{c\times d}$ such that
\begin{align}\label{9.2.4}
	A\hat{+}0_{a\times b}=B\hat{+}0_{c\times d}.
\end{align}
Assume $A\in {\cal M}_{m\times n}$, $B\in {\cal M}_{p\times q}$, $m\vee p=r$, and $n\wedge q=s$. By deleting common facts the (\ref{9.2.2}) can be rewritten as
\begin{align}\label{9.2.6}
	A\hat{+}0_{r/m\times s/n}=B\hat{+}0_{r/p\times s/q}.
\end{align}
It leads to
\begin{align}\label{9.2.7}
	A\otimes \E_{r/m\times s/n}=B\otimes \E_{r/p\times s/q}.
\end{align}

\begin{lem}\label{l9.2.3}Consider
	$\gl(\R):=({\cal M},\hat{+},\hat{\ttimes}) $. The equivalence is consistent with $\hat{+}$ and $\hat{\ttimes}$.
\end{lem}
\noindent{\it Proof.}

The consistence of the equivalence with respect to $\hat{+}$ has already been proved in Proposition \ref{p3.3.4}. We prove the consistence of equivalence with respect to $\hat{\ttimes}$.

Assume $A_1\in {\cal M}_{m\times n}$, $A_2\in {\cal M}_{p\times q}$, $B_1\in {\cal M}_{\a\times \b}$, $B_2\in {\cal M}_{\gamma\times \d}$, and
$$
\begin{array}{c}
	m\vee \a=a,~n\vee \b=b,~p\vee \gamma=c,\\
	q\vee \d=d,~m\wedge p=s,~n\wedge q=t,\\
	\a\wedge \gamma=\xi,~\b\wedge \d=\eta,~\xi \vee t=u,\\
	a\vee b=\theta=ux,~c\vee d=\lambda=uy,\\
	a=a_0s,~b-b_0\eta,~c=c_0s,~d=d_0\eta.
\end{array}
$$

Assume $A_1\sim A_2$ and $B_1\sim B_2$, we have to show that
\begin{align}\label{9.2.8}
	A_1\hat{\ttimes} B_1 \sim A_2\hat{\ttimes} B_2.
\end{align}
According to the  condition (\ref{9.2.7}) for matrix equivalence, and
using Lemma \ref{l3.3.6}, there exist $A_0\in {\cal M}_{s\times t}$ and $B_0\in {\cal M}_{\xi\times \eta}$, such that
$$
\begin{array}{cc}
	A_1=A_0\otimes \E_{m/s\times n/t},&A_2=A_0\otimes \E_{p/s\times q/t},\\
	B_1=B_0\otimes \E_{\a/\xi \times \b/\eta},&B_2=B_0\otimes \E_{\gamma/\xi\times \d/\eta}.\\
\end{array}
$$
Then we have
$$
\begin{array}{l}
	A_1\hat{\ttimes}B_1=(A_1\otimes \E{a/m\times b/n})\ttimes (B_1\otimes \E{a/\a\times b/\b})\\
	~=(A_0\otimes \E_{m/s\times n/t}\otimes \E_{a/m\times b/n})
	\ttimes (B_0\otimes \E_{\a/\xi\times \b/\eta}\otimes \E_{a/\a\times b/\b})\\
	~=(A_0\otimes \E_{a/s\times b/t})\ttimes(B_0\otimes \E_{a/\xi\times b/\eta})\\
	~=(A_0\otimes \E_{a/s\times b/t}\otimes \E^{^{\rm T}}_{\theta/b})(B_0\otimes \E_{a/\xi\times b/\eta}\otimes \E_{\theta/a})\\
	~=(A_0\otimes \E_{a/s\times \theta/t})(B_0\otimes \E_{\theta/\xi\times b/\eta})\\
	~=(A_0\otimes \E_{a_0\times u/t})(B_0\otimes \E_{u/\xi\times b_0})\\
	~=[(A_0\otimes \E_{1\times u/t})(B_0\otimes \E_{u/\xi\times 1})]\otimes
	(\E_{a_0}\E^{^{\rm T}}_{b_0}).\\
	~=[(A_0\otimes \E_{1\times u/t})(B_0\otimes \E_{u/\xi\times 1})]\otimes
	(\E_{a_0\times b_0}).\\
\end{array}
$$
Similarly, we also have
$$
A_2\hat{\ttimes}B_2=[(A_0\otimes \E_{1\times u/t})(B_0\otimes \E_{u/\xi\times 1})]\otimes
(\E_{c_0\times d_0}).\\
$$
We conclude that
$$
A_1\hat{\ttimes}B_1\sim A_2\hat{\ttimes}B_2.
$$
\hfill $\Box$

Denote the quotient space by
$$
\Sigma={\cal M}/\sim,
$$
and using Lemma \ref{l9.2.3}, the following result is obvious.
\begin{prp}\label{p9.2.4} The equivalent space
	$\gl_{\Sigma}(\R):=(\Sigma,\hat{+},\hat{\ttimes}) $ is a Lie algebra.
\end{prp}

\begin{dfn}\label{d9.2.5} Assume $V$ is a hyper vector space and there is a Lie bracket $[\cdot,\cdot]:V\times V\ra V$, such that $(V,[\cdot,\cdot])$ is a hyper Lie algebra.
	\begin{itemize}
		\item[(i)] If $N\subset V$ is a hyper vector subspace and $(N,[\cdot,\cdot])$ is also a Lie bracket, then
		$(N,[\cdot,\cdot])$ is called a sub hyper Lie algebra of $(V,[\cdot,\cdot])$.
		\item[(ii)]  If $N\subset V$  is a sub hyper Lie algebra, and
		\begin{align}\label{9.2.9}
			[A,N]\subset N,\quad \forall A\in V,
		\end{align}
		then $N$ is called an ideal of $V$.
	\end{itemize}
\end{dfn}

As an immediate consequence, the following result is straightforward verifiable.

\begin{prp}\label{p9.2.6} $({\cal M}_{\mu}, \hat{+},\hat{\ttimes})$ is a sub hyper Lie algebra of $\gl(\R)$, denoted  by $\gl_{\mu}(\R)$.
\end{prp}

\begin{rem}\label{r9.2.7} Summarizing the matrix hyper Lie algebras we obtained so far, we have the followings:
	\begin{itemize}
		\item[(i)] Consider $({\cal M}_1,\bar{+},\ltimes)$, its corresponding Lie hyper algebra is $\gl_s^1(\R)$ described in Example \ref{e9.1.4}.
		\item[(ii)] Consider $({\cal M}_1,\hat{+},\circ)$, its corresponding hyper Lie  algebra is $\gl_s^2(\R)$ described in Example \ref{e9.1.4}.
		\item[(iii)] Consider $({\cal M},\hat{+},\hat{\ttimes})$, its corresponding hyper Lie  algebra is $\gl(\R)$ described in Definition \ref{d9.2.1}.
		\item[(iv)] Consider $({\cal M}_{\mu},\hat{+},\hat{\ttimes})$, its corresponding hyper Lie  algebra is $\gl_{\mu}(\R)$, which is a sub hyper Lie algebra of $\gl(\R)$.
	\end{itemize}
\end{rem}

The following example presents an ideal of $\gl(\R)$.
\begin{exa}\label{e9.2.8}
	Define
	\begin{align}\label{9.2.10}
		T_{\Box}:=\{A\in {\cal M}\;|\;\Tr(\Box(A))=0\}.
	\end{align}
	Then $T_{\Box}$ is an ideal of $\gl(\R)$.
	
	To see this, assume $A\in {\cal M}_{m\times n}$, $B\in {\cal M}_{p\times q}$, and $m\vee n\vee p\vee q=r$.
	Then a straightforward computation shows that
	$$
		[A,B]_{\hat{\ttimes}}=(A\otimes \E_{r/m\times r/n})(B\otimes \E_{r/p\times r/q})
		-(B\otimes \E_{r/p\times r/q})(A\otimes \E_{r/m\times r/n}).
	$$
	Hence,
	$$
	\Tr([A,B]_{\hat{\ttimes}})=0.
	$$
	The conclusion is obvious.
\end{exa}

\subsection{$M$-Depending Sub Hyper Algebra}

This subsection considers a kind of sub-hyper algebras of $\gl(\R)$, which depend on any $M\in {\cal M}$. Given $M\in {\cal M}$, we define
\begin{align}\label{9.2.11}
	G_M:=\{A\in {\cal M}\;|\;A\hat{\ttimes} M\hat{+} M\hat{\ttimes}A^{^{\rm T}}=0\}.
\end{align}

\begin{thm}\label{t9.2.9} Set
	\begin{align}\label{9.2.12}
		\gl(M,\R):=(G_M,[\cdot,\cdot]_{\hat{\ttimes}}).
	\end{align}
	Then $\gl(M,\R)$ is a sub hyper Lie algebra of $\gl(\R)$.
\end{thm}

To prove this theorem, we need some preparations.

\begin{lem}\label{l9.2.10}
	\begin{align}\label{9.2.13}
		A\hat{\ttimes} M\hat{+} M\hat{\ttimes}A^{^{\rm T}}=0,
	\end{align}
	if and only if,
	\begin{align}\label{9.2.14}
		\Box(A\hat{\ttimes} M)=-\Box(M\hat{\ttimes} A^{^{\rm T}}).
	\end{align}
\end{lem}

\noindent{\it Proof.}
Assume $M\in {\cal M}_{m\times n}$, $A\in {\cal M}_{p\times q}$, and
$$
\begin{array}{c}
	m\vee p=\a,~n\vee q=\b,\\
	m\vee q=a,~n\vee p = b,\\
	\a\vee a=m\vee p\vee q=c,~\b\vee b=n\vee p\vee q=d,\\
	\a\vee \b=a\vee b=m\vee n\vee p\vee q=\xi.
\end{array}
$$
Then
$$
	A\hat{\ttimes} M=(A\otimes \E_{\a/p\times \b/q})\ttimes
	(M\otimes \E_{\a/m\times \b/n})
	=(A\otimes \E_{\a/p\times \xi/q})(M\otimes \E_{\xi/m\times \b/n})\in {\cal M}_{\a\times \b}.
$$

$$
	M\hat{\ttimes} A^{^{\rm T}}=(M\otimes \E_{a/m\times b/n})\ttimes
	(A^{^{\rm T}}\otimes \E_{a/q\times b/p})
	=(M\otimes \E_{a/m\times \xi/n})
	(A^{^{\rm T}}\otimes \E_{\xi/q\times b/p})\in {\cal M}_{a\times b}.
$$

$$
\begin{array}{l}
	A\hat{\ttimes} M\hat{+} M\hat{\ttimes} A^{^{\rm T}}=
	(A \otimes \E_{\a/p\times \xi/q})
	(M \otimes \E_{\xi/m\times \b/n})\hat{+}
	~(M\otimes \E_{a/m\times \xi/n})
	(A^{^{\rm T}} \otimes \E_{\xi/q\times b/p})\\
	~=\left[(A \otimes \E_{\a/p\times \xi/q})
	(M \otimes \E_{\xi/m\times \b/n})\right]\otimes \E_{c/\a\times d/\b}+
	~\left[(M\otimes \E_{a/m\times \xi/n})
	(A^{^{\rm T}} \otimes \E_{\xi/q\times b/p})\right]\otimes \E_{c/a\times d/b}\\
	~=(A \otimes \E_{c/p\times \xi/q})
	(M \otimes \E_{\xi/m\times d/n})+
	~(M\otimes \E_{c/m\times \xi/n})
	(A^{^{\rm T}} \otimes \E_{\xi/q\times d/p})\\
\end{array}
$$
Plugging into (\ref{9.2.13}) yields
\begin{align}\label{9.2.15}
		(A \otimes \E_{c/p\times \xi/q})
		(M \otimes \E_{\xi/m\times d/n})+
		~(M\otimes \E_{c/m\times \xi/n})
		(A^{^{\rm T}} \otimes \E_{\xi/q\times d/p})=0_{c\times d}.
\end{align}
It is clear that (\ref{9.2.15}) is equivalent to  (\ref{9.2.13}).

Right $\otimes$-multiplying both sides of  (\ref{9.2.15}) by $\E_{\xi/c\times \xi/d}$, we have
\begin{align}\label{9.2.16}
	\begin{array}{l}
		\left[(A \otimes \E_{c/p\times \xi/q})
		(M \otimes \E_{\xi/m\times d/n})\right]\otimes \E_{\xi/c\times \xi/d}+
		~\left[(M\otimes \E_{c/m\times \xi/n})
		(A^{^{\rm T}} \otimes \E_{\xi/q\times d/p})\right]\otimes \E_{\xi/c\times \xi/d}\\
		~=(A \otimes \E_{\xi/p\times \xi/q})
		(M \otimes \E_{\xi/m\times \xi/n})+
		~(M\otimes \E_{\xi/m\times \xi/n})
		(A^{^{\rm T}} \otimes \E_{\xi/q\times xi/p})\\
		~=\Box(A \hat{\ttimes} M)+
		\Box(M \otimes A^{^{\rm T}})=0_{\xi\times \xi}.
	\end{array}
\end{align}

Conversely, by deleting  the factor $\E_{\xi/c\times \xi/d}$ from both sides of (\ref{9.2.16}), (which is obviously ``legal" as long as $\E_{\xi/c\times \xi/d}\neq 0$), we have (\ref{9.2.15}).
\hfill $\Box$

\begin{lem}\label{l9.2.11}
	\begin{align}\label{9.2.17}
		\Box(A\hat{\ttimes} B)=\Box{A}\circ \Box(B),
	\end{align}
	where $\circ$ is defined by (\ref{2.2.2}).
\end{lem}

\noindent{\it Proof.}
Assume $A\in {\cal M}_{m\times n}$, $B\in {\cal M}_{p\times q}$, and
$$
\begin{array}{c}
	m\vee p=\a,~n\vee q=\b,\\
	m\vee n=a,~p\vee q = b,\\
	\a\vee \b=a\vee b=m\vee n\vee p\vee q=\xi.
\end{array}
$$
Then
$$
\begin{array}{l}
	\Box(A\hat{\ttimes} B)=\Box\left[(A\otimes \E_{\a/m\times \b/n})
	\ttimes   (B\otimes \E_{\a/p\times \b/q}) \right]\\
	~=\Box\left[ (A\otimes \E_{\a/m\times \xi/n})
	\ttimes   (B\otimes \E_{\xi/p\times \b/q})\right]\\
	~=\left[ (A\otimes \E_{\a/m\times \xi/n})
	\ttimes   (B\otimes \E_{\xi/p\times \b/q})\right]\otimes \E_{\xi/\a\times \xi/\b}\\
	~=(A\otimes \E_{\xi/m\times \xi/n})
	(B\otimes \E_{\xi/p\times \xi/q})\\
	~=\left[(A\otimes \E_{a/m\times a/n})\otimes J_{\xi/a}\right]
	\left[(B\otimes \E_{b/p\times b/q})\otimes J_{\xi/b}\right]\\
	~=\Box(A)\circ \Box(B).
\end{array}
$$
\hfill $\Box$

As an immediate consequence, we have
\begin{align}\label{9.2.18}
	\Box({\hat{\ttimes}}_{i=1}^nA_i)=\circ_{i=1}^n \Box(A_i),\quad A_i\in {\cal M}.
\end{align}

\noindent{\it Proof of Theorem \ref{t9.2.9}}
Assume $A,B\in G_M$, using Lemma \ref{l9.2.10}, (\ref{9.2.13}) is equivalent to
\begin{align}\label{9.2.19}
	\Box \left[([A,B]_{\hat{\ttimes}}\hat{\ttimes}M)\right]+\Box \left[M\hat{\ttimes}[A,B]^{\rm T}_{\hat{\ttimes}}\right] =0.
\end{align}
Using Lemma \ref{l9.2.11}, we have
$$
\begin{array}{l}
	\Box \left[([A,B]_{\hat{\ttimes}}\hat{\ttimes}M)\right]+\Box \left[M\hat{\ttimes}[A,B]^{\rm T}_{\hat{\ttimes}}\right]\\
	~=\Box(A)\circ \Box(B)\circ \Box(M)-\Box(B)\circ \Box(A)\circ \Box(M)\\
	~+\Box(M)\circ \Box(B^{\rm T})\circ \Box(A^{\rm T})-\Box(M)\circ \Box(A^{\rm T})\circ \Box(B^{\rm T})\\
	~=-\Box(A)\circ \Box(M)\circ \Box(B^{\rm T})+\Box(B)\circ \Box(M)\circ \Box(A^{\rm T})\\
	~+\Box(M)\circ \Box(B^{\rm T})\circ \Box(A^{\rm T})-\Box(M)\circ \Box(A^{\rm T})\circ \Box(B^{\rm T})\\
	~=\Box(M)\circ \Box(A^{\rm T})\circ \Box(B^{\rm T})-\Box(M)\circ \Box(B^{\rm T})\circ \Box(A^{\rm T})\\
	~+\Box(M)\circ \Box(B^{\rm T})\circ \Box(A^{\rm T})-\Box(M)\circ \Box(A^{\rm T})\circ \Box(B^{\rm T})\\
	=0.
\end{array}
$$
\hfill $\Box$

\section{Hyper Linear Group As Hyper Lie Group}

\subsection{From Matrix Hyper Rings to Matrix Hyper Group}

Observing that we have two well defined Lie group of matrices. One is well known general linear group $\GL(n,\R)$ and the other one is $\GL(m\times n,\R)$, which has been discussed in Section XII. It seems that the way to construct them can be described as follows:
Starting  from matrix ring, we first construct the Lie algebra by defining Lie bracket in a natural way, then build their corresponding group called the Lie group.  Roughly speaking, we have
$$
\begin{array}{l}
	({\cal M}_{n\times n},+,\times)\Rightarrow \gl(n,\R) \Rightarrow \GL(n,\R);\\
	({\cal M}_{m\times n},+,\ttimes)\Rightarrow \gl(m\times n,\R) \Rightarrow \GL(m\times n,\R).\\
\end{array}
$$

Recall Remark \ref{r9.2.7}, we have defined 4 kinds of matrix hyper-Lie algebras.
This section is aimed to build their corresponding hyper-Lie groups.

Following Remark \ref{r9.2.7}, we consider all 4 different matrix hyper-Lie algebras.
\begin{enumerate}
	\item[(1)]
	
	Consider case (i), which  has been discussed in \cite{che23}. But we reconsider them from a new point of view: Consider them as hyper groups and then pose an $\R^{\infty}$  differentiable structure on it.
	
	Let ${\cal T}=\bigcup_{n=1}^{\infty}{\cal T}_n$, where ${\cal T}_n$ is the set of $n\times n$ nonsingular matrices.
	Then
	\begin{align}\label{lg.1.1}%{10.1.1}
		\GL_s(\R):=({\cal T},\ltimes,{\bf e}),
	\end{align}
	where ${\bf e}=\{I_n\;|\; n\in \Z_+\}$.
	Then the  following facts are obvious;
	\begin{itemize}
		\item $\GL_s(\R)$ is a hyper group.
		
		\item Embedding $M_{n\times n}$ into $\R^{n^2}$ by $A\mapsto V_c(A)$, and using the topology ${\cal T}_d$, then $\GL_s(\R)$ becomes a topological space and a topological hyper group, with the inherited topology from $\R^{\infty}$.
		
		\item Its differentiable structure comes from $\R^{\Xi}$ with
		$$
			\Xi=\{n=m^2\in \Z_+,~\mbox{and}
			~n_1=s^2\prec n_2=t^2 ~\mbox{if and only if}~s\prec t\}.
		$$
		
	\end{itemize}
\end{enumerate}

\begin{enumerate}
	\item[(2)]
\end{enumerate}
The topological structure in (i) can hardly be used for this case, because $A\otimes \J_{k\times k}$ is no longer invertible, even if $A$ is. We are not interested in constructing a hyper-Lie group with its hyper-Lie algebra $\gl_s^2(\R)$.

\begin{enumerate}
	\item[(3)]
\end{enumerate}

Denote by
$$
	G_{m\times n}:=(\{I_{m\times n}+A_0\;|\;A_0\in {\cal M}_{m\times n},~\mbox{and}~I_{m\times n}+A_0~\mbox{is invertible.}\}, \hat{\ttimes}).
$$
We define
\begin{align}\label{lg.1.2}%{10.1.2}
	G_{M}:=\bigcup_{m=1}^{\infty}\bigcup_{n=1}^{\infty} G_{m\times n}\subset \overline{\cal M}.
\end{align}

First we need to show that $(G_M,\hat{\ttimes})$ is a hyper group.
Set
\begin{align}\label{lg.1.3}%{10.1.3}
	{\bf e}:=\left\{I_{m\times n}\;|\; (m,n)\in \Z_+\times \Z_+\right\},
\end{align}
with the lattice  MD-2. We need the following definition.

\begin{dfn}\label{dlg.1.1}%{d10.1.1}
	\begin{itemize}
		\item[(i)] The product within ${\bf e}$:
		\begin{align}\label{lg.1.4}%{10.1.4}
			I_{m\times n}\hat{\ttimes} I_{p\times q}:=I_{(m\vee p)\times (n\vee q)}.
		\end{align}
		\item[(ii)] The product of ${\bf e}$ with ${\cal M}$: Let $A_0\in {\cal M}_{p\times q}$, $m\vee p=s$, and $n\vee q=t$. Then
		\begin{align}\label{lg.1.5}%{10.1.5}
			I_{m\times n}\hat{\ttimes} A_0=A_0\hat{\ttimes} I_{m\times n}:=A_0\otimes \E_{s/m\times t/n}.
		\end{align}
	\end{itemize}
\end{dfn}

\begin{dfn}\label{dlg.1.2}%{d10.1.2}
	Let $A,~B\in G_{M}$, say, $A=I_{m\times n}+A_0$, $B=I_{p\times q}+B_0$, where $A_0\in {\cal M}_{m\times n}$, $B_0\in {\cal M}_{p\times q}$, $m\vee p=s$, and $n\vee q=t$.   Then
	\begin{align}\label{lg.1.6}%{10.1.6}
		\begin{array}{l}
			A\hat{\ttimes} B=I_{s\times t} +A_0\otimes \E_{s/m\times t/n}+B_0\otimes \E_{s/p\times t/q}
			+A_0\hat{\ttimes} B_0\in {\cal M}_{s\times t}.
		\end{array}
	\end{align}
\end{dfn}

\begin{prp}\label{p10.1.3} $(G_{M},\hat{\ttimes},{\bf e})$ is a hyper group, where $\hat{\ttimes}$ is defined by (\ref{lg.1.2})  and ${\bf e}$ is defined by (\ref{lg.1.3}).
\end{prp}

\noindent{\it Proof.}

First, we prove the associativity. By the definition of $\hat{\ttimes}$, It is enough to show that
\begin{align}\label{lg.1.7}%{10.1.7}
	(A\hat{\ttimes} B)\hat{\ttimes} C=A\hat{\ttimes}(B\hat{\ttimes} C),\quad A,B,C\in {\cal M}\bigcup {\bf e}.
\end{align}
Consider the following several cases.
\begin{itemize}
	\item[(i)]
	All $A,B,C\in {\cal M}$:
	
	Then (\ref{lg.1.7}) has been proved in Proposition \ref{p2.3.10}.
	
	\item[(ii)] There are one or two factors from ${\bf e}$.
	
	We prove one case, where $A,C\in {\cal M}$ and $B\in {\bf e}$. The proofs for all other cases are similar.
	
	Assume $A\in {\cal M}_{m\times n}$, $B=I_{p\times q}$, $C\in {\cal M}_{r\times s}$, and
	$$
	\begin{array}{l}
		m\vee p=a,~n\vee q=b,~a\vee r=c,\\
		b\vee s=d,~c\vee d=e,~n\vee r=\xi,\\
		p\vee r=u, ~q\vee s=v,~m\vee u=m\vee p\vee r=c,\\
		n\vee v=n\vee q\vee s=d.
	\end{array}
	$$
	Then
	$$
	\begin{array}{l}
		(A\hat{\ttimes} I_{p\times q})\hat{\ttimes} C\\
		~=(A \otimes \E_{a/m\times b/n})\hat{\ttimes} C\\
		~=(A  \otimes \E_{a/m\times b/n}\otimes \E_{c/a\times d/b})\ttimes(C\otimes \E_{c/r\times d/s})\\
		~=(A  \otimes \E_{c/m\times d/n})\ttimes(C\otimes \E_{c/r\times d/s})\\
		~=(A  \otimes \E_{c/m\times d/n}\otimes\E^{^{\rm T}}_{e/d})(C\otimes \E_{c/r\times d/s}\otimes\E_{e/c})\\
		~=(A  \otimes \E_{c/m\times e/n})(C\otimes \E_{e/r\times d/s})\\
		~=(A  \otimes \E_{c/m\times \xi/n})(C\otimes \E_{\xi/r\times d/s}).\\
	\end{array}
	$$
	Similarly, we have
	$$
	\begin{array}{l}
		A\hat{\ttimes} (I_{p\times q}\hat{\ttimes} C)\\
		A\hat{\ttimes} (C\otimes \E_{u/r\times v/s})\\
		~=(A \otimes \E_{c/m\times d/n})\ttimes (C\otimes \E_{u/r\times v/s}\otimes \E_{c/u\times d/v})\\
		~=(A  \otimes \E_{c/m\times d/n}\otimes \E^{^{\rm T}}_{e/d})(C\otimes \E_{c/r\times d/s}\otimes \E_{e/c})\\
		~=(A  \otimes \E_{c/m\times e/n})(C\otimes \E_{e/r\times d/s})\\
		~=(A  \otimes \E_{c/m\times \xi/n})(C\otimes \E_{\xi/r\times d/s}).\\
	\end{array}
	$$
	\item[(iii)] All $A,B,C\in {\bf e}$:
	Let $A=I_{m\times n}$, $B=I_{p\times q}$, and $C=I_{r\times s}$. Then a straightforward computation shows that
	$(A\hat{\ttimes} B)\hat{\ttimes} C=A\hat{\ttimes} (B\hat{\ttimes} C)=I_{(m\vee p\vee r)\times (n\vee q\vee s)}$.
\end{itemize}
The associativity is proved.

Next, we have to show that if $A\in {\cal M}_{m\times n}$ and $B\in {\cal M}_{p\times q}$ are invertible, then so is $A\hat{\ttimes} B$. Using associativity, we have
$$
\begin{array}{l}
	(A\hat{\ttimes} B)\hat{\ttimes} (B^{-1}\hat{\ttimes} A^{-1})=A\hat{\ttimes} (B \hat{\ttimes} B^{-1})\hat{\ttimes} A^{-1}
	=A\hat{\ttimes} (I_{p\times q})\hat{\ttimes} A^{-1}\\
	~~=I_{p\times q}\hat{\ttimes}(A\hat{\ttimes} A^{-1})
	=I_{p\times q}\hat{\ttimes} I_{m\times n})=I_{(m\vee p)\times (n\vee q)}.
\end{array}
$$
Similarly, we have
$$
(B^{-1}\hat{\ttimes} A^{-1})\hat{\ttimes} (A\hat{\ttimes} B)=I_{(m\vee p)\times (n\vee q)}.
$$
Hence, $A\hat{\ttimes}B$ is also invertible.

Finally, for each $e=I_{m\times n}$
$$
(G_{M})_{m\times n}=G_{m\times n},
$$
which is a group. We conclude that $(G_{M},\hat{\ttimes})$ is a hyper group.
\hfill $\Box$

\begin{rem}\label{rlg.1.4}
	\begin{itemize}
		\item[(i)] The invertibility of a matrix is the most important concept in group structure of $G_{M}$.  We emphasize that the Definition \ref{dn.2.5} is consistent with respect to classical one for square matrix. Let $A\in {\cal M}_{n\times n}$. Then we can express
		$A=I_{n\times n}+A_0$,
		where $A_0=A-I_n$. Then it is easy to see that $A$ is invertible in classical sense, i.e., there exists $B=A^{-1}$ such that $AB=BA=I_n$, if and only if,
		$B=I_{n\times n}+B_0$,
		where $B_0=B-I_n$ such that
		$A\hat{\ttimes} B=I_{n\times n}+A_0+B_0+A_0B_0=I_{n\times n}$.
		This shows that $I_{n\times n}=I_n$ and the two invertibilities are equivalent.
		
		\item[(ii)] Assume $A=I_{m\times n}+A_0 \in \overline{\cal M}_{m\times n}$ is invertible, where $m\neq n$. Then $A^{^{\rm T}}= I_{n\times m}+A^{^{\rm T}}_0$ is also invertible. $A\hat{\ttimes} A^{^{\rm T}}$ is also invertible. Intuitively, this seems unusual. This fact shows that the product ``$\hat{\ttimes}$" is different from the classical matrix product, and hence when $m\neq n$ the invertibility is different from the classical one.
	\end{itemize}
\end{rem}

\subsection{From Matrix Hyper Group to Hyper Linear Group}

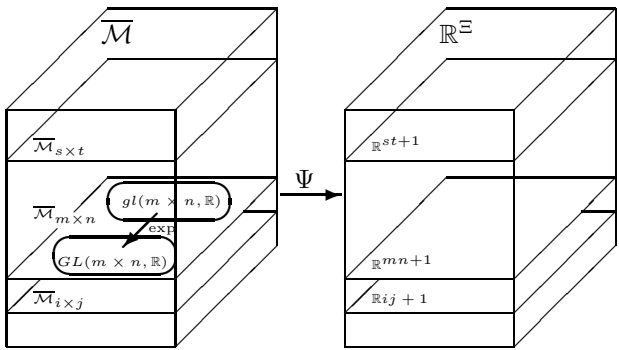
\begin{figure}
	\centering
	\setlength{\unitlength}{4.5mm}
	\begin{picture}(18,10)\thicklines
		\put(8.1,4.5){\vector(1,0){1.8}}
		\put(3.4,4.2){{\tiny$\gl(m\times n,\R)$}}
		\put(1.5,2.4){{\tiny$\GL(m\times n,\R)$}}
		\put(4.2,3.4){{\tiny$\exp$}}
		\put(0.8,3.8){{\tiny $\overline{\cal M}_{m\times n}$}}
		\put(0.8,5.8){{\tiny $\overline{\cal M}_{s\times t}$}}
		\put(0.8,1.3){{\tiny $\overline{\cal M}_{i\times j}$}}
		\put(10.8,5.8){{\tiny $\R^{st+1}$}}
		\put(10.8,2.3){{\tiny $\R^{mn+1}$}}
		\put(10.8,1.3){{\tiny $\R{ij+1}$}}
		\put(8.5,4.7){$\Psi$}
		\put(2.8,9){$\overline{\cal M}$}
		\put(12.8,9){$\R^{\Xi}$}
		\thinlines
		\put(0,0){\line(1,0){5}}
		\put(5,0){\line(1,1){3}}
		\put(0,0){\line(0,1){7}}
		\put(5,0){\line(0,1){7}}
		\put(8,3){\line(0,1){7}}
		\put(0,7){\line(1,0){5}}
		\put(0,7){\line(1,1){3}}
		\put(5,7){\line(1,1){3}}
		\put(3,10){\line(1,0){5}}
		\put(10,0){\line(1,0){5}}
		\put(15,0){\line(1,1){3}}
		\put(10,0){\line(0,1){7}}
		\put(15,0){\line(0,1){7}}
		\put(18,3){\line(0,1){7}}
		\put(10,7){\line(1,0){5}}
		\put(10,7){\line(1,1){3}}
		\put(15,7){\line(1,1){3}}
		\put(13,10){\line(1,0){5}}
		\put(0,1){\line(1,0){5}}
		\put(0,1){\line(1,1){1}}
		\put(7,4){\line(1,0){1}}
		\put(5,1){\line(1,1){3}}
		\put(0,2){\line(1,0){5}}
		\put(0,2){\line(1,1){1.5}}
		\put(2,4){\line(1,1){1}}
		\put(3,5){\line(1,0){5}}
		\put(5,2){\line(1,1){3}}
		\put(0,5.5){\line(1,0){5}}
		\put(0,5.5){\line(1,1){3}}
		\put(3,8.5){\line(1,0){5}}
		\put(5,5.5){\line(1,1){3}}
		\put(10,1){\line(1,0){5}}
		\put(10,1){\line(1,1){1}}
		\put(17,4){\line(1,0){1}}
		\put(15,1){\line(1,1){3}}
		\put(10,2){\line(1,0){5}}
		\put(10,2){\line(1,1){3}}
		\put(13,5){\line(1,0){5}}
		\put(15,2){\line(1,1){3}}
		\put(10,5.5){\line(1,0){5}}
		\put(10,5.5){\line(1,1){3}}
		\put(13,8.5){\line(1,0){5}}
		\put(15,5.5){\line(1,1){3}}
		\thicklines
		\put(4.8,4.3){\oval(3.6,1.1)}
		\put(3.2,2.7){\oval(3.6,1.1)}
		\put(4.45,3.95){\vector(-1,-1){1}}
	\end{picture}
	
	\caption{Matrix Hyper-Lie Group \label{Fig10.1}}
\end{figure}

To put topological and differentiable structure on $G_{M}$ we have to give it coordinates.
To this end, we first construct an Euclidian hyper manifold, $\R^{\Xi}$, where
$$
\Xi=\{m\times n+1,\;|\;  (m,n) \in \Z_+\times \Z_+\},
$$
and $\{(m,n)\in \Z_+\times \Z_+\}$ is the MD-2 lattice.

Then we define a mapping
$\Psi:\overline{\cal M}\ra \R^{\Xi}$ by
$$
aI_{m\times n}+A_0\mapsto (a,V_c(A_0))\in \R^{mn+1},\quad (m,n)\in \Z_+\times \Z_+.
$$
This mapping is described in Figure \ref{Fig10.1}. Then for each $m\times n+1\in \Xi$ we have
$$
\Psi|_{\overline{\cal M}_{m\times n}}:\overline{\cal M}_{m\times n}\ra \R^{mn+1}.
$$
And $V_c(A_0)$ is taken as the coordinate frame of $\GL(m\times n,\R)$.
This assignment makes each $\GL(m\times n,\R)$ a manifold of dimension $mn$. As a consequence, $\GL(\R)$ becomes a matrix hyper Lie group, called the hyper linear group.

Then each slice $\GL(m\times n,\R)$ has $\gl(m\times n.\R)$ as its Lie algebra.

\section{Concluding Remarks}

\subsection{A Mathematics over MDMVs}

To investigate linear (control) systems the classical matrix theory plays a fundamental rule. While the modern  mathematics provides a solid foundation for the matrix theory. Particularly, the concepts from abstract algebra such as group, ring, and module are used to describe the matrices/vectors under linear algebraic operators; the geometric structures such as vector space, metric, inner product, topology, manifold, etc. are essential for describing the state space of the systems. Finally the Lie group and Lie algebra structure over matrices leads to the general linear group and general linear algebra, which reveal the inside structure of the dynamics of linear (control) systems.

Recently, the STPs/STAs have been developed to overcome the dimension barrier of classical matrix theory. As the cross-dimensional operators, STPs/STAs demonstrate their dimension-free characteristic, which  makes the classical algebraic and geometric structures unsuitable for MVMDs. They, therefore, are stimulating a new matrix theory to merge. This new matrix theory demands a new mathematical foundation. This paper  explores a new foundation, which is called the CDM. The CDM consists of three parts: hyper algebra, hyper geometry, and hyper linear algebra and hyper linear group.

In this paper the hyper algebra consists of three new concepts: the hyper group, hyper ring, and hyper module.
The hyper group extends the classical group to a group with multiple identities, each identity is used to construct a component group.  The cross-dimensional operator joints all component groups together. The equivalence group shows the unified properties of component groups. Based on hyper group, the hyper ring and hyper module are also proposed in a natural way. They have similar structure: component objects and equivalence rings/modules.

The hyper geometry uses the hyper group structure to construct a hyper vector space, which has multiple zeros.  Then the hyper inner product  and hyper metric  are also obtained for cross-dimensional Euclidian space  $\R^{\infty}$. Furthermore, the metric-based topological structure and a differentiable structure have also been posed on $\R^{\infty}$. Using such topological and differentiable structure, the hyper manifold is constructed. These objects provide a cross-dimensional geometric structure for MVMDs.

Finally, the Lie algebra and Lie group structure for hyper matrix group with respect to $\hat{+}$ and $\hat{\ttimes}$ respectively, is proposed. The topological hyper group is also constructed. They eventually lead to the hyper-Lie algebra and the hyper-Lie group for MVMDs, called the hyper general linear algebra and hyper general linear group.

\subsection{Potential Applications}

Dimension-varying systems appear in nature and from engineering systems. For instance,  in a genetic regulatory network, cells may die or birth at any time; in the internet or some other service-based networks, some users may join in or withdraw out from time to time; some mechanical systems, such as docking, undocking, departure, and joining of spacecrafts \cite{yan14,jia14}, modeling of biological systems \cite{xu81,hua95}, uncertain model identification \cite{wan12}, etc.

There are also many dimension uncertain systems. For instance, a single electric power generator may be described by 2-dimensional model, or 3-dimensional, or even 5-, 6-, or 7-dimensional models \cite{mac97}. From theoretical physics, the superstring theory proposes the space-time dimension could be 4 (Instain, Relativity), 5 (Kalabi-Klein theory), 10 (Type 1 string), 11 (M-theory), or even 26 (Bosonic), etc. \cite{kak99}.

To deal with dimension varying and/or dimension uncertain systems, current fixed dimension mathematics seems not very efficient.
In recent years, we try to use STP/STA to treat these kind of systems \cite{che20,che23,che23b,chepr}. Our previous approaches are based on the dynamics over quotient space, which has nice vector space structure. We then found that this approach is not suitable for dealing with networked systems, which have higher dimensions and the dimension-varying is relatively small. This fact stimulates us to investigate the structure of real spaces with mixed dimensions itself. That is the real motivation of the current paper.

By using the real space as the state-space of time-varying systems, it may provide a practically useful frame for modeling and analyzing the time-varying or uncertain systems and for designing controls.

Since the topic of this paper is wide and essential, what we did so far is very limited and very elementary. Hence, this paper is only providing a skeleton for what CDM should be. Further exploring and investigation are expected.


\end{document}